# On an integral involving the digamma function

Donal F. Connon

dconnon@btopenworld.com

5 December 2012

**Abstract**

We consider several possible approaches to evaluating the integral $\int_0^1 \psi^2(1+x)\,dx$ and the related logarithmic series $\sum_{n=1}^{\infty} \frac{\log(n+1)}{n(n+1)}$ and, in the process, we show multifarious connections with the Hurwitz zeta function, Stieltjes constants and the generalized Euler constant function which was recently introduced by Sondow and Hadjicostas.

In Nielsen's treatise, Die Gammafunktion, there are to be found several little known formulae relating the digamma function to series involving the harmonic numbers. We show how these series play an important role in evaluating various integrals and series.

CONTENTS                                                                                   Page

## 1. Introduction

The Stieltjes constants $\gamma_n(x)$ are the coefficients of the Laurent expansion of the Hurwitz zeta function $\varsigma(s,x)$ about $s=1$

(1.1) $$\varsigma(s,x) = \sum_{n=0}^{\infty} \frac{1}{(n+x)^s} = \frac{1}{s-1} + \sum_{n=0}^{\infty} \frac{(-1)^n}{n!}\gamma_n(x)(s-1)^n$$

and we have [82]

(1.2) $\gamma_0(x) = -\psi(x)$

where $\psi(x)$ is the digamma function.

With $x = 1$ equation (1.1) reduces to the Riemann zeta function

(1.3) $\varsigma(s) = \dfrac{1}{s-1} + \sum_{n=0}^{\infty} \dfrac{(-1)^n}{n!} \gamma_n (s-1)^n$

where $\gamma_n = \gamma_n(1)$.

Since $\lim_{s \to 1}\left[\varsigma(s) - \dfrac{1}{s-1}\right] = \gamma$ we see that $\gamma_0 = \gamma$ where $\gamma$ is Euler's constant. It may be shown, as in [49, p.4], that

(1.4) $\gamma_n(x) = \lim_{N \to \infty}\left[\sum_{k=0}^{N} \dfrac{\log^n(k+x)}{k+x} - \dfrac{\log^{n+1}(N+x)}{n+1}\right]$

and

(1.4.1) $\gamma_n = \lim_{N \to \infty}\left[\sum_{k=1}^{N} \dfrac{\log^n k}{k} - \dfrac{\log^{n+1} N}{n+1}\right] = \lim_{N \to \infty}\left[\sum_{k=1}^{N} \dfrac{\log^n k}{k} - \int_1^N \dfrac{\log^n t}{t} dt\right]$

where, throughout this paper, we define $\log^0 1 = 1$.

The generalised harmonic numbers $H_n^{(r)}$ are defined by

$$H_n^{(r)} = \sum_{k=1}^{n} \dfrac{1}{k^r}$$

and, by convention, we frequently denote $H_n = H_n^{(1)}$. Since $H_{n+1}^{(r)} = H_n^{(r)} + \dfrac{1}{(n+1)^r}$ we see that it is consistent to define $H_0^{(r)} = 0$.

## 2. A logarithmic series involving the Stieltjes constants

Kanemitsu et al. [50] showed in 2004 that

(2.1) $\sum_{k=1}^{\infty} H_k \left(\log \dfrac{k+1}{k} - \dfrac{1}{k}\right) = -\dfrac{1}{2}\left[\varsigma(2) + \gamma^2 + 2\gamma_1\right]$

Their formula has been amended so as to conform to the definition of $\gamma_1$ in (1.3). A slightly expanded version of their proof is set out below.

Kanemitsu et al. [50] noted the identity

$$F(k)[H_{k+1} - H_k] = \frac{F(k)}{k+1}$$

which is valid for any function $F(x)$. We make the summation

$$\sum_{k=1}^{n} F(k)[H_{k+1} - H_k] = \sum_{k=1}^{n} \frac{F(k)}{k+1}$$

which may be written as

(2.2) $$\sum_{k=2}^{n} H_k[F(k-1) - F(k)] = \sum_{k=1}^{n} \frac{F(k)}{k+1} + F(1) - F(n)H_{n+1}$$

which is effectively Abel's summation by parts formula

$$\sum_{k=1}^{n} u_k v_k = \sum_{k=1}^{n-1} U_k(v_k - v_{k+1}) + U_n v_n$$

where $U_k = \sum_{j=1}^{n} u_k$

In particular, we consider the case $F(k) = \log(k+1)$ so that

$$\sum_{k=1}^{n} H_k \log \frac{k+1}{k} = \log(n+1)H_{n+1} - \sum_{k=1}^{n+1} \frac{\log k}{k}$$

and we have

$$\sum_{k=1}^{n} H_k \left( \log \frac{k+1}{k} - \frac{1}{k} \right) = \log(n+1)H_{n+1} - \sum_{k=1}^{n+1} \frac{\log k}{k} - \sum_{k=1}^{n} \frac{H_k}{k}$$

$$= \log(n+1)[\log(n+1) + \gamma + \varepsilon_{n+1}] - \sum_{k=1}^{n+1} \frac{\log k}{k} - \sum_{k=1}^{n} \frac{H_k}{k}$$

$$= \log^2(n+1) + \gamma \log(n+1) + \varepsilon_{n+1} \log(n+1) - \sum_{k=1}^{n+1} \frac{\log k}{k} - \sum_{k=1}^{n} \frac{H_k}{k}$$

where $\varepsilon_n = H_n - \log n - \gamma$ and it is well known that $\lim_{n \to \infty} \varepsilon_n = 0$.

This may be written as

(2.3)
$$\sum_{k=1}^{n} H_k \left( \log \frac{k+1}{k} - \frac{1}{k} \right) = -\left[ \sum_{k=1}^{n+1} \frac{\log k}{k} - \gamma \log(n+1) \right] - \left[ \sum_{k=1}^{n} \frac{H_k}{k} - \frac{1}{2} \log^2(n+1) \right] + \varepsilon_{n+1} \log(n+1)$$

and we want to consider the limit as $n \to \infty$.

Chao-Ping Chen and Feng Qi [13] showed in 2003 that

(2.4)
$$\frac{1}{2n + \beta} \leq H_n - \log n - \gamma < \frac{1}{2n + \alpha}$$

where $\alpha = \frac{1}{1-\gamma} - 2$, $\beta = \frac{1}{3}$, and the constants $\alpha$ and $\beta$ are the best possible. Hence, using L'Hôpital's rule, we easily see that

(2.5)
$$\lim_{n \to \infty} [\varepsilon_{n+1} \log(n+1)] = 0.$$

We recall Adamchik's formula [1]

(2.6)
$$\sum_{k=1}^{n} \frac{H_k}{k} = \frac{1}{2} \left( H_n^{(1)} \right)^2 + \frac{1}{2} H_n^{(2)}$$

which was also reported by M.E. Levenson in a 1938 volume of The American Mathematical Monthly [59] in a problem concerning the evaluation of

(2.7)
$$\Gamma''(1) = \int_0^\infty e^{-x} \log^2 x \, dx = \gamma^2 + \varsigma(2)$$

Substituting $\left( H_n^{(1)} \right)^2 = \left( H_n^{(1)} - \log n - \gamma + \log n + \gamma \right)^2$ in (2.6) we have

$$\sum_{k=1}^{n} \frac{H_k}{k} - \frac{1}{2} \left( H_n^{(1)} - \log n - \gamma \right)^2 - \left( H_n^{(1)} - \log n - \gamma \right)(\log n + \gamma) - \frac{1}{2}(\log n + \gamma)^2 = \frac{1}{2} H_n^{(2)}$$

Therefore using (2.5) we obtain the limit

$$\lim_{n \to \infty} \left[ \sum_{k=1}^{n} \frac{H_k}{k} - \frac{1}{2}(\log n + \gamma)^2 \right] = \frac{1}{2} \varsigma(2)$$

or equivalently

$$(2.8) \quad \lim_{n \to \infty} \left[ \sum_{k=1}^{n} \frac{H_k}{k} - \gamma \log n - \frac{1}{2} \log^2 n \right] = \frac{1}{2} \left[ \varsigma(2) + \gamma^2 \right] = \frac{1}{2} \Gamma''(1)$$

where the $\Gamma''(1)$ connection comes from (2.7).

I originally derived this formula in [26] and subsequently noted that it had been previously discovered by Kanemitsu et al. [50] in 2004.

We note from (1.2) that

$$(2.9) \quad \gamma_1 = \lim_{n \to \infty} \left[ \sum_{k=1}^{n} \frac{\log k}{k} - \frac{1}{2} \log^2 n \right]$$

and using (2.5), (2.8) and (2.9) we then see from (2.3) that

$$(2.10) \quad \sum_{k=1}^{\infty} H_k \left( \log \frac{k+1}{k} - \frac{1}{k} \right) = -\frac{1}{2} \left[ \varsigma(2) + \gamma^2 + 2\gamma_1 \right]$$

□

We may also write (2.8) as

$$(2.11) \quad \lim_{n \to \infty} \left[ \frac{1}{2} \left( H_n^{(1)} \right)^2 - \gamma \log n - \frac{1}{2} \log^2 n \right] = \frac{1}{2} \gamma^2$$

It may be noted that equation (2.8) concurs with the asymptotic formula obtained by Flajolet and Sedgewick [40]

$$(2.12) \quad \sum_{k=1}^{n} \frac{H_k}{k} = \frac{1}{2} \log^2 n + \gamma \log n + \frac{1}{2} \left[ \varsigma(2) + \gamma^2 \right] + O\left( \frac{\log n}{n} \right)$$

From (2.11) we see that

$$\lim_{n \to \infty} \left[ \frac{1}{2} \left( H_n^{(1)} \right)^2 - \gamma \log n - \sum_{k=1}^{n} \frac{\log k}{k} + \sum_{k=1}^{n} \frac{\log k}{k} - \frac{1}{2} \log^2 n \right] = \frac{1}{2} \gamma^2$$

and referring to (2.9) we get

$$(2.13) \quad \lim_{n \to \infty} \left[ \frac{1}{2} \left( H_n^{(1)} \right)^2 - \gamma \log n - \sum_{k=1}^{n} \frac{\log k}{k} \right] = \frac{1}{2} \gamma^2 - \gamma_1$$

Equation (2.8) may be written as

$$\lim_{n\to\infty}\left[\sum_{k=1}^{n}\frac{H_k}{k}-\sum_{k=1}^{n}\frac{\log k}{k}-\gamma\log n+\sum_{k=1}^{n}\frac{\log k}{k}-\frac{1}{2}\log^2 n\right]=\frac{1}{2}\left[\varsigma(2)+\gamma^2\right]$$

and using the definition (1.4) of Stieltjes constant $\gamma_1$ we have

$$\lim_{n\to\infty}\left[\sum_{k=1}^{n}\frac{H_k}{k}-\sum_{k=1}^{n}\frac{\log k}{k}-\gamma\log n\right]=\frac{1}{2}\left[\varsigma(2)+\gamma^2\right]-\gamma_1$$

or equivalently

$$\lim_{n\to\infty}\left[\sum_{k=1}^{n}\frac{H_k-\log k}{k}-\gamma\log n\right]=\frac{1}{2}\left[\varsigma(2)+\gamma^2\right]-\gamma_1$$

This may also be expressed as

$$\lim_{n\to\infty}\left[\sum_{k=1}^{n}\frac{H_k-\gamma-\log k}{k}+\sum_{k=1}^{n}\frac{\gamma}{k}-\gamma\log n\right]=\frac{1}{2}\left[\varsigma(2)+\gamma^2\right]-\gamma_1$$

or

$$\lim_{n\to\infty}\left[\sum_{k=1}^{n}\frac{H_k-\gamma-\log k}{k}+\gamma(H_n-\log n)\right]=\frac{1}{2}\left[\varsigma(2)+\gamma^2\right]-\gamma_1$$

and since $\lim_{n\to\infty}\gamma(H_n-\log n)=\gamma^2$ we therefore conclude that

(2.14) $$\sum_{k=1}^{\infty}\frac{H_k-\gamma-\log k}{k}=\frac{1}{2}\left[\varsigma(2)-\gamma^2\right]-\gamma_1$$

This formula was proposed as a problem by Furdui [41] in 2007.

We have from [70]

(2.15) $$H_k-\log k-\gamma=\int_0^{\infty}e^{-kx}\left[\frac{1}{x}-\frac{1}{e^x-1}\right]dx$$

and on summation we obtain

$$\sum_{k=1}^{\infty}\frac{H_k-\gamma-\log k}{k}=\sum_{k=1}^{\infty}\int_0^{\infty}\frac{e^{-kx}}{k}\left[\frac{1}{x}-\frac{1}{e^x-1}\right]dx$$

$$=\int_0^{\infty}\sum_{k=1}^{\infty}\frac{e^{-kx}}{k}\left[\frac{1}{x}-\frac{1}{e^x-1}\right]dx$$

$$= \int_0^\infty \log(1-e^{-x})\left[\frac{1}{e^x-1}-\frac{1}{x}\right]dx$$

Hence using (2.14) we have

(2.16) $$\int_0^\infty \log(1-e^{-x})\left[\frac{1}{e^x-1}-\frac{1}{x}\right]dx = \frac{1}{2}\left[\varsigma(2)-\gamma^2\right]-\gamma_1$$

which was previously obtained by Coffey [19]. With the substitution $u = 1-e^{-x}$ we have

(2.16.1) $$\int_0^1 \left[\frac{1}{u}+\frac{1}{(1-u)\log(1-u)}\right]\log u\, du = \frac{1}{2}\left[\varsigma(2)-\gamma^2\right]-\gamma_1$$

With $v = 1-u$ this becomes

(2.16.2) $$\int_0^1 \left[\frac{1}{1-v}+\frac{1}{v\log v}\right]\log(1-v)\, dv = \frac{1}{2}\left[\varsigma(2)-\gamma^2\right]-\gamma_1$$

We shall make use of this integral in (3.60) and (10.18). This may be compared with

(2.16.3) $$\int_0^\infty e^{-x}\left[\frac{1}{1-e^{-x}}-\frac{1}{x}\right]\log x\, dx = -\gamma^2 - \gamma_1$$

as previously determined by Coppo [32] in 1999. With the substitution $u = e^{-x}$ we have

(2.16.4) $$\int_0^1 \left[\frac{1}{1-u}+\frac{1}{\log u}\right]\log(-\log u)\, du = -\gamma^2 - \gamma_1$$

We also showed in a similar way in [26] that

(2.17) $$\int_0^\infty Li_2(e^{-x})\left[\frac{1}{x}-\frac{1}{e^x-1}\right]dx = 2\varsigma(3)+\varsigma'(2)-\gamma\varsigma(2)$$

where $Li_s(x)$ is the polylogarithm function.

## 3. Nielsen's series for the digamma function

Nielsen [67, p.52] showed that

(3.1) $$[\psi(x)+\gamma]^2 = \psi'(x)-\varsigma(2)-2\xi(x)$$

where

(3.2) $$\xi(x) = \sum_{n=1}^{\infty} H_n \left( \frac{1}{x+n} - \frac{1}{n+1} \right)$$

Integration results in

$$\int_0^u \xi(x)\,dx = \sum_{n=1}^{\infty} H_n \left( \log\frac{n+u}{n} - \frac{u}{n+1} \right)$$

$$= \sum_{n=1}^{\infty} H_n \left( \log\frac{n+u}{n} - \frac{u}{n} + \frac{u}{n} - \frac{u}{n+1} \right)$$

$$= \sum_{n=1}^{\infty} H_n \left( \log\frac{n+u}{n} - \frac{u}{n} \right) + u\sum_{n=1}^{\infty} H_n \left( \frac{1}{n} - \frac{1}{n+1} \right)$$

and we easily see that

$$\sum_{n=1}^{\infty} H_n \left( \frac{1}{n} - \frac{1}{n+1} \right) = \sum_{n=1}^{\infty} \left( \frac{H_n}{n} - \frac{H_{n+1}}{n+1} + \frac{1}{(n+1)^2} \right) = \varsigma(2)$$

Hence we obtain

(3.3) $$\int_0^u \xi(x)\,dx = \sum_{n=1}^{\infty} H_n \left( \log\frac{n+u}{n} - \frac{u}{n} \right) + \varsigma(2)u$$

and with $u=1$ we obtain

(3.4) $$\int_0^1 \xi(x)\,dx = \sum_{n=1}^{\infty} H_n \left( \log\frac{n+1}{n} - \frac{1}{n} \right) + \varsigma(2)$$

and using (2.1) this becomes

(3.5) $$2\int_0^1 \xi(x)\,dx = \varsigma(2) - \gamma^2 - 2\gamma_1$$

From (3.1) we have

$$2\int_0^1 \xi(x)\,dx = \int_0^1 \left( \psi'(x) - [\psi(x)+\gamma]^2 \right) dx - \varsigma(2)$$

so that

$$\int_0^1 \left( \psi'(x) - [\psi(x) + \gamma]^2 \right) dx - \varsigma(2) = \varsigma(2) - \gamma^2 - 2\gamma_1$$

and this gives us an apparently new integral representation for $\gamma_1$

(3.6) $$\int_0^1 \left[ \psi'(x) - \psi^2(x) - 2\gamma \psi(x) \right] dx = 2\varsigma(2) - 2\gamma_1$$

which we reported in 2009 in [28].

However, I subsequently discovered that the integral is not new when I noted that the following mild variant of it appears in Cohen's book [21, p.145] which was published in 2007

(3.7) $$\int_0^1 \left[ \psi^2(x) - \frac{1}{x^2} - \frac{2\gamma}{x} \right] dx = 1 - 2\varsigma(2) + 2\gamma_1$$

and the equivalence may be easily seen by writing

$$\int_0^1 \left[ \psi'(x) - \psi^2(x) - 2\gamma \psi(x) \right] dx = \int_0^1 \left[ \psi'(1+x) + \frac{1}{x^2} - \psi^2(x) - 2\gamma \psi(1+x) + \frac{2\gamma}{x} \right] dx$$

An equivalent representation was also given by Coffey [20] in 2011.

Using the functional equation for the digamma function we may express (3.1) as

$$\psi^2(1+x) - \frac{2}{x}\psi(1+x) + \frac{1}{x^2} + 2\gamma \psi(1+x) - \frac{2\gamma}{x} + \gamma^2 = \psi'(1+x) + \frac{1}{x^2} - \varsigma(2) - 2\xi(x)$$

so that

$$\psi^2(1+x) = 2\frac{\psi(1+x) + \gamma}{x} + \psi'(1+x) - 2\gamma \psi(1+x) - \varsigma(2) - \gamma^2 - 2\xi(x)$$

Integration then results in

$$\int_0^1 \psi^2(1+x) dx = 2\int_0^1 \frac{\psi(1+x) + \gamma}{x} dx + \int_0^1 \psi'(1+x) dx - 2\gamma \int_0^1 \psi(1+x) dx - \varsigma(2) - \gamma^2 - 2\int_0^1 \xi(x) dx$$

$$= 2\int_0^1 \frac{\psi(1+x) + \gamma}{x} dx + 1 - \varsigma(2) - \gamma^2 - 2\int_0^1 \xi(x) dx$$

Cohen [21, p.142] has reported as an exercise that

(3.8) $$\int_0^1 \frac{\psi(1+x)+\gamma}{x}\,dx = \sum_{n=1}^{\infty} \frac{\log(n+1)}{n(n+1)}$$

and employing (3.5) we obtain

(3.9) $$\int_0^1 \psi^2(1+x)\,dx = 2\sum_{n=1}^{\infty} \frac{\log(n+1)}{n(n+1)} + 1 - 2\varsigma(2) + 2\gamma_1$$

It is easily seen that this is consistent with (3.7) because substituting

$$\psi^2(x) = \psi^2(1+x) - \frac{2}{x}\psi(x) + \frac{1}{x^2}$$

we obtain

(3.9.1) $$\int_0^1 \left[\psi^2(1+x) - 2\frac{\psi(1+x)+\gamma}{x}\right]dx = 1 - 2\varsigma(2) + 2\gamma_1$$

Cohen [21, p.142] has also stated that

(3.10.1) $$\sum_{n=1}^{\infty} \frac{\log(n+1)}{n(n+1)} = \int_0^1 \frac{(1-x)\log(1-x)}{x\log x}\,dx$$

(3.10.2) $$= \sum_{n=1}^{\infty} (-1)^{n+1} \frac{\varsigma(n+1)}{n}$$

(3.10.3) $$= -\sum_{n=2}^{\infty} \varsigma'(n)$$

(3.10.4) $$= \sum_{n=1}^{\infty} \frac{1}{n}\log\frac{n+1}{n}$$

In addition, Coffey [17] has reported that

(3.10.5) $$\sum_{n=1}^{\infty} \frac{\log(n+1)}{n(n+1)} = \int_0^{\infty} \frac{(e^{-x}-1)\log(1-e^{-x})}{x}\,dx$$

(3.10.6) $$= -\gamma + \int_0^{\infty} e^{-x}\log x \log(1-e^{-x})\,dx$$

(3.10.7) $$= \frac{3}{2} - 2\gamma + \frac{\pi^2}{12} + \int_0^\infty \left[\pi \coth \pi x - \frac{1}{x}\right] \frac{dx}{e^{2\pi x} - 1}$$

It is easily seen that the substitution $x = e^{-t}$ in (3.10.1) results in (3.10.5). Using integration by parts we may deduce (3.10.6) from (3.10.5). Coffey [17] used Binet's integral representation of the digamma function to derive (3.10.7).

Some derivations of the above formulae are shown below where we have attempted to supply a number of disparate proofs.

□

With regard to (3.10.6), the substitution $u = e^{-x}$ gives us

$$\int_0^\infty e^{-x} \log x \log(1 - e^{-x}) \, dx = \int_0^1 \log(-\log u) \log(1 - u) \, du$$

and we see that

$$\int_0^1 \log(-\log u) \log(1 - u) \, du = -\sum_{n=1}^\infty \int_0^1 \frac{u^n}{n} \log(-\log u) \, du$$

Medina and Moll [63] have shown that

(3.11) $$\int_0^1 u^a \log(-\log u) \, du = -\frac{\gamma + \log(1 + a)}{1 + a}$$

and hence we have

$$\int_0^1 \log(-\log u) \log(1 - u) \, du = \sum_{n=1}^\infty \frac{\gamma + \log(n + 1)}{n(n + 1)}$$

so that

(3.12) $$\int_0^1 \log(-\log u) \log(1 - u) \, du = \gamma + \sum_{n=1}^\infty \frac{\log(n + 1)}{n(n + 1)}$$

which shows the equivalence between (3.8) and (3.10.6).

□

We have the well-known integral

$$\log t = \int_0^1 \frac{y^{t-1}-1}{\log y} dy$$

which may be easily derived by integrating the identity $\frac{1}{x} = \int_0^1 t^{x-1} dt$. Hence we see that

$$\log \frac{n+1}{n} = \int_0^1 \frac{y^n - y^{n-1}}{\log y} dy = \int_0^1 \frac{y^{n-1}(y-1)}{\log y} dy$$

and we have the summation

$$\sum_{n=1}^{\infty} \frac{1}{n} \log \frac{n+1}{n} = \sum_{n=1}^{\infty} \int_0^1 \frac{y^n}{n} \frac{(y-1)}{y \log y} dy$$

$$= \int_0^1 \sum_{n=1}^{\infty} \frac{y^n}{n} \frac{(y-1)}{y \log y} dy$$

$$= \int_0^1 \frac{(1-y)\log(1-y)}{y \log y} dy$$

This demonstrates the equivalence of (3.10.1) and (3.10.4).

□

We note the representation [79, p.14] for the digamma function

(3.13) $$\psi(1+x) + \gamma = x \sum_{n=1}^{\infty} \frac{1}{n(n+x)}$$

whereupon we easily see that

(3.14) $$\int_0^u \frac{\psi(1+x)+\gamma}{x} dx = \sum_{n=1}^{\infty} \frac{1}{n} \log \frac{n+u}{n}$$

With $u = 1$ this shows the equivalence of (3.8) and (3.10.4).

□

We have the familiar integral for the digamma function [79, p.15]

(3.15) $$\psi(1+x) + \gamma = \int_0^1 \frac{1-y^x}{1-y} dy$$

where integration by parts gives us

$$= (1-y^x)\log(1-y)\Big|_0^1 - x\int_0^1 y^{x-1}\log(1-y)\,dy$$

and thus we see that [67, p.172]

(3.16) $$\psi(1+x)+\gamma = -x\int_0^1 y^{x-1}\log(1-y)\,dy$$

We then have

$$\frac{\psi(1+x)+\gamma}{x} = -\int_0^1 y^{x-1}\log(1-y)\,dy$$

and integrating this gives us

$$\int_0^1 \frac{\psi(1+x)+\gamma}{x}\,dx = \int_0^1 \frac{(1-y)\log(1-y)}{y\log y}\,dy$$

This demonstrates the equivalence of (3.8) and (3.10.1).

$\square$

We may prove (3.10.1) by reference to the following three integrals

$$\log(n+1) = \int_0^1 \frac{t^n-1}{\log t}\,dt$$

$$\frac{1}{n} = \int_0^1 u^{n-1}\,du$$

$$\frac{1}{n+1} = \int_0^1 v^n\,dv$$

We multiply the above three factors together and make the summation to form the triple integral

$$\sum_{n=1}^{\infty} \frac{\log(n+1)}{n(n+1)} = \sum_{n=1}^{\infty} \int_0^1 \frac{t^n-1}{\log t}\,dt \int_0^1 u^{n-1}\,du \int_0^1 v^n\,dv$$

$$= \int_0^1\int_0^1\int_0^1 \left[\frac{tuv}{1-tuv} - \frac{uv}{1-uv}\right]\frac{1}{u\log t}\,dt\,du\,dv$$

$$= \int_0^1\int_0^1\int_0^1 \left[\frac{tv}{1-tuv} - \frac{v}{1-uv}\right]\frac{1}{\log t}\,dt\,du\,dv$$

$$= \int_0^1\int_0^1 \frac{\log(1-v) - \log(1-tv)}{\log t}\,dt\,dv$$

$$= \int_0^1 \frac{(1-t)\log(1-t)\,dt}{t\log t}$$

which is (3.10.1).

$\square$

Euler's gamma function is defined for $\mathrm{Re}(s) > 0$ as

$$\Gamma(s) = \int_0^\infty x^{s-1} e^{-x}\,dx$$

and we see that

$$\frac{1}{a^s} = \frac{1}{\Gamma(s)}\int_0^\infty x^{s-1} e^{-ax}\,dx$$

Differentiation with respect to $s$ results in

$$-\frac{\log a}{a^s} = \frac{1}{\Gamma(s)}\int_0^\infty x^{s-1} e^{-ax} \log x\,dx - \psi(s)\int_0^\infty x^{s-1} e^{-ax}\,dx$$

and with $a = n+1$ and $s = 1$ we have

$$\frac{\log(n+1)}{n+1} = -\int_0^\infty e^{-(n+1)x} \log x\,dx - \frac{\gamma}{n+1}$$

We have the summation

$$\sum_{n=1}^\infty \frac{\log(n+1)}{n(n+1)} = -\int_0^\infty \sum_{n=1}^\infty \frac{e^{-nx}}{n} e^{-x} \log x\,dx - \gamma\sum_{n=1}^\infty \frac{1}{n(n+1)}$$

$$= \int_0^\infty \log(1-e^{-x})e^{-x} \log x \, dx - \gamma$$

and this constitutes another derivation of Coffey's representation (3.10.6).

□

We see that

$$\sum_{n=1}^\infty \frac{\log(n+1)}{n(n+1)} = \sum_{n=1}^\infty \left[ \frac{\log(n+1)}{n} - \frac{\log(n+1)}{n+1} \right]$$

$$= \sum_{n=1}^\infty \left[ \frac{\log(n+1)}{n} - \frac{\log n}{n} - \left\{ \frac{\log(n+1)}{n+1} - \frac{\log n}{n} \right\} \right]$$

$$= \sum_{n=1}^\infty \frac{1}{n} \log \frac{n+1}{n}$$

and thus we have an elementary demonstration of the equivalence of (3.10.1) and (3.10.4).

□

It may be noted that

$$\sum_{n=1}^\infty \frac{\log(n+1)}{n(n+1)} = -\frac{d}{ds} \sum_{n=1}^\infty \frac{1}{n(n+1)^s} \bigg|_{s=1}$$

and we see that

$$\sum_{n=1}^\infty \frac{1}{n(n+1)^s} = \int_0^1 \frac{1}{z} \sum_{n=1}^\infty \frac{z^n}{(n+1)^s} \, dz$$

$$= \int_0^1 \frac{1}{z} \left[ \sum_{n=0}^\infty \frac{z^n}{(n+1)^s} - 1 \right] dz$$

$$= \int_0^1 \frac{1}{z} [\Phi(z,s,1) - 1] \, dz$$

where $\Phi(z,s,a)$ is the Hurwitz-Lerch zeta function defined by

$$\Phi(z,s,a) = \sum_{n=0}^\infty \frac{z^n}{(n+a)^s}$$

It is well known that [79, p.121]

$$\Phi(z,s,a) = \frac{1}{\Gamma(s)} \int_0^\infty \frac{t^{s-1} e^{-(a-1)t}}{e^t - z} dt$$

and we have the well known integral due to Appell [80, p.280]

$$z\,\Phi(z,s,1) = Li_s(z) = \frac{z}{\Gamma(s)} \int_0^\infty \frac{t^{s-1}}{e^t - z} dt$$

Therefore we have

$$\sum_{n=1}^\infty \frac{1}{n(n+1)^s} = \int_0^1 \frac{1}{z} \left[ \frac{1}{\Gamma(s)} \int_0^\infty \frac{t^{s-1}}{e^t - z} dt - 1 \right] dz$$

and using the integral definition of the gamma function $\Gamma(s) = \int_0^\infty x^{s-1} e^{-x} dx$ we obtain

$$S = \sum_{n=1}^\infty \frac{1}{n(n+1)^s} = \frac{1}{\Gamma(s)} \int_0^1 \frac{1}{z} \left[ \int_0^\infty \frac{t^{s-1}}{e^t - z} dt - \int_0^\infty t^{s-1} e^{-t} dt \right] dz$$

$$= \frac{1}{\Gamma(s)} \int_0^1 \int_0^\infty \left[ \frac{1}{z(e^t - z)} - \frac{e^{-t}}{z} \right] t^{s-1} dz\, dt$$

Since $\dfrac{1}{z(e^t - z)} = \dfrac{e^{-t}}{z} + \dfrac{e^{-t}}{e^t - z}$ we get

$$S = \frac{1}{\Gamma(s)} \int_0^1 \int_0^\infty \frac{t^{s-1} e^{-t}}{e^t - z} dz\, dt$$

We have the elementary integral

$$\int_0^1 \frac{dz}{e^t - z} = t - \log(e^t - 1)$$

$$= -\log(1 - e^{-t})$$

and thus

$$\sum_{n=1}^{\infty}\frac{1}{n(n+1)^s} = -\frac{1}{\Gamma(s)}\int_0^{\infty} t^{s-1}e^{-t}\log(1-e^{-t})\,dt$$

Differentiation then results in

(3.17) $$\sum_{n=1}^{\infty}\frac{\log(n+1)}{n(n+1)^s} = \frac{1}{\Gamma(s)}\int_0^{\infty} t^{s-1}e^{-t}\log t\log(1-e^{-t})\,dt - \frac{\psi(s)}{\Gamma(s)}\int_0^{\infty} t^{s-1}e^{-t}\log(1-e^{-t})\,dt$$

and with $s=1$ we have

$$\sum_{n=1}^{\infty}\frac{\log(n+1)}{n(n+1)} = \int_0^{\infty} e^{-t}\log t\log(1-e^{-t})\,dt + \gamma\int_0^{\infty} e^{-t}\log(1-e^{-t})\,dt$$

With the substitution $u = e^{-t}$ we see that $\int_0^{\infty} e^{-t}\log(1-e^{-t})\,dt = -1$ and thus we have obtained (3.10.6) again

$$\sum_{n=1}^{\infty}\frac{\log(n+1)}{n(n+1)} = \int_0^{\infty} e^{-t}\log t\log(1-e^{-t})\,dt - \gamma$$

□

We recall (3.15)

$$\psi(1+x)+\gamma = \int_0^1 \frac{u^x - 1}{u - 1}\,du$$

and we have the double integral

(3.17.1) $$\int_0^u \frac{\psi(1+x)+\gamma}{x}\,dx = \int_0^1\int_0^1 \frac{u^x - 1}{x(u-1)}\,du\,dx$$

We have the exponential integral

$$Ei(x) = \int \frac{e^x}{x}\,dx$$

so that

$$\int \frac{u^x - 1}{x}\,dx = Ei(x\log u) - \log x$$

We see that

$$Ei(x\log u) - \log x = Ei(x\log u) - \log(-x\log u) + \log(-x\log u) - \log x$$

$$= Ei(x\log u) - \log(-x\log u) + \log(-\log u)$$

It is known [44, p.877] that for $z > 0$

(3.17.2) $$Ei(z) = \gamma + \log z + \sum_{n=1}^{\infty} \frac{z^n}{n.n!}$$

and for $z < 0$ we have

(3.17.3) $$Ei(z) = \gamma + \log(-z) + \sum_{n=1}^{\infty} \frac{z^n}{n.n!}$$

If $u \in (0,1)$ and $x \in (0,1)$, then $x \log u < 0$ and hence we have

(3.17.4) $$\lim_{x \to 0+}[Ei(x\log u) - \log x] = \gamma + \log(-\log u)$$

Therefore we obtain the definite integral

(3.17.5) $$\int_0^1 \frac{u^x - 1}{x} dx = Ei(\log u) - \gamma - \log(-\log u)$$

This is equivalent to

(3.17.6) $$\int_0^1 \frac{u^x - 1}{x} dx = \sum_{n=1}^{\infty} \frac{\log^n u}{n.n!}$$

and (3.17.1) gives us

$$\int_0^1 \frac{\psi(1+x) + \gamma}{x} dx = \sum_{n=1}^{\infty} \frac{1}{n.n!} \int_0^1 \frac{\log^n u}{u-1} du$$

It is well known that [10, p.240]

(3.18) $$\int_0^1 \frac{\log^n u}{1-u} du = (-1)^n n! \varsigma(n+1)$$

(which may be obtained by differentiating (3.15) $n$ times) and accordingly we have a derivation of the equivalence of (3.8) and (3.10.2)

$$\int_0^1 \frac{\psi(1+x)+\gamma}{x} dx = \sum_{n=1}^{\infty} \frac{(-1)^{n+1}\varsigma(n+1)}{n}$$

It may be noted that differentiating (3.17.6) gives us

$$\int_0^1 u^{x-1} du = \frac{u^{x-1}}{\log u}\bigg|_0^1 = \frac{u-1}{u\log u}$$

and for the right-hand side we also have the derivative

$$\frac{1}{u\log u}\sum_{n=1}^{\infty} \frac{\log^n u}{n!} = \frac{u-1}{u\log u}$$

$\square$

Alternatively, we may expand the integrand as follows

$$\psi(1+x)+\gamma = \int_0^1 \frac{u^x-1}{u-1} du$$

$$= \int_0^1 \frac{\exp(x\log u)-1}{u-1} du$$

$$= \sum_{n=1}^{\infty} \frac{1}{n!}\int_0^1 \frac{x^n \log^n u}{1-u} du$$

and using (3.18) we obtain

(3.18.1) $\quad \psi(1+x)+\gamma = \sum_{n=1}^{\infty}(-1)^n \varsigma(n+1)x^n$

which is the familiar Maclaurin series for the digamma function.

$\square$

As shown by Nielsen [66, p.3], Bromwich [11, p.334] and Havil [48, p.106], we have for $0 < x < 1$

(3.19.1) $\quad li(x) = \gamma + \log(-\log x) + \sum_{n=1}^{\infty} \frac{\log^n x}{n.n!}$

and for $x > 1$ we have

(3.19.2) $$li(x) = \gamma + \log(\log x) + \sum_{n=1}^{\infty} \frac{\log^n x}{n.n!}$$

where $li(x)$ is the logarithm integral function defined for $x < 1$ by $li(x) = \int_0^x \frac{dt}{\log t}$.

The logarithm integral was introduced in 1809 by the self-taught Bavarian mathematician, Johann Soldner (1776-1833). Mascheroni calculated Euler's constant to 32 decimal places in 1790 but only the first 19 places were correct; the remaining places were corrected by Soldner in 1809 [48, p.89]. Soldner is also known for predicting, as early as 1801, the deflection of light by gravity (it was not until 1915 that his calculations were subsequently modified by Einstein's general theory of relativity).

With $z = \log x$ in (3.17.3) we see that $z < 0$ for $x \in (0,1)$, so that we obtain

$$Ei(\log x) = \gamma + \log(-\log x) + \sum_{n=1}^{\infty} \frac{\log^n x}{n.n!}$$

and, comparing this with (3.19.1), we obtain

(3.19.3) $\quad Ei(\log x) = li(x)$

Using integration by parts we find that

$$\int \log(-\log x)\, dx = x\log(-\log x) - \int \frac{dx}{\log x}$$

$$= x\log(-\log x) - li(x)$$

$$= (x-1)\log(-\log x) + \log(-\log x) - li(x)$$

We have

$$\lim_{x \to 0}[x\log(-\log x)] = \lim_{x \to 0} \frac{\log(-\log x)}{1/x}$$

and using L'Hôpital's rule this becomes

$$= -\lim_{x \to 0} \frac{x}{\log x} = 0$$

Hence we have

(3.19.4) $\quad \lim_{x \to 0}[x\log(-\log x)] = 0$

Therefore, since $li(0) = 0$, we have the definite integral

(3.19.5) $$\int_0^u \log(-\log x)\,dx = u\log(-\log u) - li(u)$$

$$= (u-1)\log(-\log u) + \log(-\log u) - li(u)$$

By a change of variable, we see that [79, p.2]

(3.19.6) $$\Gamma(s) = \int_0^1 \log^{s-1}(1/x)\,dx$$

whereupon differentiation results in

(3.19.7) $$\Gamma'(s) = \int_0^1 \log^{s-1}(1/x)\log(-\log x)\,dx$$

so that $s=1$ gives us

(3.19.8) $$\gamma = -\Gamma'(1) = -\int_0^1 \log(-\log x)\,dx$$

We have

$$\lim_{x\to 1}[(x-1)\log(-\log x)] = \lim_{x\to 1}\frac{\log(-\log x)}{1/(x-1)}$$

and using L'Hôpital's rule this becomes

$$= -\lim_{x\to 1}\frac{(x-1)^2}{x\log x}$$

and a further application of L'Hôpital's rule shows that this limit is zero

(3.19.8.1) $$\lim_{x\to 1}[(x-1)\log(-\log x)] = 0$$

Hence, letting $u=1$ in (3.19.5), we may deduce from (3.19.8) that

(3.19.9) $$\lim_{x\to 1}[li(u) - \log(-\log u)] = \gamma$$

which concurs with (3.19.10).

With $u=1$ in (3.17.5) we see that

(3.19.10) $$\lim_{u\to 1}[Ei(\log u) - \log(-\log u)] = \gamma$$

□

We recall (3.10.2) and (3.10.4)

$$\sum_{n=1}^{\infty} \frac{1}{n}\log\left(1+\frac{1}{n}\right) = -\sum_{n=1}^{\infty} \frac{(-1)^n}{n}\varsigma(n+1)$$

and substituting (3.18) gives us

$$\sum_{n=1}^{\infty} \frac{(-1)^n}{n}\varsigma(n+1) = \int_0^1 \sum_{n=1}^{\infty} \frac{\log^n(1-x)}{n.n!}\frac{dx}{x}$$

We may also write (3.19.1) as

$$li(1-x) - \gamma - \log(-\log(1-x)) = \sum_{n=1}^{\infty} \frac{\log^n(1-x)}{n.n!}$$

which gives us

(3.19.11) $$\sum_{n=1}^{\infty} \frac{(-1)^n}{n}\varsigma(n+1) = \int_0^1 \frac{li(1-x) - \gamma - \log(-\log(1-x))}{x}dx$$

and integration by parts results in

$$= [li(1-x) - \gamma - \log(-\log(1-x))]\log x \Big|_0^1 - \int_0^1 \log x \left[-\frac{1}{\log(1-x)} + \frac{1}{(1-x)\log(1-x)}\right]dx$$

$$= [li(1-x) - \gamma - \log(-\log(1-x))]\log x \Big|_0^1 - \int_0^1 \frac{x\log x}{(1-x)\log(1-x)}dx$$

We see that

$$\lim_{x\to 0}[li(1-x) - \gamma - \log(-\log(1-x))]\log x = \lim_{x\to 0}\frac{[li(1-x) - \gamma - \log(-\log(1-x))]}{x} x\log x$$

and, since $\lim_{x\to 0}[li(1-x) - \gamma - \log(-\log(1-x))] = 0$, we may apply L'Hôpital's rule to obtain

$$\lim_{x\to 0}\frac{[li(1-x) - \gamma - \log(-\log(1-x))]}{x} = \lim_{x\to 0}\left[-\frac{1}{\log(1-x)} + \frac{1}{(1-x)\log(1-x)}\right] = 0$$

Hence we obtain

$$\sum_{n=1}^{\infty} \frac{(-1)^n}{n} \varsigma(n+1) = -\int_0^1 \frac{x \log x}{(1-x)\log(1-x)} dx$$

which is a further demonstration of the equivalence between (3.10.1) and (3.10.2).

□

It is well known that

$$\gamma = -\int_0^{\infty} e^{-x} \log x \, dx$$

and integration by parts gives us

$$\int_u^{\infty} e^{-x} \log x \, dx = e^{-u} \log u + \int_u^{\infty} \frac{e^{-x}}{x} dx$$

We see that

$$-\gamma = \lim_{u \to 0} \left[ e^{-u} \log u + \int_u^{\infty} \frac{e^{-x}}{x} dx \right]$$

$$= \lim_{u \to 0} \left[ (e^{-u} - 1) \log u + \log u + \int_u^{\infty} \frac{e^{-x}}{x} dx \right]$$

and, since $\lim_{u \to 0}(e^{-u} - 1)\log u = \lim_{u \to 0} e^{-u}(1-e^u)\log u = \lim_{u \to 0}(1-e^u)\log u$ and $\lim_{u \to 0} u^n \log u = 0$, we obtain

$$\lim_{u \to 0}(e^{-u} - 1)\log u = 0$$

and thus we have

$$-\gamma = \lim_{u \to 0} \left[ \log u + \int_u^{\infty} \frac{e^{-x}}{x} dx \right]$$

or equivalently

(3.19.12) $\quad \gamma = \lim_{u \to 0} \left[ Ei(-u) - \log u \right]$

where the exponential integral function $Ei(-u)$ is defined by

$$-Ei(-u) = \int_u^\infty \frac{e^{-x}}{x} dx$$

We may also obtain (3.19.12) by letting $x = -\log u$ in (3.19.10).

□

We recall (2.15)

$$H_k - \log k - \gamma = \int_0^\infty e^{-kx}\left[\frac{1}{x} - \frac{1}{e^x - 1}\right] dx$$

and make the summation to obtain

$$\sum_{k=1}^\infty \frac{H_k - \log k - \gamma}{k!} u^k = \int_0^\infty \sum_{k=1}^\infty \frac{(ue^{-x})^k}{k!}\left[\frac{1}{x} - \frac{1}{e^x - 1}\right] dx$$

or equivalently

(3.19.13) $$\sum_{k=1}^\infty \frac{H_k}{k!} u^k - \sum_{k=1}^\infty \frac{\log k}{k!} u^k - \gamma(e^u - 1) = \int_0^\infty [\exp(ue^{-x}) - 1]\left[\frac{1}{x} - \frac{1}{e^x - 1}\right] dx$$

With $u = 1$ this becomes

(3.19.14) $$\sum_{k=1}^\infty \frac{H_k}{k!} - \sum_{k=1}^\infty \frac{\log k}{k!} - \gamma(e - 1) = \int_0^\infty [\exp(e^{-x}) - 1]\left[\frac{1}{x} - \frac{1}{e^x - 1}\right] dx$$

and the substitution $t = e^{-x}$ gives us

$$\int_0^\infty [\exp(e^{-x}) - 1]\left[\frac{1}{x} - \frac{1}{e^x - 1}\right] dx = \int_0^1 (e^t - 1)\left[\frac{1}{t - 1} - \frac{1}{t \log t}\right] dt$$

We now consider the definite integral $I(a,b) = \int_a^b (e^t - 1)\left[\frac{1}{t-1} - \frac{1}{t \log t}\right] dt$.

We see that

$$\int_a^b \frac{e^t}{t-1} dt = e[Ei(b-1) - Ei(a-1)]$$

and we have

$$\int_a^b e^t \frac{1}{t \log t} dt = \int_a^b e^t \frac{-1/t}{(-\log t)} dt$$

$$= e^b \log(-\log b) - e^a \log(-\log a) - \int_a^b e^t \log(-\log t) dt$$

We see that

$$\int_a^b \left[ \frac{1}{t-1} - \frac{1}{t \log t} \right] dt = \log(1-b) - \log(1-a) - \log(-\log b) + \log(-\log a)$$

and thus

$$I(a,b) = e[Ei(b-1) - Ei(a-1)] - e^b \log(-\log b) + e^a \log(-\log a) + \int_a^b e^t \log(-\log t) dt$$

$$- \log(1-b) + \log(1-a) + \log(-\log b) - \log(-\log a)$$

This gives us

$$I(0,1) = J - eEi(-1) + \int_0^1 e^t \log(-\log t) dt$$

where

$$J = \lim_{(a,b) \to (0,1)} \left[ eEi(b-1) - e^b \log(-\log b) + e^a \log(-\log a) - \log(1-b) + \log(-\log b) - \log(-\log a) \right]$$

$$= \lim_{(a,b) \to (0,1)} \left[ eEi(b-1) - \log(1-b) - (e^b - 1) \log(-\log b) + (e^a - 1) \log(-\log a) \right]$$

We note that

$$\lim_{a \to 0}(e^a - 1) \log(-\log a) = \lim_{a \to 0} \frac{(e^a - 1)}{a} a \log(-\log a) = 0$$

where we have used (3.19.4). We then have

$$J = \lim_{(a,b) \to (0,1)} \left[ eEi(b-1) - \log(1-b) - (e^b - 1) \log(-\log b) \right]$$

Noting that

$$eEi(b-1) - \log(1-b) = e[Ei(b-1) - \log(1-b)] + (e-1)\log(1-b)$$

and using (3.19.12)

$$\lim_{b \to 1}[Ei(b-1) - \log(1-b)] = \gamma$$

we obtain

$$J = e\gamma + \lim_{b \to 1}\left[(e-1)\log(1-b) - (e^b - 1)\log(-\log b)\right]$$

We see that

$$\lim_{b \to 1}\left[(e-1)\log(1-b) - (e^b - 1)\log(-\log b)\right]$$

$$= -\lim_{b \to 1}\left[\log(1-b) - \log(-\log b)\right] + \lim_{b \to 1}\left[e\log(1-b) - e^b \log(-\log b)\right]$$

and using L'Hôpital's rule we obtain

$$\lim_{b \to 1}\left[\log(1-b) - \log(-\log b)\right] = \lim_{b \to 1}\left[\log \frac{1-b}{-\log b}\right] = 0$$

We write

$$\lim_{b \to 1}\left[e\log(1-b) - e^b \log(-\log b)\right]$$

$$= e\lim_{b \to 1}\left[\log(1-b) - \log(-\log b)\right] + \lim_{b \to 1}\left[e\log(-\log b) - e^b \log(-\log b)\right]$$

$$= \lim_{b \to 1}\left[e\log(-\log b) - e^b \log(-\log b)\right]$$

$$= \lim_{b \to 1}\left[\frac{e - e^b}{1-b}(1-b)\log(-\log b)\right] = 0$$

where we have used (3.19.8.1).

We therefore have

$$I(0,1) = e\gamma - eEi(-1) + \int_0^1 e^t \log(-\log t)dt$$

Using (3.11) we obtain

(3.19.15) $$\int_0^1 e^t \log(-\log t)dt = -\sum_{k=1}^{\infty} \frac{\log k}{k!} - \gamma(e-1)$$

as shown by Medina and Moll [63] and we obtain

$$\int_0^1 (e^t - 1)\left[\frac{1}{t-1} - \frac{1}{t\log t}\right]dt = \gamma - eEi(-1) - \sum_{k=1}^{\infty} \frac{\log k}{k!}$$

Hence we obtain from (3.19.14)

(3.19.16) $$\sum_{k=1}^{\infty} \frac{H_k}{k!} = e[\gamma - Ei(-1)]$$

We see from (3.17.3) that

$$Ei(-1) = \gamma + \sum_{k=1}^{\infty} \frac{(-1)^k}{k.k!}$$

The Mathworld website for harmonic numbers reports Gosper's formula

(3.19.17) $$\sum_{k=0}^{\infty} \frac{H_k}{k!} x^k = -e^x \sum_{k=1}^{\infty} \frac{(-1)^k}{k.k!} x^k = e^x[\log x + \Gamma(0,x) + \gamma]$$

where $\Gamma(a,x)$ is the incomplete gamma function defined by

$$\Gamma(a,x) = \int_x^{\infty} t^{a-1} e^{-t} dt$$

so that

$$\Gamma(0,x) = \int_x^{\infty} \frac{e^{-t}}{t} dt$$

and we see that

$$\Gamma(0,x) = -Ei(-x)$$

With $x = 1$ (3.19.17) we obtain

$$\sum_{k=1}^{\infty}\frac{H_k}{k!} = -e\sum_{k=1}^{\infty}\frac{(-1)^k}{k.k!} = e[\gamma + \Gamma(0,1)]$$

Another derivation of (3.19.17) is given in [23]. The formula $\sum_{k=0}^{\infty}\frac{H_k}{k!}x^k = -e^x\sum_{k=1}^{\infty}\frac{(-1)^k}{k.k!}x^k$ is also contained in Ramanujan's Notebooks (see Berndt [6], Part I, p.46). This may also be easily obtained by using the Cauchy product formula together with Euler's identity

$$\sum_{k=1}^{n}\binom{n}{k}\frac{(-1)^{k+1}}{k} = H_n^{(1)}$$

Havil shows a proof of $-\sum_{k=1}^{\infty}\frac{(-1)^k}{k.k!}x^k = \log x + \Gamma(0,x) + \gamma$ in [48, p.108].

We showed in [23] that

$$\sum_{k=1}^{\infty}\frac{H_k}{k!} = e[\gamma + \Gamma'(1,1)]$$

where $\Gamma'(a,x) = \frac{\partial}{\partial a}\Gamma(a,x) = \int_x^{\infty}t^{a-1}e^{-t}\log t\,dt$ and thus

$$\Gamma'(1,1) = \int_1^{\infty}e^{-t}\log t\,dt$$

This gives us

(3.19.18) $$\sum_{k=1}^{\infty}\frac{H_k}{k!} = e\gamma + e\int_1^{\infty}e^{-t}\log t\,dt$$

We have

$$-\gamma = \int_0^{\infty}e^{-x}\log x\,dx$$

$$= \int_0^{z}e^{-x}\log x\,dx + \int_z^{\infty}e^{-x}\log x\,dx$$

$$= z\int_0^{1}e^{-zx}\log(zx)\,dx + \int_0^{\infty}e^{-(z+x)}\log(z+x)\,dx$$

and another integration by parts results in

(3.19.19) $$-\gamma = \log z + z\int_0^1 e^{-zx}\log x\,dx + e^{-z}\int_0^\infty \frac{e^{-x}}{z+x}dx$$

With $z = 1$ we obtain

(3.19.19.1) $$-\gamma = \int_0^1 e^{-x}\log x\,dx + e^{-1}\int_0^\infty \frac{e^{-x}}{1+x}dx$$

We have

$$Ei(-1) = -\int_1^\infty \frac{e^{-u}}{u}du$$

and an obvious change of variable gives us

$$Ei(-1) = -\int_0^\infty \frac{e^{-x}}{1+x}dx$$

Hence we see that

(3.19.19.2) $$-\gamma = \int_0^1 e^{-x}\log x\,dx - e^{-1}Ei(-1)$$

and combining this with (3.19.18) results in another derivation of (3.19.16).

The Gompertz constant $G$ is defined by Mathworld as

$$G = \int_0^\infty \frac{e^{-x}}{1+x}dx$$

and we see from (3.19.3) that

$$= -eEi(-1)$$

$$= 0.59634762...$$

As reported in Nielsen's book [66, p.23] we see that

$$-e\,li(e^{-1}) = \int_0^\infty \frac{e^{-x}}{1+x}\,dx$$

Differentiating (3.19.19) gives us

$$0 = \frac{1}{z} - z\int_0^1 e^{-zx} x \log x\,dx + \int_0^1 e^{-zx} \log x\,dx - e^{-z}\int_0^\infty \frac{e^{-x}}{(z+x)^2}\,dx - e^{-z}\int_0^\infty \frac{e^{-x}}{z+x}\,dx$$

and with $z=1$ we obtain

$$0 = 1 - \int_0^1 e^{-x} x \log x\,dx + \int_0^1 e^{-x} \log x\,dx - e^{-1}\int_0^\infty \frac{e^{-x}}{(1+x)^2}\,dx - e^{-1}\int_0^\infty \frac{e^{-x}}{1+x}\,dx$$

Then (3.19.19.1) gives us

$$e^{-1}\int_0^\infty \frac{e^{-x}}{(1+x)^2}\,dx = 1 + \gamma + \int_0^1 e^{-x} \log x\,dx$$

□

We note that Bailey et al. [5] have reported that

(3.20) $$\sum_{j=k}^\infty \log \frac{(j+1)^2}{j(j+2)} = \log\left(1+\frac{1}{k}\right)$$

and we make the summation

(3.21) $$\sum_{k=1}^\infty \frac{1}{k}\sum_{j=k}^\infty \log \frac{(j+1)^2}{j(j+2)} = \sum_{k=1}^\infty \frac{1}{k}\log\left(1+\frac{1}{k}\right)$$

The following identity holds whenever the series is absolutely convergent [51, p.138]

(3.22) $$\sum_{n=1}^\infty a_n \sum_{k=1}^n b_k = \sum_{n=1}^\infty b_n \sum_{k=n}^\infty a_k$$

and thus we have

$$\sum_{k=1}^\infty \frac{1}{k}\sum_{j=k}^\infty \log \frac{(1+j)^2}{j(j+2)} = \sum_{j=1}^\infty \log \frac{(1+j)^2}{j(j+2)} \sum_{k=1}^j \frac{1}{k}$$

$$= \sum_{j=1}^{\infty} H_j \log \frac{(j+1)^2}{j(j+2)}$$

$$= \sum_{j=1}^{\infty} \left[ \log \frac{j+1}{j} + \log \frac{j+1}{j+2} \right] H_j$$

$$= \sum_{j=1}^{\infty} H_j \log \frac{j+1}{j} + \sum_{j=1}^{\infty} H_j \log \frac{j+1}{j+2}$$

$$= \sum_{j=1}^{\infty} H_j \log \frac{j+1}{j} + \sum_{j=1}^{\infty} \log \frac{j+1}{j+2} \left[ H_{j+1} - \frac{1}{j+1} \right]$$

We see that

$$\sum_{j=1}^{\infty} \log \frac{j+1}{j+2} \left[ H_{j+1} - \frac{1}{j+1} \right] = \sum_{j=1}^{\infty} H_{j+1} \log \frac{j+1}{j+2} - \sum_{j=1}^{\infty} \frac{1}{j+1} \log \frac{j+1}{j+2}$$

$$= \sum_{m=2}^{\infty} H_m \log \frac{m}{m+1} - \sum_{m=2}^{\infty} \frac{1}{m} \log \frac{m}{m+1}$$

$$= -\sum_{m=2}^{\infty} H_m \log \frac{m+1}{m} + \sum_{m=2}^{\infty} \frac{1}{m} \log \frac{m+1}{m}$$

$$= -\sum_{m=1}^{\infty} H_m \log \frac{m+1}{m} + \sum_{m=1}^{\infty} \frac{1}{m} \log \frac{m+1}{m}$$

and we simply recover (3.21).

As a by-product we note that

$$\sum_{n=1}^{\infty} \frac{\log(n+1)}{n(n+1)} = \sum_{j=1}^{\infty} H_j \log \frac{(j+1)^2}{j(j+2)}$$

which corrects the sign in Coffey's formula (49) in [17].

□

We may also apply (3.22) to (2.1) but, as shown below, this does not impart any significantly new knowledge.

We write

$$\sum_{k=1}^{\infty} H_k \left( \log \frac{k+1}{k} - \frac{1}{k} \right) = \sum_{k=1}^{\infty} \left( \log \frac{k+1}{k} - \frac{1}{k} \right) \sum_{j=1}^{k} \frac{1}{j}$$

$$= \sum_{k=1}^{\infty} \frac{1}{k} \sum_{j=k}^{\infty} \left( \log \frac{j+1}{j} - \frac{1}{j} \right)$$

$$= \sum_{k=1}^{\infty} \frac{1}{k} \left[ \sum_{j=1}^{\infty} \left( \log \frac{j+1}{j} - \frac{1}{j} \right) - \sum_{j=1}^{k-1} \left( \log \frac{j+1}{j} - \frac{1}{j} \right) \right]$$

and using (4.1) this becomes

$$= -\sum_{k=1}^{\infty} \frac{1}{k} \left[ \gamma + \sum_{j=1}^{k-1} \left( \log \frac{j+1}{j} - \frac{1}{j} \right) \right]$$

$$= -\sum_{k=1}^{\infty} \frac{1}{k} \left[ \gamma + \sum_{j=1}^{k-1} \log \frac{j+1}{j} - H_{k-1} \right]$$

$$= -\sum_{k=1}^{\infty} \frac{1}{k} \left[ \gamma + \sum_{j=1}^{k-1} \log \frac{j+1}{j} - H_k + \frac{1}{k} \right]$$

$$= -\varsigma(2) - \sum_{k=1}^{\infty} \frac{1}{k} \left[ \gamma + \sum_{j=1}^{k-1} \log \frac{j+1}{j} - H_k \right]$$

We consider the finite sum

$$\sum_{k=1}^{n} \frac{1}{k} \left[ \gamma + \sum_{j=1}^{k-1} \log \frac{j+1}{j} - H_k \right] = -\sum_{k=1}^{n} \frac{H_k}{k} + \sum_{k=1}^{n} \frac{1}{k} \left[ \gamma + \sum_{j=1}^{k-1} \log \frac{j+1}{j} \right]$$

$$= -\sum_{k=1}^{n} \frac{H_k}{k} + \sum_{k=1}^{n} \frac{1}{k} \left[ \gamma + \sum_{j=1}^{k-1} \log \frac{j+1}{j} \right]$$

$$= -\sum_{k=1}^{n} \frac{H_k}{k} + \sum_{k=1}^{n} \frac{1}{k} [\gamma + \log k]$$

$$= -\left[ \sum_{k=1}^{n} \frac{H_k}{k} - \gamma \log n - \frac{1}{2} \log^2 n \right] + \sum_{k=1}^{n} \frac{\log k}{k} - \gamma \log n - \frac{1}{2} \log^2 n + \gamma H_n$$

$$= -\left[\sum_{k=1}^{n}\frac{H_k}{k} - \gamma \log n - \frac{1}{2}\log^2 n\right] + \left[\sum_{k=1}^{n}\frac{\log k}{k} - \frac{1}{2}\log^2 n\right] + \gamma[H_n - \log n]$$

Then using (2.8) and (2.9) we obtain (2.1) again

$$\sum_{k=1}^{\infty} H_k\left(\log\frac{k+1}{k} - \frac{1}{k}\right) = -\varsigma(2) + \frac{1}{2}\left[\varsigma(2) + \gamma^2\right] - \gamma_1 - \gamma^2$$

but, in this format, we have the additional benefit of seeing where the individual components come from.

$\square$

Coffey [17] also showed that

(3.23) $$\sum_{n=1}^{\infty}\frac{\log(n+1)}{n(n+1)} = -\sum_{n=1}^{\infty} H_n\left[\frac{\log(n+2)}{n+2} - \frac{\log(n+1)}{n+1}\right]$$

and an alternative derivation is shown below.

We write

$$\sum_{k=1}^{n} H_k\left[\frac{\log(k+1)}{k+1} - \frac{\log(k+2)}{k+2}\right] = \sum_{k=1}^{n+1} H_{k-1}\left[\frac{\log k}{k} - \frac{\log(k+1)}{k+1}\right]$$

$$= \sum_{k=1}^{n+1} H_k\left[\frac{\log k}{k} - \frac{\log(k+1)}{k+1}\right] - \sum_{k=1}^{n+1}\frac{1}{k}\left[\frac{\log k}{k} - \frac{\log(k+1)}{k+1}\right]$$

$$= \sum_{k=2}^{n} H_k\left[\frac{\log k}{k} - \frac{\log(k+1)}{k+1}\right] - \frac{1}{2}\log 2 + H_{n+1}\left[\frac{\log(n+1)}{n+1} - \frac{\log(n+2)}{n+2}\right]$$

$$-\sum_{k=1}^{n+1}\frac{\log k}{k^2} + \sum_{k=1}^{n+1}\frac{\log(k+1)}{k(k+1)}$$

Using (2.2) we obtain

$$\sum_{k=1}^{n} H_k\left[\frac{\log(k+1)}{k+1} - \frac{\log(k+2)}{k+2}\right]$$

$$= \sum_{k=1}^{n} \frac{\log(k+1)}{(k+1)^2} + \frac{1}{2}\log 2 - H_{n+1}\frac{\log(n+1)}{n+1}$$

$$-\frac{1}{2}\log 2 + H_{n+1}\left[\frac{\log(n+1)}{n+1} - \frac{\log(n+2)}{n+2}\right]$$

$$-\sum_{k=1}^{n+1} \frac{\log k}{k^2} + \sum_{k=1}^{n+1} \frac{\log(k+1)}{k(k+1)}$$

$$= -H_{n+1}\frac{\log(n+2)}{n+2} + \sum_{k=1}^{n+1} \frac{\log(k+1)}{k(k+1)}$$

$$= -[\log(n+1) + \gamma + \varepsilon_{n+1}]\frac{\log(n+2)}{n+2} + \sum_{k=1}^{n+1} \frac{\log(k+1)}{k(k+1)}$$

Using L'Hôpital's rule it is easy to show that

$$\lim_{n\to\infty} \frac{\log(n+1)\log(n+2)}{n+2} = 0$$

and we therefore rediscover Coffey's formula (3.23).

$\square$

A more direct derivation follows. We let $F(k) = \frac{\log(k+2)}{k+2}$ in (2.2) to obtain

$$-\sum_{k=2}^{n} H_k\left[\frac{\log(k+2)}{k+2} - \frac{\log(k+1)}{k+1}\right] = \sum_{k=1}^{n} \frac{\log(k+2)}{(k+1)(k+2)} + \frac{1}{3}\log 3 - \frac{\log(n+2)}{n+2}H_{n+1}$$

We see that

$$\sum_{k=1}^{n} \frac{\log(k+2)}{(k+1)(k+2)} = \sum_{m=2}^{n+1} \frac{\log(m+1)}{m(m+1)}$$

$$= \sum_{k=1}^{n+1} \frac{\log(k+1)}{k(k+1)} - \frac{1}{2}\log 2$$

and we have

$$\frac{\log(n+2)}{n+2}H_{n+1} = \frac{\log(n+2)}{n+2}H_{n+2} - \frac{\log(n+2)}{(n+2)^2}$$

$$= \frac{\log(n+2)}{n+2}[\varepsilon_{n+2} + \log(n+2) + \gamma] - \frac{\log(n+2)}{(n+2)^2}$$

We then obtain

$$-\sum_{k=1}^{n} H_k \left[ \frac{\log(k+2)}{k+2} - \frac{\log(k+1)}{k+1} \right] = \sum_{k=1}^{n+1} \frac{\log(k+1)}{k(k+1)} - \frac{\log(n+2)}{n+2}[\varepsilon_{n+2} + \log(n+2) + \gamma]$$

$$+ \frac{\log(n+2)}{(n+2)^2}$$

Using L'Hôpital's rule it is easy to show that

$$\lim_{n \to \infty} \frac{\log^2(n+2)}{n+2} = 0$$

and we therefore see that

$$\lim_{n \to \infty} \frac{\log(n+2)}{n+2} H_{n+1} = 0$$

Taking the limit $n \to \infty$, we therefore obtain Coffey's formula (3.23).

□

An alternative proof follows. We have

$$S = \sum_{n=1}^{\infty} H_n \left[ \frac{\log(n+2)}{n+2} - \frac{\log(n+1)}{n+1} \right] = \sum_{n=1}^{\infty} \left[ \frac{\log(n+2)}{n+2} - \frac{\log(n+1)}{n+1} \right] \sum_{k=1}^{n} \frac{1}{k}$$

$$= \sum_{n=1}^{\infty} \frac{1}{n} \sum_{k=n}^{\infty} \left[ \frac{\log(k+2)}{k+2} - \frac{\log(k+1)}{k+1} \right]$$

$$= \sum_{n=1}^{\infty} \frac{1}{n} \left\langle \sum_{k=1}^{\infty} \left[ \frac{\log(k+2)}{k+2} - \frac{\log(k+1)}{k+1} \right] - \sum_{k=1}^{n-1} \left[ \frac{\log(k+2)}{k+2} - \frac{\log(k+1)}{k+1} \right] \right\rangle$$

We have

$$\gamma_1(x) - \gamma_1(1) = \sum_{n=0}^{\infty} \left[ \frac{\log(n+x)}{n+x} - \frac{\log(n+1)}{n+1} \right]$$

and with $x = 2$, since the series telescopes, we see that

$$\gamma_1(2) - \gamma_1(1) = \sum_{n=0}^{\infty}\left[\frac{\log(n+2)}{n+2} - \frac{\log(n+1)}{n+1}\right] = 0$$

and hence we have $\gamma_1(2) = \gamma_1(1)$.

Thus we see that

$$\sum_{k=1}^{\infty}\left[\frac{\log(k+2)}{k+2} - \frac{\log(k+1)}{k+1}\right] = \sum_{k=0}^{\infty}\left[\frac{\log(k+2)}{k+2} - \frac{\log(k+1)}{k+1}\right] - \frac{1}{2}\log 2$$

$$= -\frac{1}{2}\log 2$$

and hence

$$S = -\sum_{n=1}^{\infty}\frac{1}{n}\left\langle\frac{1}{2}\log 2 + \sum_{k=1}^{n-1}\left[\frac{\log(k+2)}{k+2} - \frac{\log(k+1)}{k+1}\right]\right\rangle$$

$$= -\sum_{n=1}^{\infty}\frac{1}{n}\sum_{k=0}^{n-1}\left[\frac{\log(k+2)}{k+2} - \frac{\log(k+1)}{k+1}\right]$$

$$= -\sum_{n=1}^{\infty}\frac{1}{n}\sum_{m=1}^{n}\left[\frac{\log(m+1)}{m+1} - \frac{\log m}{m}\right]$$

$$= -\sum_{n=1}^{\infty}\frac{\log(n+1)}{n(n+1)}$$

$\square$

We now try something completely different! Guillera and Sondow [47] have shown that

(3.24) $$\psi(x) = \sum_{k=0}^{\infty}\frac{1}{k+1}\sum_{j=0}^{k}\binom{k}{j}(-1)^j \log(x+j)$$

This result was also independently obtained in a different way in [23]. We then have

(3.25) $$-\gamma = \psi(1) = \sum_{k=0}^{\infty}\frac{1}{k+1}\sum_{j=0}^{k}\binom{k}{j}(-1)^j \log(1+j)$$

Hence we have upon integration

$$\int_0^1 \frac{\psi(1+x)+\gamma}{x}\,dx = \sum_{k=0}^{\infty} \frac{1}{k+1} \sum_{j=0}^{k} \binom{k}{j}(-1)^j \int_0^1 \frac{\log[1+x/(j+1)]}{x}\,dx$$

We have the indefinite integral in terms of the dilogarithm function $Li_2(z)$

$$\int \frac{\log[1+x/(j+1)]}{x}\,dx = -Li_2\left(-\frac{x}{j+1}\right)$$

so that

$$\int_0^1 \frac{\log[1+x/(j+1)]}{x}\,dx = -Li_2\left(-\frac{1}{j+1}\right)$$

and hence we obtain

(3.25.1) $$\int_0^1 \frac{\psi(1+x)+\gamma}{x}\,dx = -\sum_{k=0}^{\infty} \frac{1}{k+1} \sum_{j=0}^{k} \binom{k}{j}(-1)^j Li_2\left(-\frac{1}{j+1}\right)$$

We see that

$$\sum_{k=0}^{\infty} \frac{1}{k+1} \sum_{j=0}^{k} \binom{k}{j}(-1)^j Li_2\left(-\frac{1}{j+1}\right) = \sum_{k=0}^{\infty} \frac{1}{k+1} \sum_{j=0}^{k} \binom{k}{j}(-1)^j \sum_{m=1}^{\infty} \frac{(-1)^m}{m^2} \frac{1}{(j+1)^m}$$

$$= \sum_{k=0}^{\infty} \frac{1}{k+1} \sum_{m=1}^{\infty} \frac{(-1)^m}{m^2} \sum_{j=0}^{k} \binom{k}{j} \frac{(-1)^j}{(j+1)^m}$$

We note from [30] that

(3.26) $$\sum_{k=0}^{n} \binom{n}{k} \frac{(-1)^k}{(k+x)^{r+1}} = \frac{\Gamma(n+1)}{(x)_{n+1}} \frac{(-1)^r}{r!} Y_r\left(-0!H_{n+1}^{(1)}(x), 1!H_{n+1}^{(2)}(x), \ldots, (-1)^r(r-1)!H_{n+1}^{(r)}(x)\right)$$

where $Y_r(x_1, \ldots, x_r)$ are the (exponential) complete Bell polynomials defined by $Y_0 = 1$ and for $r \geq 1$

$$Y_r(x_1, \ldots, x_r) = \sum_{\pi(r)} \frac{r!}{k_1! k_2! \ldots k_r!} \left(\frac{x_1}{1!}\right)^{k_1} \left(\frac{x_2}{2!}\right)^{k_2} \cdots \left(\frac{x_r}{r!}\right)^{k_r}$$

where the sum is taken over all partitions $\pi(r)$ of $r$, i.e. over all sets of integers $k_j$ such that

$$k_1 + 2k_2 + 3k_3 + \cdots + rk_r = r$$

and the ascending factorial symbol $(x)_n$, also known as the Pochhamer symbol, is defined by [79, p.16] as

$$(x)_n = x(x+1)(x+2)\cdots(x+n-1) \text{ if } n > 0$$

$$(x)_0 = 1$$

Coppo [32] has expressed (3.26) in a slightly different form

(3.27) $$\sum_{k=0}^{n}\binom{n}{k}\frac{(-1)^k}{(k+x)^{r+1}} = \frac{\Gamma(n+1)}{(x)_{n+1}}\frac{1}{r!}Y_r\left(0!H_{n+1}^{(1)}(x), 1!H_{n+1}^{(2)}(x), \ldots, (r-1)!H_{n+1}^{(r)}(x)\right)$$

and, because $Y_r(x_1, x_2, \ldots, x_r) = (-1)^r Y_r(-x_1, x_2, \ldots, (-1)^r x_r)$, we note that these statements are equivalent.

This gives us

$$\sum_{j=0}^{k}\binom{k}{j}\frac{(-1)^j}{(j+1)^m} = \frac{\Gamma(k+1)}{(1)_{k+1}}\frac{(-1)^{m-1}}{(m-1)!}Y_{m-1}\left(-0!H_{k+1}^{(1)}, 1!H_{k+1}^{(2)}, \ldots, (-1)^m(m-2)!H_{k+1}^{(m-1)}\right)$$

and we have

$$\sum_{k=0}^{\infty}\frac{1}{k+1}\sum_{m=1}^{\infty}\frac{(-1)^m}{m^2}\sum_{j=0}^{k}\binom{k}{j}\frac{(-1)^j}{(j+1)^m}$$

$$= \sum_{k=0}^{\infty}\frac{1}{k+1}\sum_{m=1}^{\infty}\frac{(-1)^m}{m^2}\frac{\Gamma(k+1)}{(1)_{k+1}}\frac{(-1)^{m-1}}{(m-1)!}Y_{m-1}\left(-0!H_{k+1}^{(1)}, 1!H_{k+1}^{(2)}, \ldots, (-1)^m(m-2)!H_{k+1}^{(m-1)}\right)$$

$$= \sum_{k=0}^{\infty}\frac{1}{k+1}\sum_{m=1}^{\infty}\frac{(-1)^m}{m^2}\frac{k!}{(k+1)!}\frac{(-1)^{m-1}}{(m-1)!}Y_{m-1}\left(-0!H_{k+1}^{(1)}, 1!H_{k+1}^{(2)}, \ldots, (-1)^m(m-2)!H_{k+1}^{(m-1)}\right)$$

$$= -\sum_{k=0}^{\infty}\frac{1}{(k+1)^2}\sum_{m=1}^{\infty}\frac{1}{m^2}\frac{1}{(m-1)!}Y_{m-1}\left(-0!H_{k+1}^{(1)}, 1!H_{k+1}^{(2)}, \ldots, (-1)^m(m-2)!H_{k+1}^{(m-1)}\right)$$

$$= -\sum_{k=0}^{\infty}\frac{1}{(k+1)^2}\sum_{m=1}^{\infty}\frac{1}{m.m!}Y_{m-1}\left(-0!H_{k+1}^{(1)}, 1!H_{k+1}^{(2)}, \ldots, (-1)^m(m-2)!H_{k+1}^{(m-1)}\right)$$

$$= -\sum_{k=1}^{\infty}\frac{1}{k^2}\sum_{m=1}^{\infty}\frac{1}{m.m!}Y_{m-1}\left(-0!H_k^{(1)}, 1!H_k^{(2)}, \ldots, (-1)^m(m-2)!H_k^{(m-1)}\right)$$

Hence we obtain

(3.28) $$\int_0^1 \frac{\psi(1+x)+\gamma}{x} dx = \sum_{k=1}^{\infty} \frac{1}{k^2} \sum_{m=1}^{\infty} \frac{1}{m.m!} Y_{m-1}\left(-0!H_k^{(1)}, 1!H_k^{(2)}, \ldots, (-1)^m(m-2)!H_k^{(m-1)}\right)$$

and we shall see a different representation of this in (6.15) below.

□

Alternatively, we consider

$$\sum_{k=0}^{\infty} \frac{1}{k+1} \sum_{m=1}^{\infty} \frac{(-1)^m}{m^2} \sum_{j=0}^{k} \binom{k}{j} \frac{(-1)^j}{(j+1)^m} = \sum_{m=1}^{\infty} \frac{(-1)^m}{m^2} \sum_{k=0}^{\infty} \frac{1}{k+1} \sum_{j=0}^{k} \binom{k}{j} \frac{(-1)^j}{(j+1)^m}$$

and we have the Hasse formula for the Hurwitz zeta function

(3.29) $$\varsigma(s,a) = \frac{1}{s-1} \sum_{k=0}^{\infty} \frac{1}{k+1} \sum_{j=0}^{k} \binom{k}{j} \frac{(-1)^j}{(j+a)^{s-1}}$$

so that

$$\sum_{k=0}^{\infty} \frac{1}{k+1} \sum_{j=0}^{k} \binom{k}{j} \frac{(-1)^j}{(j+1)^m} = m\varsigma(m+1)$$

This gives us

$$\sum_{m=1}^{\infty} \frac{(-1)^m}{m^2} \sum_{k=0}^{\infty} \frac{1}{k+1} \sum_{j=0}^{k} \binom{k}{j} \frac{(-1)^j}{(j+1)^m} = \sum_{m=1}^{\infty} \frac{(-1)^m \varsigma(m+1)}{m}$$

and we end up with

$$\int_0^1 \frac{\psi(1+x)+\gamma}{x} dx = \sum_{m=1}^{\infty} \frac{(-1)^{m+1} \varsigma(m+1)}{m}$$

which yet again shows the equivalence of (3.8) and (3.10.2)

□

Moll [64] has shown that for $a > 0$ and $n \geq 2$

(3.30) $$\int_0^1 \frac{\log^{n-1} x}{x+a} dx = (-1)^n (n-1)! Li_n\left(-\frac{1}{a}\right)$$

so that

(3.31) $$\int_0^1 \frac{\log x}{x+a}\,dx = Li_2\left(-\frac{1}{a}\right)$$

We have previously shown that

(3.32) $$\int_0^1 \frac{\psi(1+x)+\gamma}{x}\,dx = -\sum_{k=0}^{\infty}\frac{1}{k+1}\sum_{j=0}^{k}\binom{k}{j}(-1)^j Li_2\left(-\frac{1}{j+1}\right)$$

$$= -\sum_{k=0}^{\infty}\frac{1}{k+1}\sum_{j=0}^{k}\binom{k}{j}(-1)^j \int_0^1 \frac{\log x}{x+j+1}\,dx$$

$$= -\int_0^1 \sum_{k=0}^{\infty}\frac{1}{k+1}\sum_{j=0}^{k}\binom{k}{j}(-1)^j \frac{1}{x+j+1}\log x\,dx$$

The Hasse formula (3.29) gives us

$$\varsigma(2,1+x) = \sum_{k=0}^{\infty}\frac{1}{k+1}\sum_{j=0}^{k}\binom{k}{j}(-1)^j \frac{1}{x+j+1}$$

and using [79, p.22] $\varsigma(2,1+x) = \psi'(1+x)$ we obtain

$$\int_0^1 \frac{\psi(1+x)+\gamma}{x}\,dx = -\int_0^1 \psi'(1+x)\log x\,dx$$

which may also be deduced more directly using integration by parts.

We see that

$$\int_0^1 \frac{\psi(1+x)+\gamma}{x}\,dx = -\int_0^1 \psi'(1+x)\log x\,dx = -\frac{d}{ds}\int_0^1 x^s \psi'(1+x)\,dx\bigg|_{s=0}$$

and we designate $I(s)$ as

$$I(s) = \int_0^1 x^s \psi'(1+x)\,dx = x^s \psi(1+x)\bigg|_0^1 - s\int_0^1 x^{s-1}\psi(1+x)\,dx$$

$$= 1 - \gamma - s\int_0^1 x^{s-1}\psi(1+x)\,dx$$

Further analysis here may be worthwhile.

□

Alternatively, integration by parts gives us

$$\int_0^1 \frac{\psi(1+x)+\gamma}{x}dx = \frac{\log\Gamma(1+x)+\gamma x}{x}\bigg|_0^1 + \int_0^1 \frac{\log\Gamma(1+x)+\gamma x}{x^2}dx$$

L'Hôpital's rule shows that

$$\lim_{x\to 0}\frac{\log\Gamma(1+x)+\gamma x}{x} = \lim_{x\to 0}\frac{\psi(1+x)+\gamma}{1} = 0$$

and, since the integrated part is equal to $\gamma$ at $x=1$, we therefore have

$$\int_0^1 \frac{\psi(1+x)+\gamma}{x}dx = \gamma + \int_0^1 \frac{\log\Gamma(1+x)+\gamma x}{x^2}dx$$

We showed in [24] that

$$\log\Gamma(x) = \sum_{n=0}^{\infty}\frac{1}{n+1}\sum_{k=0}^{n}\binom{n}{k}(-1)^k(x+k)\log(x+k) - x + \frac{1}{2}[\log(2\pi)+1]$$

and using (3.24) we see that

(3.32.1) $\log\Gamma(1+x)+\gamma x = \sum_{n=0}^{\infty}\frac{1}{n+1}\sum_{k=0}^{n}\binom{n}{k}(-1)^k(1+k+x)\log(1+k+x) - x$

$$+\frac{1}{2}[\log(2\pi)-1] - x\sum_{n=0}^{\infty}\frac{1}{n+1}\sum_{k=0}^{n}\binom{n}{k}(-1)^k\log(1+k)$$

Since

$$(1+k+x)\log(1+k+x) - x\log(1+k)$$

$$= (1+k+x)[\log(1+k+x)-\log(1+k)] + (1+k)\log(1+k)$$

we see that

$$\log\Gamma(1+x)+\gamma x = \sum_{n=0}^{\infty}\frac{1}{n+1}\sum_{k=0}^{n}\binom{n}{k}(-1)^k(1+k+x)\log\left(1+\frac{x}{k+1}\right) - x + \frac{1}{2}[\log(2\pi)-1]$$

$$+\sum_{n=0}^{\infty}\frac{1}{n+1}\sum_{k=0}^{n}\binom{n}{k}(-1)^k(1+k)\log(1+k)$$

We have

$$\int\frac{(a+x)}{x^2}\log\left(1+\frac{x}{a}\right)dx=-Li_2\left(-\frac{x}{a}\right)-\frac{(a+x)}{x}\log(a+x)+\frac{1}{x}a\log a+\log x$$

and thus

$$\int_\varepsilon^1\frac{\log\Gamma(1+x)+\gamma x}{x^2}dx=-\sum_{n=0}^{\infty}\frac{1}{n+1}\sum_{k=0}^{n}\binom{n}{k}(-1)^k\left[Li_2\left(-\frac{1}{k+1}\right)-Li_2\left(-\frac{\varepsilon}{k+1}\right)\right]$$

$$-\sum_{n=0}^{\infty}\frac{1}{n+1}\sum_{k=0}^{n}\binom{n}{k}(-1)^k(1+k+1)\log(1+k+1)$$

$$+\frac{1}{\varepsilon}\sum_{n=0}^{\infty}\frac{1}{n+1}\sum_{k=0}^{n}\binom{n}{k}(-1)^k(1+k+\varepsilon)\log(1+k+\varepsilon)$$

$$+\left[1-\frac{1}{\varepsilon}\right]\sum_{n=0}^{\infty}\frac{1}{n+1}\sum_{k=0}^{n}\binom{n}{k}(-1)^k(1+k)\log(1+k)$$

$$-\left[1-\frac{1}{\varepsilon}\right]\sum_{n=0}^{\infty}\frac{1}{n+1}\sum_{k=0}^{n}\binom{n}{k}(-1)^k(1+k)\log(1+k)$$

$$-\log\varepsilon\sum_{n=0}^{\infty}\frac{1}{n+1}\sum_{k=0}^{n}\binom{n}{k}(-1)^k$$

$$+\frac{1}{2}[\log(2\pi)-1]\left(\frac{1}{\varepsilon}-1\right)+\log\varepsilon$$

We see that the terms involving $\log\varepsilon$ cancel and we note from (3.32.1) that

$$\frac{\sum_{n=0}^{\infty}\frac{1}{n+1}\sum_{k=0}^{n}\binom{n}{k}(-1)^k(1+k+\varepsilon)\log(1+k+\varepsilon)+\frac{1}{2}[\log(2\pi)-1]}{\varepsilon}$$

$$=\frac{\log\Gamma(1+\varepsilon)+\varepsilon(1+\gamma)}{\varepsilon}-\gamma$$

and using L'Hôpital's rule we see that as $\varepsilon \to 0$ the limit becomes $\Gamma'(1)+1$. Hence we obtain

$$\int_0^1 \frac{\log \Gamma(1+x) + \gamma x}{x^2} dx = -\sum_{n=0}^{\infty} \frac{1}{n+1} \sum_{k=0}^{n} \binom{n}{k} (-1)^k Li_2\left(-\frac{1}{k+1}\right) + 1 - \gamma$$

$$-\sum_{n=0}^{\infty} \frac{1}{n+1} \sum_{k=0}^{n} \binom{n}{k} (-1)^k (1+k+1) \log(1+k+1)$$

$$-\frac{1}{2}[\log(2\pi) - 1]$$

which simplifies to

$$\int_0^1 \frac{\log \Gamma(1+x) + \gamma x}{x^2} dx = -\sum_{n=0}^{\infty} \frac{1}{n+1} \sum_{k=0}^{n} \binom{n}{k} (-1)^k Li_2\left(-\frac{1}{k+1}\right) - \gamma$$

This is simply equivalent to (3.25.1).

$\square$

We recall Euler's identity for the dilogarithm function [79, p.108]

(3.33) $\quad -Li_2\left(-\frac{1}{z}\right) = \varsigma(2) + \frac{1}{2} \log^2 z + Li_2(-z)$

which enables us to write (3.32) as

$$\int_0^1 \frac{\psi(1+x) + \gamma}{x} dx = \sum_{k=0}^{\infty} \frac{1}{k+1} \sum_{j=0}^{k} \binom{k}{j} (-1)^j \left[\varsigma(2) + \frac{1}{2} \log^2(j+1) + Li_2(-[j+1])\right]$$

Noting that $\sum_{j=0}^{k} \binom{k}{j}(-1)^j = \delta_{0,k}$, where $\delta_{j,k}$ is the Kronecker delta, we obtain

$$\sum_{k=0}^{\infty} \frac{1}{k+1} \sum_{j=0}^{k} \binom{k}{j} (-1)^j = \sum_{k=0}^{\infty} \frac{1}{k+1} \delta_{0,k} = 1$$

We showed in [24] that the Stieltjes constants may be represented by

$$\gamma_p(u) = -\frac{1}{p+1} \sum_{n=0}^{\infty} \frac{1}{n+1} \sum_{k=0}^{n} \binom{n}{k} (-1)^k \log^{p+1}(u+k)$$

and hence we get

$$(3.34) \quad \int_0^1 \frac{\psi(1+x)+\gamma}{x}\,dx = \varsigma(2) - \gamma_1 + \sum_{k=0}^{\infty} \frac{1}{k+1} \sum_{j=0}^{k} \binom{k}{j} (-1)^j Li_2(-[j+1])$$

We see that

$$\sum_{k=0}^{\infty} \frac{1}{k+1} \sum_{j=0}^{k} \binom{k}{j} (-1)^j Li_2(-[j+1]) = \sum_{k=0}^{\infty} \frac{1}{k+1} \sum_{j=0}^{k} \binom{k}{j} (-1)^j \sum_{m=1}^{\infty} \frac{(-1)^m (j+1)^m}{m^2}$$

$$= \sum_{m=1}^{\infty} \frac{(-1)^m}{m^2} \sum_{k=0}^{\infty} \frac{1}{k+1} \sum_{j=0}^{k} \binom{k}{j} (-1)^j (j+1)^m$$

From (3.29) we have

$$-m\varsigma(1-m,a) = \sum_{k=0}^{\infty} \frac{1}{k+1} \sum_{j=0}^{k} \binom{k}{j} (-1)^j (j+a)^m$$

It is well known that [3, p.264] for integers $m \geq 1$

$$\varsigma(1-m,a) = -\frac{B_m(a)}{m}$$

where $B_m(a)$ are the Bernoulli polynomials. Hence, with $a=1$, we have in terms of the Bernoulli numbers

$$\varsigma(1-m) = -\frac{B_m}{m}$$

and thus it appears that

$$(3.35) \quad \sum_{m=1}^{\infty} \frac{(-1)^m}{m^2} \sum_{k=0}^{\infty} \frac{1}{k+1} \sum_{j=0}^{k} \binom{k}{j} (-1)^j (j+1)^m = \sum_{m=1}^{\infty} \frac{(-1)^m B_m}{m^2}$$

However, this analysis cannot be correct because the series on the right-hand side is not convergent (the error presumably arises because of the restricted range of validity of equation (3.33)).

□

Let us try a different approach. We have

$$\int_0^1 \frac{\psi(1+x)+\gamma}{x}\,dx = \varsigma(2) - \gamma_1 + \sum_{k=0}^{\infty} \frac{1}{k+1} \sum_{j=0}^{k} \binom{k}{j} (-1)^j Li_2(-[j+1])$$

and we note that

$$Li_2(z) = -\int_0^z \frac{\log(1-t)}{t} dt$$

$$= -\int_0^1 \frac{\log(1-zu)}{u} du$$

This gives us

$$S = \sum_{k=0}^{\infty} \frac{1}{k+1} \sum_{j=0}^{k} \binom{k}{j} (-1)^j Li_2(-[j+1]) = -\sum_{k=0}^{\infty} \frac{1}{k+1} \sum_{j=0}^{k} \binom{k}{j} (-1)^j \int_0^1 \frac{\log(1+[1+j]u)}{u} du$$

Since $\log(1+[1+j]u) = \log u + \log\left(j+1+\frac{1}{u}\right)$ we see from (3.24) that

$$S = -\int_0^1 \frac{\psi\left(1+\frac{1}{u}\right) + \log u}{u} du$$

and with the substitution $x = \frac{1}{u}$ we obtain

$$S = -\int_1^{\infty} \frac{\psi(1+x) - \log x}{x} dx$$

Hence we obtain

(3.35.1) $$\int_0^1 \frac{\psi(1+x) + \gamma}{x} dx + \int_1^{\infty} \frac{\psi(1+x) - \log x}{x} dx = \varsigma(2) - \gamma_1$$

This is reminiscent of two integrals derived by de Bruijn [79, p.102] in 1937 when he showed that for $0 < \text{Re}(s) < 1$ the Riemann zeta function may be expressed as

(3.35.2) $$\varsigma(s) = \frac{\sin(\pi s)}{\pi} \int_0^{\infty} \frac{\log x - \psi(1+x)}{x^s} dx$$

and also by

(3.35.3) $$\varsigma(s) = \frac{1}{s-1} + \frac{\sin(\pi s)}{\pi} \int_0^{\infty} \frac{\log(1+x) - \psi(1+x)}{x^s} dx$$

□

We slightly digress at this point by noting that Lewin [60, p.144] reports that

$$(3.36) \qquad Li_2(-m) = -\frac{1}{2}\sum_{r=1}^{m} Li_2\left(\frac{1}{r^2}\right) + \sum_{r=2}^{m} \log r \log \frac{r-1}{r}$$

and thus we have

$$\sum_{k=0}^{\infty} \frac{1}{k+1} \sum_{j=0}^{k} \binom{k}{j}(-1)^j Li_2(-[j+1]) = \sum_{k=0}^{\infty} \frac{1}{k+1} \sum_{j=0}^{k} \binom{k}{j}(-1)^j \left[ -\frac{1}{2}\sum_{r=1}^{j+1} Li_2\left(\frac{1}{r^2}\right) + \sum_{r=2}^{j+1} \log r \log \frac{r-1}{r} \right]$$

Richmond and Szekeres [73] have shown that

$$\sum_{r=2}^{\infty} L\left(\frac{1}{r^2}\right) = \varsigma(2)$$

where $L(x) = Li_2(x) + \frac{1}{2}\log x \log(1-x)$. It is easily seen that

$$L\left(\frac{1}{r^2}\right) = Li_2\left(\frac{1}{r^2}\right) + \log r \log \frac{r^2}{r^2-1}$$

and thus

$$\sum_{r=2}^{\infty} Li_2\left(\frac{1}{r^2}\right) + \sum_{r=2}^{\infty} \log r \log \frac{r^2}{r^2-1} = \varsigma(2)$$

but I do not know if this digression is germane to this paper.

□

Coffey [15] showed in 2006 that for integers $j \geq 0$ and $|z| < 1$

$$(3.37) \qquad \sum_{n=1}^{\infty} z^n \sum_{k=j}^{\infty} \frac{(-1)^k}{(k-j)!} n^{k-j} \gamma_k(u) = \sum_{n=1}^{\infty} \varsigma^{(j)}(n+1,u) z^n - (-1)^j j! Li_{j+1}(z)$$

so that with $z = j = 1$ we have

$$(3.38) \qquad \sum_{n=1}^{\infty} \frac{1}{n} \sum_{k=1}^{\infty} \frac{(-1)^k}{(k-1)!} n^k \gamma_k = \sum_{n=1}^{\infty} \varsigma'(n+1) + \varsigma(2)$$

With $j = 0$ in (3.37) we get

(3.39) $$\sum_{n=1}^{\infty} z^n \sum_{k=0}^{\infty} \frac{(-1)^k}{k!} n^k \gamma_k(u) = \sum_{n=1}^{\infty} \varsigma(n+1,u) z^n - \sum_{n=1}^{\infty} \frac{z^n}{n}$$

Dividing this by $z$ and integrating with respect to $z$ gives us

(3.40) $$\sum_{n=1}^{\infty} \frac{v^n}{n} \sum_{k=0}^{\infty} \frac{(-1)^k}{k!} n^k \gamma_k(u) = \sum_{n=1}^{\infty} \frac{\varsigma(n+1,u)}{n} v^n - \sum_{n=1}^{\infty} \frac{v^n}{n^2}$$

With $u=1$ and $v=-1$ we get

(3.41) $$\sum_{n=1}^{\infty} \frac{(-1)^n}{n} \sum_{k=0}^{\infty} \frac{(-1)^k}{k!} n^k \gamma_k = \sum_{n=1}^{\infty} \frac{(-1)^n \varsigma(n+1)}{n} - \sum_{n=1}^{\infty} \frac{(-1)^n}{n^2}$$

□

All that now remains is for someone to find a closed form expression for the series $\sum_{n=1}^{\infty} \frac{\log(n+1)}{n(n+1)}$ and success here would also result in a closed form evaluation of the aesthetically interesting integral $\int_0^1 \psi^2(1+x)\,dx$. Some other representations of this integral are considered below:

**Sine and cosine integrals**

It was shown in [31] that

(3.43) $$2\int_0^1 \psi^2(1+x)\,dx = \sum_{n=1}^{\infty} [Ci^2(2n\pi) + si^2(2n\pi)]$$

where $si(x)$ and $Ci(x)$ are the sine and cosine integrals defined [44, p.878] by

$$si(x) = -\int_x^{\infty} \frac{\sin t}{t}\,dt$$

and for $x > 0$

$$Ci(x) = -\int_x^{\infty} \frac{\cos t}{t}\,dt = \gamma + \log x + \int_0^x \frac{\cos t - 1}{t}\,dt$$

where $\gamma$ is Euler's constant.

Therefore, using (3.9) we obtain

$$(3.44) \quad \frac{1}{2}\sum_{n=1}^{\infty}[Ci^2(2n\pi)+si^2(2n\pi)] = 2\sum_{n=1}^{\infty}\frac{\log(n+1)}{n(n+1)}+1-2\varsigma(2)+2\gamma_1$$

We have from Bateman's compendium (p.146)

$$\int_0^{\infty} e^{-ax} x^{-1} \log(1+x^2)\,dx = Ci^2(a)+si^2(a)$$

and, using the geometric series to make the summation, we obtain

$$(3.45) \quad \int_0^{\infty} \frac{\log(1+x^2)}{x(e^{2\pi x}-1)}\,dx = \sum_{n=1}^{\infty}[Ci^2(2n\pi)+si^2(2n\pi)]$$

Therefore we see from (3.43) that

$$(3.46) \quad \int_0^1 \psi^2(1+x)\,dx = \frac{1}{2}\int_0^{\infty}\frac{\log(1+x^2)}{x(e^{2\pi x}-1)}\,dx$$

and we note that there is a passing resemblance to the integral representation of the Stieltjes constant $\gamma_1(u)$

$$(3.47) \quad \gamma_1(u) = \frac{1}{2u}\log u - \frac{1}{2}\log^2 u + \int_0^{\infty}\frac{x\log(u^2+x^2)}{(u^2+x^2)(e^{2\pi x}-1)}\,dx - 2u\int_0^{\infty}\frac{\tan^{-1}(x/u)}{(u^2+x^2)(e^{2\pi x}-1)}\,dx$$

which was derived by Coffey [16] (a different proof appears in [28]). This may perhaps help to explain why $\gamma_1$ appears in (3.44) above.

We also have [66, p.32]

$$-4\int_0^{\pi/2}\frac{e^{-a\tan x}\log\cos x}{\sin 2x}\,dx = Ci^2(a)+si^2(a)$$

and, using the geometric series to make the summation, we obtain

$$-4\int_0^{\pi/2}\frac{\log\cos x}{(e^{2\pi\tan x}-1)\sin 2x}\,dx = \sum_{n=1}^{\infty}[Ci^2(2n\pi)+si^2(2n\pi)]$$

It is easily seen that the substitution $u = \tan x$ results in (3.45).

The digamma function may be expressed as [31]

$$\psi(x) = \log x - \frac{1}{2x} + 2\sum_{n=1}^{\infty}[\cos(2n\pi x)Ci(2n\pi x) + \sin(2n\pi x)si(2n\pi x)]$$

and this gives us the integral

$$\int_a^u \left[\frac{\psi(1+x) - \log x}{x} - \frac{1}{2x^2}\right] dx = 2\sum_{n=1}^{\infty} \int_a^u \frac{\cos(2n\pi x)Ci(2n\pi x) + \sin(2n\pi x)si(2n\pi x)}{x} dx$$

We have

$$2\int_a^u \frac{\sin(2n\pi x)si(2n\pi x)}{x} dx = si^2(2n\pi u) - si^2(2n\pi a)$$

and

$$2\int_a^u \frac{\cos(2n\pi x)Ci(2n\pi x)}{x} dx = Ci^2(2n\pi u) - Ci^2(2n\pi a)$$

so that

$$\int_a^u \left[\frac{\psi(1+x) + \gamma}{x}\right] dx = \frac{1}{2}[\log^2 u - \log^2 a] + \gamma[\log u - \log a] - \frac{1}{2}\left(\frac{1}{u} - \frac{1}{a}\right)$$

$$+ \sum_{n=1}^{\infty}[Ci^2(2n\pi u) + si^2(2n\pi u)] - \sum_{n=1}^{\infty}[Ci^2(2n\pi a) + si^2(2n\pi a)]$$

Hence we have

$$\int_0^1 \left[\frac{\psi(1+x) + \gamma}{x}\right] dx = \sum_{n=1}^{\infty}[Ci^2(2n\pi) + si^2(2n\pi)]$$

$$+ \lim_{a \to 0}\left\{\frac{1}{2}\left(\frac{1}{a} - 1\right) - \log^2 a - \gamma \log a - \sum_{n=1}^{\infty}\left[Ci^2(2n\pi a) + si^2(2n\pi a)\right]\right\}$$

We also have (valid at infinity?)

$$\int_1^\infty \left[\frac{\psi(1+x) - \log x}{x} - \frac{1}{2x^2}\right] dx = -\sum_{n=1}^{\infty}[Ci^2(2n\pi) + si^2(2n\pi)]$$

so that

$$\int_1^\infty \frac{\psi(1+x) - \log x}{x} dx = \frac{1}{2} - \sum_{n=1}^\infty [Ci^2(2n\pi) + si^2(2n\pi)]$$

We shall see later in (3.35.1) that

$$\int_0^1 \frac{\psi(1+x) + \gamma}{x} dx + \int_1^\infty \frac{\psi(1+x) - \log x}{x} dx = \varsigma(2) - \gamma_1$$

but these formulae do not reconcile with (3.9.1).

**Bernoulli numbers of the second kind**

Coffey [18] has shown that

(3.50) $$\int_0^1 \psi^2(1+x) dx = 1 - \gamma^2 - \varsigma(2) + 2 \sum_{n=1}^\infty \frac{(-1)^{n+1} b_n}{n^2 n!}$$

where $b_n$ are the Bernoulli numbers of the second kind which appear in the generating function

$$\frac{x}{\log(1+x)} = \sum_{n=0}^\infty \frac{b_n}{n!} x^n$$

With L'Hôpital's rule we see that $b_0 = 1$. We see that

$$\frac{x}{\log(1-x)} = \sum_{n=0}^\infty \frac{(-1)^{n+1} b_n}{n!} x^n$$

and thus

(3.51) $$\frac{1}{\log(1-x)} + \frac{1}{x} = \sum_{n=1}^\infty \frac{(-1)^{n+1} b_n}{n!} x^{n-1}$$

We have the integral

$$\int_0^t \left[ \frac{1}{\log(1-x)} + \frac{1}{x} \right] dx = \sum_{n=1}^\infty \frac{(-1)^{n+1} b_n}{n.n!} t^n$$

and

$$\int_0^u \frac{dt}{t} \int_0^t \left[ \frac{1}{\log(1-x)} + \frac{1}{x} \right] dx = \sum_{n=1}^\infty \frac{(-1)^{n+1} b_n}{n^2 n!} u^n$$

Letting $x = ty$ we have

$$\int_0^t \left[\frac{1}{\log(1-x)} + \frac{1}{x}\right] dx = \int_0^1 \left[\frac{1}{\log(1-ty)} + \frac{1}{ty}\right] t \, dy$$

and thus

$$\int_0^u dt \int_0^1 \left[\frac{1}{\log(1-ty)} + \frac{1}{ty}\right] dy = \sum_{n=1}^{\infty} \frac{(-1)^{n+1} b_n}{n^2 n!} u^n$$

We then have using (3.50)

(3.52) $$\int_0^1 \psi^2(1+x) \, dx = 1 - \gamma^2 - \varsigma(2) + 2 \int_0^1 \int_0^1 \left[\frac{1}{\log(1-ty)} + \frac{1}{ty}\right] dy \, dt$$

It may be noted that

(3.52.1) $$\gamma = \int_0^1 \left[\frac{1}{\log(1-x)} + \frac{1}{x}\right] dx = \sum_{n=1}^{\infty} \frac{(-1)^{n+1} b_n}{n.n!}$$

Some further applications of the Bernoulli numbers of the second kind are given in [33].

Glasser and Manna [43] have evaluated the following series which is <u>similar</u> to the one that appears in Coffey's formula (3.50)

(3.53) $$\sum_{n=1}^{\infty} \frac{(-1)^n b_{n+1}}{n(n+1)!} = \frac{1}{2}[\log(2\pi) - 1 - \gamma]$$

but, as we all know, evaluating $\varsigma(2)$ is not the same as evaluating $\varsigma(3)$ !

We see from (3.51) that

$$\frac{1}{x \log(1-x)} + \frac{1}{x^2} - \frac{1}{2x} = \sum_{n=2}^{\infty} \frac{(-1)^{n+1} b_n}{n!} x^{n-2}$$

whereupon integration gives us

$$\int_0^1 \left[\frac{1}{\log(1-x)} + \frac{1}{x^2} - \frac{1}{2x}\right] dx = \sum_{n=2}^{\infty} \frac{(-1)^{n+1} b_n}{(n-1).n!}$$

so that we have

$$\int_0^1 \left[\frac{1}{\log(1-x)} + \frac{1}{x^2} - \frac{1}{2x}\right] dx = \frac{1}{2}[\log(2\pi) - 1 - \gamma]$$

We multiply (3.51) by $\log(1-x)$ to get

(3.54) $$1+\frac{\log(1-x)}{x} = \sum_{n=1}^{\infty} \frac{(-1)^{n+1} b_n}{n!} x^{n-1} \log(1-x)$$

and integration gives us

$$1+\int_0^1 \frac{\log(1-x)}{x} dx = \sum_{n=1}^{\infty} \frac{(-1)^{n+1} b_n}{n!} \int_0^1 x^{n-1} \log(1-x) dx$$

Using $\int_0^1 x^{n-1} \log(1-x) dx = -\frac{H_n}{n}$ gives us the known result [33]

(3.55) $$\varsigma(2) - 1 = \sum_{n=1}^{\infty} \frac{(-1)^{n+1} b_n H_n}{n \cdot n!}$$

Letting $x \to tx$ in (3.54) gives us

$$1+\int_0^1 \frac{\log(1-tx)}{tx} dx = \sum_{n=1}^{\infty} \frac{(-1)^{n+1} b_n t^{n-1}}{n!} \int_0^1 x^{n-1} \log(1-tx) dx$$

We have

$$\int_0^1 \frac{\log(1-tx)}{tx} dx = \int_0^t \frac{\log(1-u)}{u} du = -Li_2(t)$$

and *Mathematica* gives us

$$\int x^{n-1} \log(1-tx) dx = \frac{x^n}{n(n+1)} [{}_2F_1(1, n+1; n+2; tx) + (n+1)\log(1-tx)]$$

so that

$$\int_0^1 x^{n-1} \log(1-tx) dx = \frac{1}{n(n+1)} [{}_2F_1(1, n+1; n+2; t) + (n+1)\log(1-t)]$$

Hence we obtain

$$1 - Li_2(t) = \sum_{n=1}^{\infty} \frac{(-1)^{n+1} b_n t^{n-1}}{n!} \frac{1}{n(n+1)} [{}_2F_1(1, n+1; n+2; t) + (n+1)\log(1-t)]$$

Following [42] we note that

(3.56) $$\int_0^1 (1+x)^t \, dt = \int_0^1 \exp[t \log(1+x)] \, dt$$

$$= \frac{x}{\log(1+x)}$$

We have

$$(1+x)^t = \sum_{n=0}^{\infty} \binom{t}{n} x^n = \sum_{n=0}^{\infty} \frac{\Gamma(1+t)}{\Gamma(1+t-n)} \frac{x^n}{n!}$$

which implies that

$$b_n = n! \int_0^1 \binom{t}{n} dt = \int_0^1 t(t-1)\ldots(t-n+1) \, dt$$

With the substitution $t = 1 - u$, we see that for $n \geq 1$

$$b_n = (-1)^{n-1} \int_0^1 (1-u)u(u+1)\ldots(u+n-2) \, du$$

Since the integrand is positive, this implies that $(-1)^{n-1} b_n > 0$ and thus $b_n$ alternate in sign.

In terms of the Stirling numbers of the first kind we have

$$u(u+1)\ldots(u+n-2) = (-1)^n \sum_{k=0}^{n-1} (-1)^k s(n-1,k) u^k$$

and using

$$\int_0^1 (1-u) u^k \, du = \frac{1}{(k+1)(k+2)}$$

we obtain

$$b_n = -\sum_{k=0}^{n-1} (-1)^k \frac{s(n-1,k)}{(k+1)(k+2)}$$

We also have

$$b_n = \int_0^1 t(t-1)...(t-n+1)dt$$

$$= \int_0^1 \sum_{j=0}^n s(k,j)t^j dt$$

and thus we have

$$b_n = \sum_{j=0}^n \frac{s(k,j)}{j+1}$$

Coppo [33] refers to the non-alternating Cauchy numbers (of the first kind) defined by

$$\lambda_n = \int_0^1 t(1-t)...(n-1-t)dt$$

$$\lambda_1 = 1/2$$

and it is easily seen that $\lambda_n = (-1)^{n-1} b_n$.

Yingying [81] showed in 2002 that

$$\sum_{n=k+1}^\infty \int_0^1 \frac{x}{n(n+x)} dx = \sum_{n=1}^\infty \frac{a_n}{(k+1)...(k+n)}$$

where

$$a_1 = \frac{1}{2} \qquad a_n = \frac{1}{n}\int_0^1 t(1-t)...(n-1-t)dt = \frac{(-1)^{n+1}}{n}\sum_{j=1}^n \frac{s(n,j)}{j+1}$$

□

Using (3.56) we see that

$$\frac{1}{\log(1-x)} + \frac{1}{x} = \int_0^1 \left[-\frac{(1-x)^t}{x} + t^{x-1}\right] dt$$

and we have the double integral

$$\int_a^1 \left[\frac{1}{\log(1-x)}+\frac{1}{x}\right]dx = \int_a^1\int_0^1\left[-\frac{(1-x)^t}{x}+t^{x-1}\right]dt\,dx$$

We see that

$$\int_a^1 \frac{(1-x)^t}{x}dx = \int_a^1 \exp[t\log(1-x)]\frac{dx}{x}$$

$$= \sum_{n=0}^\infty \frac{t^n}{n!}\int_a^1 \frac{\log^n(1-x)}{x}dx$$

$$= -\log a + \sum_{n=1}^\infty \frac{t^n}{n!}\int_a^1 \frac{\log^n(1-x)}{x}dx$$

We have

$$\int_a^1 t^{x-1}dx = \frac{1-t^{a-1}}{\log t}$$

and it is known that [10, p.98]

$$\int_0^1 \frac{1-t^{a-1}}{\log t}dt = -\log a$$

This gives us

(3.57) $$\int_a^1\left[\frac{1}{\log(1-x)}+\frac{1}{x}\right]dx = -\sum_{n=1}^\infty \frac{1}{(n+1)n!}\int_a^1 \frac{\log^n(1-x)}{x}dx$$

and using (3.18) we obtain the limit as $a \to 0$

(3.58) $$\int_0^1\left[\frac{1}{\log(1-x)}+\frac{1}{x}\right]dx = \sum_{n=1}^\infty (-1)^{n+1}\frac{\varsigma(n+1)}{n+1}$$

Referring to (6.1) we obtain the known result

(3.59) $$\gamma = \int_0^1\left[\frac{1}{\log(1-x)}+\frac{1}{x}\right]dx = \sum_{n=2}^\infty \frac{(-1)^n \varsigma(n)}{n}$$

Differentiating (3.57) with respect to $a$ gives us

$$\frac{1}{\log(1-x)} + \frac{1}{x} = -\sum_{n=1}^{\infty} \frac{1}{(n+1)!} \frac{\log^n(1-x)}{x}$$

$$= -\frac{1}{x\log(1-x)} \sum_{n=1}^{\infty} \frac{\log^{n+1}(1-x)}{(n+1)!}$$

$$= -\frac{1}{x\log(1-x)} \sum_{n=2}^{\infty} \frac{\log^n(1-x)}{n!}$$

$$= -\frac{1}{x\log(1-x)} \left[ \sum_{n=0}^{\infty} \frac{\log^n(1-x)}{n!} - 1 - \log(1-x) \right]$$

$$= -\frac{1}{x\log(1-x)} [1 - x - 1 - \log(1-x)]$$

$$= \frac{1}{x\log(1-x)} [x + \log(1-x)]$$

$$= \frac{1}{\log(1-x)} + \frac{1}{x}$$

and this indirectly verifies (3.57).

$\square$

Another derivation of (3.50) is shown below. We write (2.16.2) in the form

$$\int_0^1 \left[ \frac{\log(1-x)}{1-x} + \frac{\log(1-x)}{x\log x} \right] dx = \int_0^1 \left[ \frac{\log(1-x)}{1-x} + \frac{\log(1-x)}{\log x} + \frac{(1-x)\log(1-x)}{x\log x} \right] dx$$

$$= \frac{1}{2} \left[ \varsigma(2) - \gamma^2 \right] - \gamma_1$$

and using (3.10.1) this becomes

$$\int_0^1 \left[ \frac{1}{1-x} + \frac{1}{\log x} \right] \log(1-x)\, dx = \frac{1}{2} \left[ \varsigma(2) - \gamma^2 \right] - \gamma_1 - \sum_{n=1}^{\infty} \frac{\log(n+1)}{n(n+1)}$$

It is easily seen that

$$\int_0^1 \left[\frac{1}{1-x} + \frac{1}{\log x}\right] \log(1-x)\,dx = \int_0^1 \left[\frac{1}{x} + \frac{1}{\log(1-x)}\right] \log x\,dx$$

$$= \sum_{n=1}^{\infty} \frac{(-1)^{n+1} b_n}{n!} \int_0^1 x^{n-1} \log x\,dx$$

where we have used (3.51).

Integration by parts shows that

$$\int_0^1 x^{n-1} \log x\,dx = -\frac{1}{n^2}$$

and hence we have

(3.60) $\quad \dfrac{1}{2}\left[\varsigma(2) - \gamma^2\right] - \gamma_1 - \displaystyle\sum_{n=1}^{\infty} \frac{\log(n+1)}{n(n+1)} = -\sum_{n=1}^{\infty} \frac{(-1)^{n+1} b_n}{n^2 n!}$

One may then obtain Coffey's formula by substituting (3.9).

□

We recall (2.16.4)

$$\int_0^1 \left[\frac{1}{1-x} + \frac{1}{\log x}\right] \log(-\log x)\,dx = -\gamma^2 - \gamma_1$$

and substituting

$$\frac{1}{\log x} + \frac{1}{1-x} = \sum_{n=1}^{\infty} \frac{(-1)^{n+1} b_n}{n!} (1-x)^{n-1}$$

gives us

$$\int_0^1 \left[\frac{1}{\log x} + \frac{1}{1-x}\right] \log(-\log x)\,dx = \sum_{n=1}^{\infty} \frac{(-1)^{n+1} b_n}{n!} \int_0^1 (1-x)^{n-1} \log(-\log x)\,dx$$

$$= \sum_{n=1}^{\infty} \frac{(-1)^{n+1} b_n}{n!} \sum_{k=0}^{n-1} \binom{n-1}{k} (-1)^k \int_0^1 x^k \log(-\log x)\,dx$$

Using (3.11) this becomes

$$= \sum_{n=1}^{\infty} \frac{(-1)^{n+1} b_n}{n!} \sum_{k=0}^{n-1} \binom{n-1}{k} (-1)^{k+1} \frac{\gamma + \log(k+1)}{k+1}$$

We see that

$$\sum_{k=0}^{n-1} \binom{n-1}{k} (-1)^{k+1} \frac{\gamma + \log(k+1)}{k+1} = \sum_{m=1}^{n} \binom{n-1}{m-1} (-1)^m \frac{\gamma + \log m}{m}$$

$$= \sum_{m=1}^{n} \binom{n}{m} \frac{m}{n} (-1)^m \frac{\gamma + \log m}{m}$$

$$= \frac{1}{n} \sum_{m=1}^{n} \binom{n}{m} (-1)^m (\gamma + \log m)$$

and hence we have

$$\int_0^1 \left[ \frac{1}{1-x} + \frac{1}{\log x} \right] \log(-\log x)\, dx = \sum_{n=1}^{\infty} \frac{(-1)^{n+1} b_n}{n \cdot n!} \sum_{m=1}^{n} \binom{n}{m} (-1)^m (\gamma + \log m)$$

$$= -\gamma^2 - \gamma_1$$

We note that

$$\sum_{m=1}^{n} \binom{n}{m} (-1)^m (\gamma + \log m) = \gamma \sum_{m=1}^{n} \binom{n}{m} (-1)^m + \sum_{m=1}^{n} \binom{n}{m} (-1)^m \log m$$

$$= -\gamma + \sum_{m=1}^{n} \binom{n}{m} (-1)^m \log m$$

and hence using (3.52.1) we obtain

(3.61) $$\sum_{n=1}^{\infty} \frac{(-1)^{n+1} b_n}{n \cdot n!} \sum_{m=1}^{n} \binom{n}{m} (-1)^m \log m = -\gamma_1$$

It may be noted that

(3.62) $$\sum_{m=1}^{n} \binom{n}{m} (-1)^m \log m = \frac{d}{ds} \sum_{m=1}^{n} \binom{n}{m} (-1)^m m^s \bigg|_{s=0}$$

and the loose connection with the Stirling numbers of the second kind

$$S(l,n) = \frac{1}{n!} \sum_{m=1}^{n} \binom{n}{m} (-1)^{n-m} m^l$$

Having regard to (3.52.1), we conjecture that

(3.63) $$\sum_{n=1}^{\infty} \frac{(-1)^{n+1} b_n}{n.n!} \sum_{m=1}^{n} \binom{n}{m} (-1)^m \log^p m = -\gamma_p$$

$\square$

Lewin [60, p.22] showed with integration by parts that

$$\int_0^t x^\alpha Li_2(x)\, dx = \frac{t^{\alpha+1}}{\alpha+1} Li_2(t) + \int_0^t \frac{x^{\alpha+1}}{\alpha+1} \frac{\log(1-x)}{x}\, dx$$

We have

$$\int_0^t \frac{x^\alpha}{\alpha+1} \log(1-x)\, dx = \frac{t^{\alpha+1}}{(\alpha+1)^2} \log(1-t) + \frac{1}{(\alpha+1)^2} \int_0^t \frac{x^{\alpha+1}}{1-x}\, dx$$

and noting that

$$\int_0^t \frac{x^{\alpha+1}}{1-x}\, dx = -\int_0^t \frac{1-x^{\alpha+1}}{1-x}\, dx + \log(1-t)$$

we may write this as

(3.65) $$\int_0^t x^\alpha Li_2(x)\, dx = \frac{t^{\alpha+1}}{\alpha+1} Li_2(t) + \frac{t^{\alpha+1}-1}{(\alpha+1)^2} \log(1-t) - \frac{1}{(\alpha+1)^2} \int_0^t \frac{1-x^{\alpha+1}}{1-x}\, dx$$

With $t=1$ we have

(3.66) $$\int_0^1 x^\alpha Li_2(x)\, dx = \frac{\varsigma(2)}{\alpha+1} - \frac{\psi(\alpha+2)+\gamma}{(\alpha+1)^2}$$

where in the final part we have used

$$\psi(\alpha) + \gamma = \int_0^1 \frac{1-x^{\alpha-1}}{1-x}\, dx$$

The formula (3.65) was used extensively in [25] to evaluate integrals involving polylogarithms.

We may write (3.66) as

$$\frac{\psi(1+x)+\gamma}{x} = \varsigma(2) - x\int_0^1 t^{x-1} Li_2(t)\,dt$$

and we have via integration by parts

$$\int_0^u xt^{x-1}\,dx = \frac{ut^{u-1}}{\log t} - \int_0^u \frac{t^{x-1}}{\log t}\,dx$$

$$= \frac{ut^{u-1}}{\log t} - \frac{ut^{u-1} - t^{-1}}{\log^2 t}$$

This gives us

$$\int_0^1 xt^{x-1}\,dx = \frac{1}{\log t} - \frac{1-t^{-1}}{\log^2 t} = \frac{t\log t - t + 1}{t\log^2 t}$$

and thus we have

(3.67) $$\int_0^1 \frac{\psi(1+x)+\gamma}{x}\,dx = \varsigma(2) - \int_0^1 \frac{t\log t - t + 1}{t\log^2 t} Li_2(t)\,dt$$

□

Following Lewin's procedure [60, p.22], an alternative integration by parts shows that

$$\int_0^t x^\alpha Li_p(x)\,dx = \frac{t^{\alpha+1}}{\alpha+1} Li_p(t) - \frac{1}{\alpha+1}\int_0^t x^{\alpha+1} \frac{Li_{p-1}(x)}{x}\,dx$$

$$= \frac{t^{\alpha+1}}{\alpha+1} Li_p(t) - \frac{1}{\alpha+1}\int_0^t x^\alpha Li_{p-1}(x)\,dx$$

$$= \frac{t^{\alpha+1}}{\alpha+1} Li_p(t) - \frac{t^{\alpha+1}}{(\alpha+1)^2} Li_{p-1}(t) + \frac{1}{(\alpha+1)^2}\int_0^t x^\alpha Li_{p-2}(x)\,dx$$

and hence we have

$$\int_0^t x^\alpha Li_p(x)\,dx = t^{\alpha+1}\sum_{k=0}^{p-2}(-1)^k \frac{Li_{p-k}(t)}{(\alpha+1)^{k+1}} + \frac{(-1)^p}{(\alpha+1)^{p-1}}\int_0^t x^{\alpha+1} \frac{\log(1-x)}{x}\,dx$$

$$= t^{\alpha+1} \sum_{k=0}^{p-2} (-1)^k \frac{Li_{p-k}(t)}{(\alpha+1)^{k+1}} + (-1)^p \frac{t^{\alpha+1}-1}{(\alpha+1)^p} \log(1-t) - \frac{(-1)^p}{(\alpha+1)^p} \int_0^t \frac{1-x^{\alpha+1}}{1-x} dx$$

Therefore, letting $t=1$ we get

(3.68) $$\int_0^1 x^\alpha Li_p(x) \, dx = \sum_{k=0}^{p-2} (-1)^k \frac{\varsigma(p-k)}{(\alpha+1)^{k+1}} - \frac{(-1)^p}{(\alpha+1)^p} [\psi(\alpha+2)+\gamma]$$

With $\alpha = n-1$ we have

$$\int_0^1 x^n \frac{Li_p(x)}{x} dx = \sum_{k=0}^{p-2} (-1)^k \frac{\varsigma(p-k)}{n^{k+1}} - (-1)^p \frac{H_n^{(1)}}{n^p}$$

and with $\alpha = 0$ we obtain

(3.69) $$\int_0^1 Li_p(x) \, dx = \sum_{k=0}^{p-2} (-1)^k \varsigma(p-k) + (-1)^{p+1}$$

Integrating (4.13) gives us

$$\sum_{n=1}^\infty \frac{x^{n+1}}{n(n+1)} \log\left(1+\frac{1}{n}\right) = -\sum_{n=1}^\infty \frac{(-1)^n}{n} \int_0^x Li_{n+1}(u) \, du$$

By trial and error we see that for $n \geq 2$

(3.70) $$\int_0^x Li_n(u) \, du = (-1)^n [(x-1)\log(1-x) - x] + (-1)^n x \sum_{j=2}^n (-1)^j Li_j(x)$$

and we obtain

$$\sum_{n=1}^\infty \frac{x^{n+1}}{n(n+1)} \log\left(1+\frac{1}{n}\right) = \sum_{n=1}^\infty \frac{1}{n} \left[ (x-1)\log(1-x) - x \right] + x \sum_{j=2}^{n+1} (-1)^j Li_j(x) \right]$$

but, prima facie, this does not appear to be convergent.

□

Integration by parts gives us

$$\int_0^u \frac{\psi(1+x)+\gamma}{x} dx = [\psi(1+x)+\gamma] \log x \Big|_0^u - \int_0^u \psi'(1+x) \log x \, dx$$

We see that

$$\lim_{x\to 0}[\psi(1+x)+\gamma]\log x = \lim_{x\to 0}\frac{\psi(1+x)+\gamma}{x}\lim_{x\to 0} x\log x$$

and L'Hôpital's rule shows that the limit is zero. Therefore we have

$$\int_0^u \frac{\psi(1+x)+\gamma}{x}\,dx = [\psi(1+u)+\gamma]\log u - \int_0^u \psi'(1+x)\log x\,dx$$

$$= [\psi(1+u)+\gamma]\log u - \sum_{n=1}^{\infty}\int_0^u \frac{\log x}{(x+n)^2}\,dx$$

We have

$$\int \frac{\log x}{(x+n)^2}\,dx = \frac{x\log x - (x+n)\log(x+n)}{n(x+n)}$$

so that

$$\int_0^u \frac{\log x}{(x+n)^2}\,dx = \frac{u\log u}{n(u+n)} - \frac{1}{n}\log\left(1+\frac{u}{n}\right)$$

Hence we have

$$\int_0^u \frac{\psi(1+x)+\gamma}{x}\,dx = [\psi(1+u)+\gamma]\log u - \log u\sum_{n=1}^{\infty}\frac{u}{n(u+n)} + \sum_{n=1}^{\infty}\frac{1}{n}\log\left(1+\frac{u}{n}\right)$$

and, using (3.13), this becomes

$$\int_0^u \frac{\psi(1+x)+\gamma}{x}\,dx = \sum_{n=1}^{\infty}\frac{1}{n}\log\left(1+\frac{u}{n}\right)$$

so that

$$\int_0^1 \frac{\psi(1+x)+\gamma}{x}\,dx = \sum_{n=1}^{\infty}\frac{1}{n}\log\left(1+\frac{1}{n}\right)$$

□

We have the Maclaurin expansion

$$\psi(x+1)+\gamma = \frac{1}{x}\sum_{k=2}^{\infty}(-1)^k \varsigma(k)x^k$$

and using the Cauchy product formula we obtain

$$[\psi(x+1)+\gamma]^2 = \frac{1}{x^2}\sum_{n=4}^{\infty}\sum_{k=2}^{n-2}\varsigma(k)\varsigma(n-k)(-1)^n x^n$$

as previously noted by Coffey [20]. Integration gives us

$$\int_0^u [\psi(x+1)+\gamma]^2 dx = \sum_{n=4}^{\infty}\frac{(-1)^n u^{n-1}}{n-1}\sum_{k=2}^{n-2}\varsigma(k)\varsigma(n-k)$$

and, in particular, we have

$$\int_0^1 [\psi(x+1)+\gamma]^2 dx = \sum_{n=4}^{\infty}\frac{(-1)^n}{n-1}\sum_{k=2}^{n-2}\varsigma(k)\varsigma(n-k)$$

A multitude of integrals and series related to $\sum_{n=1}^{\infty}\frac{\log(n+1)}{n(n+1)}$ has been displayed above. However, I fear that the evaluation of this series may prove to be as elusive as the related series which results from the integral

$$\int_0^1 \psi(x)\sin\pi x\, dx = -\frac{2}{\pi}\left[\log(2\pi)+\gamma+2\sum_{n=1}^{\infty}\frac{\log n}{4n^2-1}\right]$$

Kölbig [55] states that the infinite series in (4.12) "does not seem to be expressible in terms of well-known functions".

We have

$$2\sum_{n=1}^{\infty}\frac{\log n}{4n^2-1} = \sum_{n=1}^{\infty}\left[\frac{1}{2n-1}-\frac{1}{2n+1}\right]\log n$$

and consider the finite sum

$$\sum_{n=1}^{N}\frac{\log n}{2n-1} = \sum_{m=0}^{N-1}\frac{\log(m+1)}{2m+1}$$

$$= \sum_{n=1}^{N}\frac{\log(n+1)}{2n+1} - \frac{\log(N+1)}{2N+1}$$

Hence we have

$$\sum_{n=1}^{N}\left[\frac{1}{2n-1}-\frac{1}{2n+1}\right]\log n = \sum_{n=1}^{N}\left[\frac{\log(n+1)}{2n+1}-\frac{\log n}{2n+1}\right] - \frac{\log(N+1)}{2N+1}$$

Therefore, as $N \to \infty$ we see that

$$2\sum_{n=1}^{\infty}\frac{\log n}{4n^2-1} = \sum_{n=1}^{\infty}\frac{1}{2n+1}\log\frac{n+1}{n}$$

$$= \sum_{n=1}^{\infty}\frac{1}{2n+1}\log\left(1+\frac{1}{n}\right)$$

Using

$$\log\frac{n+1}{n} = \int_0^1 \frac{y^n - y^{n-1}}{\log y}dy = \int_0^1 \frac{y^{n-1}(y-1)}{\log y}dy$$

we get

$$\sum_{n=1}^{\infty}\frac{1}{2n+1}\log\frac{n+1}{n} = \sum_{n=1}^{\infty}\int_0^1 \frac{y^n}{2n+1}\frac{(y-1)}{y\log y}dy$$

$$= \int_0^1 \sum_{n=1}^{\infty}\frac{y^n}{2n+1}\frac{(y-1)}{y\log y}dy$$

We have

$$\sum_{n=1}^{\infty}\frac{y^n}{2n+1} = \sum_{n=0}^{\infty}\frac{y^n}{2n+1} - 1$$

$$= \frac{1}{\sqrt{y}}\sum_{n=0}^{\infty}\frac{\left(\sqrt{y}\right)^{2n+1}}{2n+1} - 1$$

and thus

$$\sum_{n=1}^{\infty}\frac{1}{2n+1}\log\frac{n+1}{n} = \int_0^1 \left[\frac{1}{2\sqrt{y}}\log\frac{1+\sqrt{y}}{1-\sqrt{y}} - 1\right]\frac{y-1}{y\log y}dy$$

$$= \int_0^1 \left[\frac{1}{2x}\log\frac{1+x}{1-x} - 1\right]\frac{x^2-1}{x\log x}dx$$

**4. The generalised Euler constant function $\gamma(x)$**

Euler's constant $\gamma$ may be defined as the limit of the sequence

$$\gamma = \lim_{n\to\infty}[H_n - \log n]$$

$$= \lim_{n\to\infty}[H_n - \log(n+1) + \log(n+1) - \log n]$$

$$= \lim_{n\to\infty}[H_n - \log(n+1)] + \lim_{n\to\infty} \log\left(1 + \frac{1}{n}\right)$$

$$= \lim_{n\to\infty}[H_n - \log(n+1)]$$

We have

$$\sum_{k=1}^{n} \log\left(1 + \frac{1}{k}\right) = \sum_{k=1}^{n} [\log(k+1) - \log k]$$

and it is easily seen that this telescopes to

$$\sum_{k=1}^{n} \log\left(1 + \frac{1}{k}\right) = \log(n+1)$$

We can therefore write Euler's constant $\gamma$ as the infinite series

(4.1) $$\gamma = \sum_{k=1}^{\infty} \left[\frac{1}{k} - \log\left(1 + \frac{1}{k}\right)\right]$$

This may also be obtained by letting $x = 0$ in the expansion

(4.1.1) $$\log \Gamma(1+x) = -\log(1+x) - \gamma(1+x) + \sum_{k=1}^{\infty}\left[\log k - \log(k+1+x) + \frac{1+x}{k}\right]$$

obtained from the Weierstrass canonical form of the gamma function [79, p.1]

Sondow [77] discovered an alternating series similar to (4.1)

(4.2) $$\log\frac{4}{\pi} = \sum_{k=1}^{\infty}(-1)^{k-1}\left[\frac{1}{k} - \log\left(1 + \frac{1}{k}\right)\right]$$

$$= \sum_{k=1}^{\infty}\frac{(-1)^{k-1}}{k} - \sum_{k=1}^{\infty}(-1)^{k-1}\log\left(1+\frac{1}{k}\right)$$

$$= \log 2 - \log\frac{\pi}{2} = \log\frac{4}{\pi}$$

where, in the last line, we have employed

$$\log\frac{\pi}{2} = \sum_{k=1}^{\infty}(-1)^{k-1}\log\left(1+\frac{1}{k}\right)$$

being the logarithm of Wallis's product formula for $\pi/2$ which was published in Arithmetica Infinitorum in 1659

$$\frac{\pi}{2} = \prod_{n=1}^{\infty}\frac{(2n)^2}{(2n-1)(2n+1)}$$

It is well known that Wallis's product formula may be obtained from Euler's infinite product for the $\sin z$ function.

Sondow [77] also noted that $\gamma$ and $\log\frac{4}{\pi}$ are related by Euler's formula

(4.3) $$\gamma - \log\frac{4}{\pi} = 2\sum_{n=2}^{\infty}(-1)^n\frac{\varsigma(n)}{n2^n}$$

and this concurs with the Maclaurin expansion (6.1) of $\log\Gamma(1+x)$ in the case where $x = 1/2$.

This may also be obtained by letting $x = -1/2$ in (4.1.1)

$$\log\Gamma(1/2) = \log 2 - \frac{1}{2}\gamma + \sum_{k=1}^{\infty}\left[-\log\left(1+\frac{1}{2k}\right)+\frac{1}{2k}\right]$$

$$= \log 2 - \frac{1}{2}\gamma + \sum_{k=1}^{\infty}\left[\sum_{m=1}^{\infty}\frac{(-1)^m}{m(2k)^m}+\frac{1}{2k}\right]$$

$$= \log 2 - \frac{1}{2}\gamma + \sum_{k=1}^{\infty}\left[\sum_{m=2}^{\infty}\frac{(-1)^m}{m(2k)^m}\right]$$

$$= \log 2 - \frac{1}{2}\gamma + \left[\sum_{m=2}^{\infty}\frac{(-1)^m\varsigma(m)}{m2^m}\right]$$

Sondow's formula (4.2) may also be obtained from Lerch's series expansion for the digamma function for $0 < x < 1$ (see for example [46, p.105], [58] and [67, p.204])

(4.4) $$\psi(x)\sin\pi x + \frac{\pi}{2}\cos\pi x + (\gamma + \log 2\pi)\sin\pi x = -\sum_{n=1}^{\infty}\sin(2n+1)\pi x.\log\frac{n+1}{n}$$

Letting $x = 1/2$ in (4.4) we obtain

$$\psi(1/2) + \gamma + \log 2\pi = \sum_{n=1}^{\infty} (-1)^{n-1} \log \frac{n+1}{n}$$

and, since [79, p.20] $\psi(1/2) = -\gamma - 2\log 2$, (4.2) follows automatically.

Having regard to the structure of (4.2) and (4.3), Sondow and Hadjicostas [78] considered the generalised Euler constant function $\gamma(x)$ in 2006 which they defined as

(4.5) $$\gamma(x) = \sum_{n=1}^{\infty} x^{n-1} \left[ \frac{1}{n} - \log\left(1 + \frac{1}{n}\right) \right]$$

where we see that $\gamma(1) = \gamma$ and using (4.2) we have $\gamma(-1) = \log(4/\pi)$.

It is easily seen that

$$x\gamma(x) = x(1 - \log 2) + \sum_{n=2}^{\infty} x^n \left[ \frac{1}{n} - \log\left(1 + \frac{1}{n}\right) \right]$$

and using

$$\sum_{n=2}^{\infty} x^n \left[ \frac{1}{n} - \log\left(1 + \frac{1}{n}\right) \right] = \sum_{n=2}^{\infty} x^n \left[ \frac{1}{n} + \sum_{k=1}^{\infty} \frac{(-1)^k}{kn^k} \right]$$

$$= \sum_{n=2}^{\infty} x^n \sum_{k=2}^{\infty} \frac{(-1)^k}{kn^k}$$

we have

(4.6) $$x\gamma(x) = \sum_{n=2}^{\infty} \frac{(-1)^n}{n} Li_n(x)$$

as originally shown by Sondow and Hadjicostas [78].

With reference to (4.6) and, using (4.4.38b) from [24]

$$Li_s(x) = \frac{x(-1)^{s-1}}{\Gamma(s)} \int_0^1 \frac{\log^{s-1} y}{(1 - xy)} dy$$

we have

$$\sum_{n=2}^{\infty} \frac{(-1)^n}{n} Li_n(x) = -x \sum_{n=2}^{\infty} \int_0^1 \frac{\log^n y}{n!(1 - xy)\log y} dy$$

so that

$$(4.7) \qquad \sum_{n=2}^{\infty} \frac{(-1)^n}{n} Li_n(x) = x\int_0^1 \frac{1-y+\log y}{(1-xy)\log y} dy$$

$$= x\int_0^1 \frac{1-y}{(1-xy)\log y} dy - \log(1-x)$$

We therefore see from (4.6) that

$$(4.8) \qquad \gamma(x) = \int_0^1 \frac{1-y+\log y}{(1-xy)\log y} dy$$

which was also derived by Sondow and Hadjicostas [78, equation (11)] by a different method. See also equation (4.4.112b) in [24].

Letting $x = 1$ in (4.7) results in

$$(4.9) \qquad \sum_{n=2}^{\infty} \frac{(-1)^n}{n} \varsigma(n) = \int_0^1 \frac{1-y+\log y}{(1-y)\log y} dy = \gamma$$

and with $x = -1$ we get

$$(4.10) \qquad \log\frac{4}{\pi} = \int_0^1 \frac{1-y+\log y}{(1+y)\log y} dy$$

$$= \int_0^1 \frac{1-y}{(1+y)\log y} dy + \log 2$$

Hence we see that [76]

$$(4.11) \qquad \log\frac{\pi}{2} = \int_0^1 \frac{y-1}{(1+y)\log y} dy$$

Three other proofs of this result are given in Sondow's paper [76]. I subsequently discovered that a derivation of this integral was also recorded in 1864 in Bertrand's treatise [7, Book II, p.149], albeit with a trivial sign error being made in the last step of the proof. Bertrand showed that

(4.11.1) $$\int_0^1 \frac{y^\alpha - y^\beta}{y(1+y^n)\log y} dy = \log\left[\frac{\alpha}{\beta} \cdot \frac{\beta+n}{\alpha+n} \cdot \frac{\alpha+2n}{\beta+2n} \cdots \right]$$

and (4.11) follows by letting $\alpha = 1$, $\beta = 2$ and $n = 1$ and then applying the Wallis identity.

It is reported in Whittaker & Watson [80, p.262] that Kummer discovered the following identity (for $\alpha > 0$, $\beta > 0$)

(4.11.2) $$\int_0^1 \frac{t^{\alpha-1} - t^{\beta-1}}{(1+t)\log t} dt = \log \frac{\Gamma\left(\frac{1+\alpha}{2}\right)\Gamma\left(\frac{\beta}{2}\right)}{\Gamma\left(\frac{1+\beta}{2}\right)\Gamma\left(\frac{\alpha}{2}\right)}$$

(this integral is also contained in G&R [44, p.541]). Therefore, for $\alpha = 2$ and $\beta = 1$ we have

$$\int_0^1 \frac{t-1}{(1+t)\log t} dt = \log \frac{\Gamma(3/2)\Gamma(1/2)}{\Gamma(1)\Gamma(1)} = \log \frac{\pi}{2}$$

and this particular integral is also contained in G&R [44, p.540]. We also showed in equation (3.86h) in [22] that

(4.11.3) $$\sum_{n=0}^\infty \frac{1}{2^{n+1}} \sum_{k=0}^n \binom{n}{k}(-1)^k \log \frac{b+k}{a+k} = \int_0^1 \frac{t^{b-1} - t^{a-1}}{(1+t)\log t} dt$$

From (4.11.1) we see that

$$\int_0^1 \frac{y^{\alpha-1} - y^{\beta-1}}{(1+y)\log y} dy = \log\left[\frac{\alpha}{\beta} \cdot \frac{\beta+1}{\alpha+1} \cdot \frac{\alpha+2}{\beta+2} \cdots \right]$$

and hence we obtain

(4.11.4) $$\frac{\Gamma\left(\frac{1+\alpha}{2}\right)\Gamma\left(\frac{\beta}{2}\right)}{\Gamma\left(\frac{1+\beta}{2}\right)\Gamma\left(\frac{\alpha}{2}\right)} = \frac{\alpha}{\beta} \cdot \frac{\beta+1}{\alpha+1} \cdot \frac{\alpha+2}{\beta+2} \cdots$$

Integrating (4.5) gives us

$$\int_0^u \gamma(x) dx = \sum_{n=1}^\infty \frac{u^n}{n}\left[\frac{1}{n} - \log\left(1+\frac{1}{n}\right)\right]$$

$$= Li_2(u) - \sum_{n=1}^{\infty} \frac{u^n}{n} \log\left(1 + \frac{1}{n}\right)$$

Alternatively, using (4.6) we obtain

$$\int_0^u \gamma(x)\,dx = \sum_{n=2}^{\infty} \frac{(-1)^n}{n} \int_0^u \frac{Li_n(x)}{x}\,dx$$

and thus we get

(4.12) $$\int_0^u \gamma(x)\,dx = \sum_{n=2}^{\infty} \frac{(-1)^n}{n} Li_{n+1}(u)$$

Hence we have

(4.13) $$\sum_{n=1}^{\infty} \frac{u^n}{n} \log\left(1 + \frac{1}{n}\right) = \sum_{n=1}^{\infty} \frac{(-1)^{n+1}}{n} Li_{n+1}(u)$$

where the summation on the right-hand side now starts at $n = 1$.

Letting $u = 1$ results in

(4.14) $$\sum_{n=1}^{\infty} \frac{1}{n} \log\left(1 + \frac{1}{n}\right) = \sum_{n=1}^{\infty} \frac{(-1)^{n+1}}{n} \varsigma(n+1)$$

which concurs with (3.10.2) and (3.10.4).

We may also derive (4.13) as follows

$$\sum_{n=1}^{\infty} \frac{(-1)^n}{n} Li_{n+1}(u) = \sum_{n=1}^{\infty} \frac{(-1)^n}{n} \sum_{k=1}^{\infty} \frac{u^k}{k^{n+1}}$$

$$= \sum_{n=1}^{\infty} \frac{(-1)^n}{nk^n} \sum_{k=1}^{\infty} \frac{u^k}{k}$$

$$= \sum_{k=1}^{\infty} \frac{u^k}{k} \sum_{n=1}^{\infty} \frac{(-1)^n}{nk^n}$$

$$= -\sum_{k=1}^{\infty} \frac{u^k}{k} \log\left(1 + \frac{1}{k}\right)$$

□

We have from equation (4.4.38e) in [24] for $n \geq 1$

$$Li_{n+1}(u) = \frac{(-1)^n}{(n-1)!} \int_0^1 \frac{\log^{n-1} y \log(1-uy)}{y} dy$$

and therefore we get

$$\sum_{n=2}^{\infty} \frac{(-1)^n}{n} Li_{n+1}(u) x^n = \sum_{n=2}^{\infty} \int_0^1 \frac{x^n \log^n y}{n!} \frac{\log(1-uy)}{y \log y} dy$$

$$= \int_0^1 \sum_{n=2}^{\infty} \frac{x^n \log^n y}{n!} \frac{\log(1-uy)}{y \log y} dy$$

$$= \int_0^1 \left[ e^{x \log y} - 1 - x \log y \right] \frac{\log(1-uy)}{y \log y} dy$$

$$= \int_0^1 \left[ y^x - 1 - x \log y \right] \frac{\log(1-uy)}{y \log y} dy$$

$$= \int_0^1 \frac{(y^x - 1) \log(1-uy)}{y \log y} dy - x \int_0^1 \frac{\log(1-uy)}{y} dy$$

$$= \int_0^1 \frac{(y^x - 1) \log(1-uy)}{y \log y} dy + x Li_2(u)$$

Hence we see that

(4.15) $$\sum_{n=1}^{\infty} \frac{(-1)^n}{n} Li_{n+1}(u) x^n = \int_0^1 \frac{(y^x - 1) \log(1-uy)}{y \log y} dy$$

and with $x = 1$ we have

(4.16) $$\sum_{n=1}^{\infty} \frac{(-1)^n}{n} Li_{n+1}(u) = \int_0^1 \frac{(y-1) \log(1-uy)}{y \log y} dy$$

which corresponds with (3.10.1) and (3.10.2) in the case where $u = 1$. It appears that *Mathematica* cannot evaluate this integral.

$\square$

In his book, "The Art of Computer Programming", Knuth [52, p.78] states that if $f(x) = \sum_{n=1}^{\infty} a_n x^n$ is convergent, then we have

(4.17) $$\sum_{n=1}^{\infty} a_n H_n^{(1)} x^n = \int_0^1 \frac{f(x) - f(tx)}{1-t} dt$$

Letting $y = tx$ in the above integral we obtain the symmetrical form

(4.18) $$\sum_{n=1}^{\infty} a_n H_n^{(1)} x^n = \int_0^x \frac{f(x) - f(y)}{x - y} dy$$

With $f(x) = Li_s(x) = \sum_{n=1}^{\infty} \frac{x^n}{n^s}$ we therefore obtain

$$\sum_{n=1}^{\infty} \frac{H_n^{(1)}}{n^s} x^n = \int_0^1 \frac{Li_s(x) - Li_s(tx)}{1-t} dt$$

or the symmetrical form

(4.19) $$\sum_{n=1}^{\infty} \frac{H_n^{(1)}}{n^s} x^n = \int_0^x \frac{Li_s(x) - Li_s(y)}{x - y} dy$$

and with $x = \pm 1$ we have

(4.20) $$\sum_{n=1}^{\infty} \frac{H_n^{(1)}}{n^s} = \int_0^1 \frac{\varsigma(s) - Li_s(y)}{1-y} dy$$

(4.21) $$\sum_{n=1}^{\infty} (-1)^n \frac{H_n^{(1)}}{n^s} = \int_0^1 \frac{(2^{1-s}-1)\varsigma(s) - Li_s(-y)}{1-y} dy$$

We now consider the case where $f(x) = x\gamma(x)$ so that

(4.22) $$f(x) = -\sum_{n=1}^{\infty} \left[ \log\left(1 + \frac{1}{n}\right) - \frac{1}{n} \right] x^n$$

and reference to Knuth's formula (4.17) gives us

(4.23) $$\sum_{n=1}^{\infty} H_n \left( \log \frac{n+1}{n} - \frac{1}{n} \right) x^n = -x \int_0^1 \frac{\gamma(x) - \gamma(tx)}{1-t} dt$$

Referring to (2.1) with $x = 1$ we obtain

(4.24) $$\int_0^1 \frac{\gamma - \gamma(t)}{1-t} dt = \frac{1}{2}\left[\varsigma(2) + \gamma^2 + 2\gamma_1\right]$$

We note from the definition (4.5) that

$$\gamma - \gamma(t) = \sum_{n=1}^{\infty} (1 - t^{n-1})\left[\frac{1}{n} - \log\left(1 + \frac{1}{n}\right)\right]$$

so that

$$\int_0^1 \frac{\gamma - \gamma(t)}{1-t} dt = \sum_{n=1}^{\infty} \int_0^1 \frac{1 - t^{n-1}}{1-t}\left[\frac{1}{n} - \log\left(1 + \frac{1}{n}\right)\right] dt$$

It is well known that (which is a direct consequence of (3.15))

(4.25) $$\int_0^1 \frac{1 - t^{n-1}}{1-t} dt = H_n$$

and we simply recover (4.23).

We now divide (4.23) by $x$ and integrate to obtain

$$\sum_{n=1}^{\infty} \frac{H_n}{n}\left(\log\frac{n+1}{n} - \frac{1}{n}\right) u^n = -\int_0^u \int_0^1 \frac{\gamma(x) - \gamma(tx)}{1-t} dt\, dx$$

and we see that

$$\sum_{n=1}^{\infty} \frac{H_n}{n}\left(\log\frac{n+1}{n} - \frac{1}{n}\right) u^n = \sum_{n=1}^{\infty} \frac{H_n}{n} \log\frac{n+1}{n} u^n - \sum_{n=1}^{\infty} \frac{H_n}{n^2} u^n$$

We showed in Eq. (3.105d) in [22] that

$$\sum_{n=1}^{\infty} \frac{H_n}{n^2} u^n = \frac{1}{2}\log^2(1-u)\log u + \log(1-u)Li_2(1-u) + Li_3(u) - Li_3(1-u) + \varsigma(3)$$

Using (4.12) we have

$$\int_0^u [\gamma(x) - \gamma(tx)] dx = \sum_{n=2}^{\infty} \frac{(-1)^n}{n}\left[Li_{n+1}(u) - \frac{1}{t} Li_{n+1}(tu)\right]$$

Hence we obtain

$$(4.26) \quad \sum_{n=1}^{\infty} \frac{H_n}{n} \log \frac{n+1}{n} u^n = -\sum_{n=2}^{\infty} \frac{(-1)^n}{n} \int_0^1 \left[ \frac{t Li_{n+1}(u) - Li_{n+1}(u)}{t(t-1)} \right] dt$$

$$+ \frac{1}{2} \log^2(1-u) \log u + \log(1-u) Li_2(1-u) + Li_3(u) - Li_3(1-u) + \varsigma(3)$$

and with $u = 1$ this becomes

$$\sum_{n=1}^{\infty} \frac{H_n}{n} \log \frac{n+1}{n} = -\sum_{n=2}^{\infty} \frac{(-1)^n}{n} \int_0^1 \left[ \frac{\varsigma(n+1)t - Li_{n+1}(u)}{t(t-1)} \right] dt + 2\varsigma(3)$$

Coffey [19] showed that

$$(4.27) \quad -2 \int_0^1 \left[ \frac{\varsigma(n+1)t - Li_{n+1}(u)}{t(t-1)} \right] dt = (n+1)\varsigma(n+2) - \sum_{k=1}^{n-1} \varsigma(k+1)\varsigma(n+1-k)$$

so that we have

$$(4.28) \quad \sum_{n=1}^{\infty} \frac{H_n}{n} \log \frac{n+1}{n} = \frac{1}{2} \sum_{n=2}^{\infty} \frac{(-1)^n}{n} \left[ (n+1)\varsigma(n+2) - \sum_{k=1}^{n-1} \varsigma(k+1)\varsigma(n+1-k) \right] + 2\varsigma(3)$$

Substituting (7.3.1) gives us

$$(4.29) \quad \sum_{n=1}^{\infty} \frac{H_n}{n} \log \frac{n+1}{n} = \sum_{n=2}^{\infty} \frac{(-1)^n}{n} \sum_{m=1}^{\infty} \frac{H_m}{(m+1)^{n+1}} + 2\varsigma(3)$$

The Maclaurin expansion of the logarithm gives us

$$\sum_{n=1}^{\infty} \frac{H_n}{n} \log\left(1 + \frac{1}{n}\right) = -\sum_{n=1}^{\infty} \frac{H_n}{n} \sum_{m=1}^{\infty} \frac{(-1)^m}{mn^m}$$

$$= -\sum_{m=1}^{\infty} \frac{(-1)^m}{m} \sum_{n=1}^{\infty} \frac{H_n}{n^{m+1}}$$

$$= -\sum_{m=1}^{\infty} \frac{(-1)^m}{m} \left[ \sum_{n=1}^{\infty} \frac{H_{n-1}}{n^{m+1}} + \frac{1}{n^{m+2}} \right]$$

$$= -\sum_{m=1}^{\infty} \frac{(-1)^m}{m} \sum_{n=1}^{\infty} \frac{H_{n-1}}{n^{m+1}} - \sum_{m=1}^{\infty} \frac{(-1)^m \varsigma(m+2)}{m}$$

$$= -\sum_{m=1}^{\infty}\frac{(-1)^m}{m}\sum_{n=1}^{\infty}\frac{H_k}{(k+1)^{m+1}} - \sum_{m=1}^{\infty}\frac{(-1)^m \varsigma(m+2)}{m}$$

$$= \sum_{n=1}^{\infty}\frac{H_k}{(k+1)^2} - \sum_{m=2}^{\infty}\frac{(-1)^m}{m}\sum_{n=1}^{\infty}\frac{H_k}{(k+1)^{m+1}} - \sum_{m=1}^{\infty}\frac{(-1)^m \varsigma(m+2)}{m}$$

and hence we have

(4.30) $$\sum_{n=1}^{\infty}\frac{H_n}{n}\log\frac{n+1}{n} = \varsigma(3) - \sum_{m=2}^{\infty}\frac{(-1)^m}{m}\sum_{n=1}^{\infty}\frac{H_k}{(k+1)^{m+1}} - \sum_{m=1}^{\infty}\frac{(-1)^m \varsigma(m+2)}{m}$$

Equating this with (4.29) gives us

$$\sum_{m=1}^{\infty}\frac{(-1)^m \varsigma(m+2)}{m} = -\varsigma(3) - 2\sum_{m=2}^{\infty}\frac{(-1)^m}{m}\sum_{n=1}^{\infty}\frac{H_k}{(k+1)^{m+1}}$$

We see that

$$\sum_{m=1}^{\infty}\frac{(-1)^m \varsigma(m+2)}{m} = \sum_{j=3}^{\infty}\frac{(-1)^j \varsigma(j)}{j-2}$$

and we note [79, p.160] that

$$\sum_{j=2}^{\infty}(-1)^j \varsigma(j) t^{j-1} = \psi(1+t) + \gamma \qquad |t| < 1$$

or equivalently

$$\sum_{j=3}^{\infty}(-1)^j \varsigma(j) t^{j-1} = \psi(1+t) + \gamma - \varsigma(2)t$$

This gives us

$$\sum_{j=3}^{\infty}(-1)^j \varsigma(j) t^{j-3} = \frac{\psi(1+t) + \gamma - \varsigma(2)t}{t^2}$$

Integration results in

$$\sum_{j=3}^{\infty}\frac{(-1)^j \varsigma(j) u^{j-2}}{j-2} = \int_0^u \frac{\psi(1+t) + \gamma - \varsigma(2)t}{t^2}\,dt$$

so that

$$\sum_{j=3}^{\infty}\frac{(-1)^{j}\varsigma(j)}{j-2}=\int_{0}^{1}\frac{\psi(1+t)+\gamma-\varsigma(2)t}{t^{2}}dt$$

We therefore have

(4.31) $$\int_{0}^{1}\frac{\psi(1+t)+\gamma-\varsigma(2)t}{t^{2}}dt=-\varsigma(3)-2\sum_{m=2}^{\infty}\frac{(-1)^{m}}{m}\sum_{n=1}^{\infty}\frac{H_{k}}{(k+1)^{m+1}}$$

□

With $u \to 1-u$ in (4.13) we get

(4.32) $$\sum_{n=1}^{\infty}\frac{(1-u)^{n-1}}{n}\log\left(1+\frac{1}{n}\right)=-\sum_{r=1}^{\infty}\frac{(-1)^{r}}{r}\frac{Li_{r+1}(1-u)}{1-u}$$

We also have the known integral expression for the harmonic numbers $H_n$ (see for example [25])

$$H_{n}=-n\int_{0}^{1}(1-u)^{n-1}\log u\, du$$

so that multiplying (4.32) by $\log u$ followed by integration gives us

(4.33) $$\sum_{n=1}^{\infty}\frac{H_{n}}{n^{2}}\log\left(1+\frac{1}{n}\right)=\sum_{r=1}^{\infty}\frac{(-1)^{r}}{r}\int_{0}^{1}\frac{Li_{r+1}(1-u)\log u}{1-u}du$$

It is easily shown that

$$\int\frac{Li_{2}(1-u)\log u}{1-u}du=\frac{1}{2}[Li_{2}(1-u)]^{2}$$

but it appears that *Mathematica* cannot evaluate the higher order integrals.

However, we showed in Equation (4.4.24y) in [24] that

(4.34) $$\int_{0}^{1}\frac{u^{\beta-1}Li_{p}[(1-u)t]\log u}{1-u}du=-\sum_{n=0}^{\infty}\frac{t^{n}}{(n+1)^{p}}\sum_{k=0}^{n}\binom{n}{k}\frac{(-1)^{k}}{(k+\beta)^{2}}$$

and with $\beta=t=1$ we have

(4.35) $$\int_0^1 \frac{Li_{r+1}(1-u)\log u}{1-u}\,du = -\sum_{n=0}^{\infty}\frac{1}{(n+1)^{r+1}}\sum_{k=0}^{n}\binom{n}{k}\frac{(-1)^k}{(k+1)^2}$$

Larcombe et al. [56] showed that

$$m\binom{m+n}{n}\sum_{k=0}^{n}\binom{n}{k}\frac{(-1)^k}{(m+k)^2} = \sum_{k=m}^{m+n}\frac{1}{k}$$

so that letting $m=1$ results in

$$\sum_{k=0}^{n}\binom{n}{k}\frac{(-1)^k}{(k+1)^2} = \frac{H_{n+1}}{n+1}$$

Hence we obtain

$$\int_0^1 \frac{Li_{r+1}(1-u)\log u}{1-u}\,du = -\sum_{n=0}^{\infty}\frac{H_{n+1}}{(n+1)^{r+2}} = -\sum_{n=1}^{\infty}\frac{H_n}{n^{r+2}}$$

and thus we deduce that

(4.36) $$\sum_{n=1}^{\infty}\frac{H_n}{n^2}\log\left(1+\frac{1}{n}\right) = \sum_{r=1}^{\infty}\frac{(-1)^{r+1}}{r}\sum_{n=1}^{\infty}\frac{H_n}{n^{r+2}}$$

We shall see some applications of Lerch's theorem with the function $\gamma(x)$ in the next section.

## 5. Various applications of Lerch's theorem

Lerch [58] proved the following theorem in 1895:

If the series

(5.1) $$f(x) = \sum_{n=1}^{\infty}\frac{c_n}{n}\sin 2\pi nx$$

is convergent for $0 < x < 1$, then the derivative of $f(x)$ is given in this interval by

(5.2) $$f'(x)\cdot\frac{\sin \pi x}{\pi} = \sum_{n=0}^{\infty}(c_n - c_{n+1})\sin(2n+1)\pi x$$

where $c_0 = 0$, provided that the last series converges uniformly for $\varepsilon \leq x \leq 1-\varepsilon$ for all $\varepsilon > 0$. Equation (5.2) may be written as

(5.3) $$f'(x) \cdot \frac{\sin \pi x}{\pi} = -c_1 \sin \pi x + \sum_{n=1}^{\infty}(c_n - c_{n+1})\sin(2n+1)\pi x$$

Subject to the same conditions, if the series

(5.4) $$g(x) = \sum_{n=1}^{\infty} \frac{c_n}{n} \cos 2\pi n x$$

is convergent for $0 < x < 1$, then we have

(5.5) $$g'(x) \cdot \frac{\sin \pi x}{\pi} = \sum_{n=0}^{\infty}(c_n - c_{n+1})\cos(2n+1)\pi x$$

or equivalently

(5.6) $$g'(x) \cdot \frac{\sin \pi x}{\pi} = -c_1 \cos \pi x + \sum_{n=1}^{\infty}(c_n - c_{n+1})\cos(2n+1)\pi x$$

□

Deninger [34] showed that

(5.7) $$\sum_{n=1}^{\infty} \frac{\log n}{n} \cos 2n\pi x = \frac{1}{2}[\varsigma''(0,x) + \varsigma''(0,1-x)] + [\gamma + \log(2\pi)]\log(2\sin \pi x)$$

and another proof appears in [29]. This is a companion to Kummer's Fourier series formula for $\log \Gamma(x)$

(5.8) $$\log \Gamma(x) = \frac{1}{2}\log \frac{\pi}{\sin \pi x} + \frac{1}{2}(1-2x)[\gamma + \log(2\pi)] + \frac{1}{\pi}\sum_{n=1}^{\infty} \frac{\log n}{n} \sin 2n\pi x$$

Using Euler's reflection formula $\Gamma(x)\Gamma(1-x) = \frac{\pi}{\sin \pi x}$ we may write this as

$$\sum_{n=1}^{\infty} \frac{\log n}{n} \sin 2n\pi x = \frac{\pi}{2}[\log \Gamma(x) - \log \Gamma(1-x)] - [\gamma + \log(2\pi)]\frac{\pi}{2}(1-2x)$$

With Lerch's identity [79, p.92]

$$\log \Gamma(x) = \varsigma'(0,x) + \frac{1}{2}\log(2\pi)$$

this becomes

$$\sum_{n=1}^{\infty}\frac{\log n}{n}\sin 2n\pi x = \frac{\pi}{2}[\varsigma'(0,x)-\varsigma'(0,1-x)]-[\gamma+\log(2\pi)]\frac{\pi}{2}(1-2x)$$

and, in this format, we can see the structural similarity with (5.7). Using (5.11) and (5.12) we may write

$$\sum_{n=1}^{\infty}\frac{\gamma+\log(2\pi n)}{n}\cos 2n\pi x = \frac{1}{2}[\varsigma''(0,x)+\varsigma''(0,1-x)]$$

$$\sum_{n=1}^{\infty}\frac{\gamma+\log(2\pi n)}{n}\sin 2n\pi x = \frac{\pi}{2}[\varsigma'(0,x)-\varsigma'(0,1-x)]$$

We have
$$\varsigma''(s,x) = \sum_{n=0}^{\infty}\frac{\log^2(n+x)}{(n+x)^s}$$

and
$$\frac{\partial}{\partial x}\varsigma''(s,x) = 2\sum_{n=0}^{\infty}\frac{\log(n+x)}{(n+x)^{s+1}} - s\sum_{n=0}^{\infty}\frac{\log^2(n+x)}{(n+x)^{s+1}}$$

$$\frac{\partial}{\partial x}\varsigma''(s,1-x) = -2\sum_{n=0}^{\infty}\frac{\log(n+1-x)}{(n+1-x)^{s+1}} + s\sum_{n=0}^{\infty}\frac{\log^2(n+1-x)}{(n+1-x)^{s+1}}$$

These are valid in the first instance for $\operatorname{Re}(s)>1$ but, by analytic continuation, the following equation is valid for $s=0$

(5.8.1) $$\frac{\partial}{\partial x}[\varsigma''(0,x)+\varsigma''(0,1-x)] = 2\left[\sum_{n=0}^{\infty}\frac{\log(n+x)}{n+x}-\sum_{n=0}^{\infty}\frac{\log(n+1-x)}{n+1-x}\right]$$

$$= -2[\varsigma'(1,x)-\varsigma'(1,1-x)]$$

Therefore, using Lerch's theorem to differentiate (5.7) results in

(5.9) $$\sum_{n=1}^{\infty}\log\left(1+\frac{1}{n}\right)\cos(2n+1)\pi x = [\varsigma'(1,x)-\varsigma'(1,1-x)]\frac{\sin \pi x}{\pi} - [\gamma+\log(2\pi)]\cos \pi x$$

or equivalently

$$\cos \pi x\sum_{n=1}^{\infty}\log\left(1+\frac{1}{n}\right)\cos 2n\pi x - \sin \pi x\sum_{n=1}^{\infty}\log\left(1+\frac{1}{n}\right)\sin 2n\pi x$$

$$= [\varsigma'(1,x) - \varsigma'(1,1-x)]\frac{\sin \pi x}{\pi} - [\gamma + \log(2\pi)]\cos \pi x$$

We recall (4.4)

(5.10) $$\sum_{n=1}^{\infty} \log\left(1+\frac{1}{n}\right)\sin(2n+1)\pi x = -\left[\psi(x)\sin \pi x + \frac{\pi}{2}\cos \pi x + (\gamma + \log 2\pi)\sin \pi x\right]$$

or equivalently

$$\sin \pi x \sum_{n=1}^{\infty} \log\left(1+\frac{1}{n}\right)\cos 2n\pi x + \cos \pi x \sum_{n=1}^{\infty} \log\left(1+\frac{1}{n}\right)\sin 2n\pi x$$

$$= -\left[\psi(x)\sin \pi x + \frac{\pi}{2}\cos \pi x + [\gamma + \log(2\pi)]\sin \pi x\right]$$

We may combine (5.9) and (5.10) as follows

$$\sum_{n=1}^{\infty} \log\left(1+\frac{1}{n}\right) e^{(2n+1)\pi i x} = [\varsigma'(1,x) - \varsigma'(1,1-x)]\frac{\sin \pi x}{\pi} - [\gamma + \log(2\pi)]\cos \pi x$$

$$-i\left[\psi(x)\sin \pi x + \frac{\pi}{2}\cos \pi x + [\gamma + \log(2\pi)]\sin \pi x\right]$$

or as

$$e^{i\pi x}\sum_{n=1}^{\infty} \log\left(1+\frac{1}{n}\right) e^{2n\pi i x} = [\varsigma'(1,x) - \varsigma'(1,1-x)]\frac{\sin \pi x}{\pi} - [\gamma + \log(2\pi)]\cos \pi x$$

$$-i\left[\psi(x)\sin \pi x + \frac{\pi}{2}\cos \pi x + [\gamma + \log(2\pi)]\sin \pi x\right]$$

We recall the Fourier series shown in Carslaw's book [12, p.241]

(5.11) $$-\log(2\sin \pi x) = \sum_{n=1}^{\infty} \frac{\cos 2n\pi x}{n}$$

(5.12) $$\frac{1}{2}\pi(1-2x) = \sum_{n=1}^{\infty} \frac{\sin 2n\pi x}{n}$$

Indeed, Lerch [58] showed how the above two series may be directly derived from his theorem by letting $c_n \equiv 1$ in (5.3) and (5.6).

Letting $x \to e^{i2\pi x}$ in (4.5) we then see that

$$e^{i2\pi x}\gamma(e^{i2\pi x}) = \sum_{n=1}^{\infty}\left[\frac{1}{n} - \log\left(1 + \frac{1}{n}\right)\right]e^{i2n\pi x}$$

$$= -\log(2\sin\pi x) + i\frac{1}{2}\pi(1-2x)$$

$$-\left[[\varsigma'(1,x) - \varsigma'(1,1-x)]\frac{\sin\pi x}{\pi} - [\gamma + \log(2\pi)]\cos\pi x\right]e^{-i\pi x}$$

$$+ i\left[\psi(x)\sin\pi x + \frac{\pi}{2}\cos\pi x + [\gamma + \log(2\pi)]\sin\pi x\right]e^{-i\pi x}$$

We recall the Laurent expansion (1.1) of the Hurwitz zeta function

$$\varsigma(s,u) = \frac{1}{s-1} + \sum_{p=0}^{\infty}\frac{(-1)^p}{p!}\gamma_p(u)(s-1)^p$$

and we have

$$\varsigma(s,u) - \varsigma(s,1-u) = \sum_{p=0}^{\infty}\frac{(-1)^p}{p!}[\gamma_p(u) - \gamma_p(1-u)](s-1)^p$$

Differentiation with respect to $s$ results in

$$\varsigma'(s,u) - \varsigma'(s,1-u) = \sum_{p=0}^{\infty}\frac{(-1)^p}{p!}p[\gamma_p(u) - \gamma_p(1-u)](s-1)^{p-1}$$

and in the limit as $s \to 1$ we have

(5.13)  $\varsigma'(1,u) - \varsigma'(1,1-u) = -[\gamma_1(u) - \gamma_1(1-u)]$

In fact, as previously noted in equation (4.3.228b) in [24], more generally we have

(5.14)  $\gamma_p(x) - \gamma_p(y) = \lim_{s\to 1}(-1)^p \frac{\partial^p}{\partial s^p}[\varsigma(s,x) - \varsigma(s,y)]$

Therefore using (5.13) we obtain

$$e^{i2\pi x}\gamma(e^{i2\pi x}) = -\log(2\sin\pi x) + i\frac{1}{2}\pi(1-2x)$$

$$+\left[[\gamma_1(x)-\gamma_1(1-x)]\frac{\sin \pi x}{\pi}+[\gamma+\log(2\pi)]\cos \pi x\right]e^{-i\pi x}$$

$$+i\left[\psi(x)\sin \pi x+\frac{\pi}{2}\cos \pi x+[\gamma+\log(2\pi)]\sin \pi x\right]e^{-i\pi x}$$

or equivalently

(5.15) $$e^{i3\pi x}\gamma(e^{i2\pi x})=-e^{i\pi x}\log(2\sin \pi x)+ie^{i\pi x}\frac{1}{2}\pi(1-2x)$$

$$+\left[[\gamma_1(x)-\gamma_1(1-x)]\frac{\sin \pi x}{\pi}+[\gamma+\log(2\pi)]\cos \pi x\right]$$

$$+i\left[\psi(x)\sin \pi x+\frac{\pi}{2}\cos \pi x+[\gamma+\log(2\pi)]\sin \pi x\right]$$

As a spot check on the algebra, we note that letting $x=1/2$ results in $\gamma(-1)=\log(4/\pi)$.

$$e^{i\pi}\gamma(e^{i\pi})=-\gamma(-1)=-\log 2-[\psi(1/2)+[\gamma+\log(2\pi)]]=\log(4/\pi)$$

Adamchik [2] has reported that

(5.16) $$\varsigma'\left(1,\frac{p}{q}\right)-\varsigma'\left(1,1-\frac{p}{q}\right)=\pi[\log(2\pi q)+\gamma]\cot\left(\frac{\pi p}{q}\right)-2\pi\sum_{j=1}^{q-1}\log\Gamma\left(\frac{j}{q}\right)\sin\left(\frac{2\pi jp}{q}\right)$$

where $p$ and $q$ are positive integers and $p < q$. Adamchik [2] notes that this formula was first proved by Almkvist and Meurman. Adamchik's derivation [2] is rather terse and an expanded exposition may be found in [29]

We therefore obtain from (5.14) and (5.16)

(5.17) $$-\left[\gamma_1\left(\frac{p}{q}\right)-\gamma_1\left(1-\frac{p}{q}\right)\right]=\pi[\log(2\pi q)+\gamma]\cot\left(\frac{\pi p}{q}\right)-2\pi\sum_{j=1}^{q-1}\log\Gamma\left(\frac{j}{q}\right)\sin\left(\frac{2\pi jp}{q}\right)$$

With $x=1/3$ in (5.15) and using

$$\gamma_1\left(\tfrac{1}{3}\right)-\gamma_1\left(\tfrac{2}{3}\right)=-\frac{\pi}{2\sqrt{3}}\left[2\gamma-\log 3+8\log(2\pi)-12\log\Gamma\left(\tfrac{1}{3}\right)\right]$$

$$\psi\left(\tfrac{1}{3}\right) = -\left(\gamma + \frac{1}{6}\pi\sqrt{3} + \frac{3}{2}\log 3\right)$$

we obtain

$$(5.18) \quad \gamma(e^{i2\pi/3}) = \frac{\sqrt{3}}{12}\pi - 3\log\Gamma\left(\tfrac{1}{3}\right) + \frac{3}{2}\log(2\pi) + i\left[\sqrt{3}\log 3 - \frac{\sqrt{3}}{2}\log(2\pi) - \frac{1}{12}\pi\right]$$

which concurs with Example 8 in [78].

It should be noted that Sondow and Hadjicostas [78] have also derived several formulae for $\gamma(e^{i\pi p/q})$ by very different methods. These formulae are significantly more compact than the one obtained above.

$\square$

Letting $u = e^{2ix}$ in (4.13) gives us

$$(5.19) \quad \sum_{n=1}^{\infty} \frac{\sin 2n\pi x}{n}\log\left(1+\frac{1}{n}\right) = \sum_{n=1}^{\infty} \frac{(-1)^{n+1}}{n} Li_{n+1}(\sin 2\pi x)$$

$$(5.20) \quad \sum_{n=1}^{\infty} \frac{\cos 2n\pi x}{n}\log\left(1+\frac{1}{n}\right) = \sum_{n=1}^{\infty} \frac{(-1)^{n+1}}{n} Li_{n+1}(\cos 2\pi x)$$

and differentiation via Lerch's theorem results in

$$(5.21) \quad 2\sin\pi x \cot 2\pi x \sum_{n=1}^{\infty} \frac{(-1)^{n+1}}{n} Li_n(\sin 2\pi x)$$

$$= -\sin\pi x \log 2 + \sum_{n=1}^{\infty} \log \frac{(n+1)^2}{n(n+2)} \sin(2n+1)\pi x$$

$$(5.22) \quad -2\sin\pi x \tan 2\pi x \sum_{n=1}^{\infty} \frac{(-1)^{n+1}}{n} Li_n(\cos 2\pi x)$$

$$= -\cos\pi x \log 2 + \sum_{n=1}^{\infty} \log \frac{(n+1)^2}{n(n+2)} \cos(2n+1)\pi x$$

With $x = 0$ in (5.22) we obtain a particular case of (3.20)

$$(5.23) \quad \sum_{n=1}^{\infty} \log \frac{(n+1)^2}{n(n+2)} = \log 2$$

or equivalently

$$\sum_{n=1}^{\infty} \log\left(1 + \frac{1}{n(n+2)}\right) = \log 2$$

Dividing (21) by $\sin \pi x$ gives us

$$2\cot 2\pi x \sum_{n=1}^{\infty} \frac{(-1)^{n+1}}{n} Li_n(\sin 2\pi x) = -\log 2 + \sum_{n=1}^{\infty} \log \frac{(n+1)^2}{n(n+2)} \frac{\sin(2n+1)\pi x}{\sin \pi x}$$

and using [44, p.35, 1.342.2]

$$\frac{\sin(2n+1)\pi x}{\sin \pi x} = 1 + 2\sum_{k=1}^{n} \cos 2k\pi x$$

results in

(5.24) $$\cot 2\pi x \sum_{n=1}^{\infty} \frac{(-1)^{n+1}}{n} Li_n(\sin 2\pi x) = \sum_{n=1}^{\infty} \log \frac{(n+1)^2}{n(n+2)} \sum_{k=1}^{n} \cos 2k\pi x$$

where we have employed (5.23).

One may also apply Lerch's theorem to

$$\frac{1}{2}\log^2(1-t) + Li_2(t) = \sum_{n=1}^{\infty} \frac{H_n}{n} t^n$$

in the form

$$\frac{1}{2}\log^2(1-e^{ix}) + Li_2(e^{ix}) = \sum_{n=1}^{\infty} \frac{H_n}{n}(\cos nx + i\sin nx)$$

and also to the series

$$-\frac{\log(1-t)}{1-t} = \sum_{n=1}^{\infty} H_n t^n$$

## 6. The Barnes double gamma function

The following was posed as a Quickie Problem Q974 in the October 2007 issue of Mathematics Magazine: prove that

$$\sum_{k=2}^{\infty}(-1)^k \frac{\varsigma(k)}{k+1} = 1 + \frac{1}{2}\gamma - \frac{1}{2}\log(2\pi)$$

A different derivation is shown below.

We have Euler's formula for the gamma function [79, p.2]

$$\Gamma(x) = \frac{1}{x}\prod_{n=1}^{\infty}\left[\left(1+\frac{1}{n}\right)^x\left(1+\frac{x}{n}\right)^{-1}\right]$$

and multiplying this by $x$ and taking the logarithm of both sides gives us

$$\log\Gamma(1+x) = \sum_{n=1}^{\infty}\left[x\log\left(1+\frac{1}{n}\right) - \log\left(1+\frac{x}{n}\right)\right]$$

This may be written as

$$\log\Gamma(1+x) = \sum_{n=1}^{\infty}\left[x\log\left(1+\frac{1}{n}\right) + \sum_{k=1}^{\infty}\frac{(-1)^k x^k}{kn^k}\right]$$

$$= \sum_{n=1}^{\infty}\left[x\left[\log\left(1+\frac{1}{n}\right) - \frac{1}{n}\right] + \sum_{k=2}^{\infty}\frac{(-1)^k x^k}{kn^k}\right]$$

$$= x\sum_{n=1}^{\infty}\left[\log\left(1+\frac{1}{n}\right) - \frac{1}{n}\right] + \sum_{n=1}^{\infty}\sum_{k=2}^{\infty}\frac{(-1)^k x^k}{kn^k}$$

$$= x\sum_{n=1}^{\infty}\left[\log\left(1+\frac{1}{n}\right) - \frac{1}{n}\right] + \sum_{k=2}^{\infty}\frac{(-1)^k x^k}{k}\sum_{n=1}^{\infty}\frac{1}{n^k}$$

and hence, in a rather convoluted manner, we have derived the Maclaurin series for $\log\Gamma(1+x)$

(6.1) $$\log\Gamma(1+x) = -\gamma x + \sum_{k=2}^{\infty}\frac{(-1)^k \varsigma(k)}{k}x^k$$

This series may be obtained much more directly by noting that

$$\frac{d^{n+1}}{dx^{n+1}}\log\Gamma(1+x) = \psi^{(n)}(1+x)$$

Differentiation results in

(6.2) $$\psi(x+1) = -\gamma + \sum_{k=2}^{\infty}(-1)^k \varsigma(k) x^{k-1}$$

and we now multiply this equation by $x$ and integrate to obtain

$$\int_0^u x\psi(x+1)\,dx = -\frac{1}{2}\gamma u^2 + \sum_{k=2}^{\infty}(-1)^k \frac{\varsigma(k)}{k+1} u^{k+1}$$

Integration by parts gives us

$$\int_0^u x\psi(x+1)\,dx = x\log\Gamma(x+1)\Big|_0^u - \int_0^u \log\Gamma(x+1)\,dx$$

$$= u\log\Gamma(u+1) - \int_0^u \log\Gamma(x+1)\,dx$$

Hence we obtain

$$\int_0^u \log\Gamma(x+1)\,dx = u\log\Gamma(u+1) + \frac{1}{2}\gamma u^2 - \sum_{k=2}^{\infty}(-1)^k \frac{\varsigma(k)}{k+1} u^{k+1}$$

or alternatively

$$\int_0^u \log\Gamma(x)\,dx = u\log\Gamma(u) + u + \frac{1}{2}\gamma u^2 - \sum_{k=2}^{\infty}(-1)^k \frac{\varsigma(k)}{k+1} u^{k+1}$$

We have the well-known Raabe's integral

$$\int_0^1 \log\Gamma(x)\,dx = \frac{1}{2}\log(2\pi)$$

and hence with $u=1$ we have

$$\frac{1}{2}\log(2\pi) = 1 + \frac{1}{2}\gamma - \sum_{k=2}^{\infty}(-1)^k \frac{\varsigma(k)}{k+1}$$

The result may be generalised by reference to Alexeiewsky's theorem: in 1894 Alexeiewsky [79, p.32] showed that

$$\int_0^u \log\Gamma(x+1)\,dx = \frac{1}{2}u(u-1) + \frac{1}{2}\log(2\pi) - \log G(u+1) + u\log\Gamma(u)$$

where $G(x)$ is the Barnes double gamma function defined by [79, p.25]

$$G(1+x) = (2\pi)^{x/2} \exp\left[-\frac{1}{2}(\gamma x^2 + x^2 + x)\right] \prod_{k=1}^{\infty}\left\{\left(1+\frac{x}{k}\right)^k \exp\left(\frac{x^2}{2k}-x\right)\right\}$$

It is obvious from the definition that $G(1)=1$ and, since [79, p.25] $G(1+x)=G(x)\Gamma(x)$, we have $G(2)=1$.

We then obtain the following identity originally derived by Srivastava [79, p.210] in 1988

(6.3) $$\sum_{k=2}^{\infty}(-1)^k \frac{\varsigma(k)}{k+1}u^{k+1} = [1-\log(2\pi)]\frac{u}{2} + (1+\gamma)\frac{u^2}{2} + \log G(u+1)$$

The alert reader will wonder why this section on the Barnes double gamma function appears in this paper. The reason is rather simple: in a temporary mental aberration I erroneously confused the series $\sum_{k=2}^{\infty}(-1)^k \frac{\varsigma(k)}{k+1}$ with the more fiendishly difficult one $\sum_{k=1}^{\infty}(-1)^k \frac{\varsigma(k+1)}{k}$ !

Differentiating (6.3) gives us

$$\sum_{k=2}^{\infty}(-1)^k \varsigma(k)u^k = [1-\log(2\pi)]\frac{1}{2} + (1+\gamma)u + \frac{d}{du}\log G(u+1)$$

and dividing by $u^2$ followed by integration results in

$$\sum_{k=2}^{\infty}(-1)^k \frac{\varsigma(k)}{k-1}t^{k-1} = \int_0^t \left[\frac{1}{2}[1-\log(2\pi)] + (1+\gamma)u + \frac{d}{du}\log G(u+1)\right]\frac{du}{u^2}$$

It is reported in [80, p.264] that

$$\frac{G'(1+u)}{G(1+u)} = \frac{1}{2}\log(2\pi) + \frac{1}{2} - u + u\psi(u)$$

and we therefore obtain the equivalence of (3.8) and (3.10.2).

□

We also have the well-known result [79, p.159] for $|x| < |u|$

$$\sum_{n=2}^{\infty}(-1)^n \varsigma(n,u)\frac{x^n}{n} = \log\Gamma(u+x) - \log\Gamma(u) - x\psi(u)$$

and differentiating this results in

$$\sum_{n=2}^{\infty}(-1)^n \varsigma(n,u)x^{n-1} = \psi(u+x) - \psi(u)$$

We now divide this by $x$ and integrate to obtain

$$\sum_{n=2}^{\infty}(-1)^n \frac{\varsigma(n,u)}{n-1}t^{n-1} = \int_0^t \frac{\psi(u+x) - \psi(u)}{x}dx$$

or equivalently

(6.4) $$\sum_{k=1}^{\infty}(-1)^{k+1}\frac{\varsigma(k+1,u)}{k}t^k = \int_0^t \frac{\psi(u+x) - \psi(u)}{x}dx$$

With $t = u = 1$ we obtain the equivalence of (3.8) and (3.10.2) again.

□

We now recall the series

(6.5) $$\psi(x+u) - \psi(u) = \sum_{k=1}^{\infty}\frac{(-1)^{k+1}}{k}\frac{x(x-1)...(x-k+1)}{u(u+1)...(u+k-1)}$$

which converges for $\mathrm{Re}(x+u) > 0$. According to Raina and Ladda [69], this summation formula is due to Nörlund [68, p.168].

Combining (6.4) with (6.5) gives us

$$\sum_{k=1}^{\infty}(-1)^{k+1}\frac{\varsigma(k+1,u)}{k}t^k = \sum_{k=1}^{\infty}\frac{(-1)^{k+1}}{k}\frac{1}{u(u+1)...(u+k-1)}\int_0^t(x-1)...(x-k+1)dx$$

The Stirling numbers $s(n,k)$ of the first kind [79, p.56] may be defined by the following generating function

(6.6) $$x(x-1)...(x-k+1) = \sum_{j=0}^{k} s(k,j)x^j$$

or by the Maclaurin expansion for $|x| < 1$

(6.7) $$\log^k(1+x) = k!\sum_{n=k}^{\infty} s(n,k)\frac{x^n}{n!}$$

We see from (6.6) that

$$\int_0^t (x-1)\ldots(x-k+1)\,dx = \sum_{j=1}^{k} \frac{s(k,j)}{j} t^j$$

where the summation starts at $j=1$ because $s(k,0)=0$ for all $k \in \mathbf{N}$. We then have

(6.8) $$\sum_{k=1}^{\infty} (-1)^{k+1} \frac{\varsigma(k+1,u)}{k} t^k = \sum_{k=1}^{\infty} \frac{(-1)^{k+1}}{k} \frac{1}{u(u+1)\ldots(u+k-1)} \sum_{j=1}^{k} \frac{s(k,j)}{j} t^j$$

Letting $u=1$ gives us

(6.9) $$\sum_{k=1}^{\infty} (-1)^{k+1} \frac{\varsigma(k+1)}{k} = \sum_{k=1}^{\infty} \frac{(-1)^{k+1}}{k.k!} \sum_{j=1}^{k} \frac{s(k,j)}{j}$$

where the left-hand side corresponds with (3.10.2).

Differentiation of (6.8) with respect to $t$ results in

$$\sum_{k=1}^{\infty} (-1)^{k+1} \varsigma(k+1,u) t^{k-1} = \sum_{k=1}^{\infty} \frac{(-1)^{k+1}}{k} \frac{1}{u(u+1)\ldots(u+k-1)} \sum_{j=1}^{k} s(k,j) t^{j-1}$$

and with $t=0$ we obtain

$$\varsigma(2,u) = \sum_{k=1}^{\infty} \frac{(-1)^{k+1}}{k} \frac{s(k,1)}{u(u+1)\ldots(u+k-1)}$$

Since [79, p.56] $s(k,1) = (-1)^{k+1}(k-1)!$ we obtain

$$\varsigma(2,u) = \sum_{k=1}^{\infty} \frac{1}{k} \frac{(k-1)!}{u(u+1)\ldots(u+k-1)}$$

and with $u=1$ this simply reverts to the Riemann zeta function

$$\varsigma(2,1) = \varsigma(2) = \sum_{k=1}^{\infty} \frac{1}{k^2}$$

Differentiation of (6.8) $p$ times with respect to $t$ results in

$$\sum_{k=1}^{\infty}(-1)^{k+1}\varsigma(k+1,u)(k-1)\cdots(k-p+1)t^{k-p}$$

$$=\sum_{k=1}^{\infty}\frac{(-1)^{k+1}}{k}\frac{1}{u(u+1)\ldots(u+k-1)}\sum_{j=1}^{k}s(k,j)(k-1)\cdots(k-p+1)t^{j-p}$$

and with $t=0$ we obtain

(6.10) $\qquad \varsigma(p+1,u)=(-1)^p\sum_{k=1}^{\infty}\frac{(-1)^k}{k}\frac{s(k,p)}{u(u+1)\ldots(u+k-1)}$

With $u=1$ we get

(6.11) $\qquad \varsigma(p+1)=(-1)^p\sum_{k=1}^{\infty}\frac{(-1)^k}{k.k!}s(k,p)$

This result was previously obtained in 1995 by Shen [74] by employing a different method.

We may also differentiate (6.10) with respect to $u$. We designate $A(u)$ as

$$A(u)=\frac{1}{u(u+1)\ldots(u+k-1)}=\frac{\Gamma(u)}{\Gamma(u+k)}$$

and we see that the derivative is

(6.12) $\qquad A'(u)=A(u)[\psi(u)-\psi(u+k)]$

Hence we have

(6.13) $\qquad (p+1)\varsigma(p+2,u)=(-1)^{p+1}\sum_{k=1}^{\infty}\frac{(-1)^k}{k}\frac{s(k,p)[\psi(u)-\psi(u+k)]}{u(u+1)\ldots(u+k-1)}$

which may be compared with (6.10) with $p\to p+1$

$$\varsigma(p+2,u)=(-1)^{p+1}\sum_{k=1}^{\infty}\frac{(-1)^k}{k}\frac{s(k,p+1)}{u(u+1)\ldots(u+k-1)}$$

In view of the relation (6.12), the $p$ th derivative will be given in terms of the exponential Bell polynomials [30]

It appears that differentiating (6.8) with respect to $u$ is not valid.

We see from (6.7) that

$$\log^k(1-x) = k! \sum_{n=1}^{\infty} \frac{(-1)^n}{n!} s(n,k) x^n$$

and we obtain by integration

$$\int_0^t \frac{\log^k(1-x)}{x} dx = k! \sum_{n=1}^{\infty} \frac{(-1)^n}{n \cdot n!} s(n,k) t^n$$

In particular we have

$$\int_0^1 \frac{\log^k(1-x)}{x} dx = k! \sum_{n=1}^{\infty} \frac{(-1)^n}{n \cdot n!} s(n,k)$$

and using (3.19) we see that

$$\varsigma(k+1) = (-1)^k \sum_{n=1}^{\infty} \frac{(-1)^n}{n \cdot n!} s(n,k)$$

which corresponds with (6.11).

□

We note the known relationship for Stirling numbers of the first kind for $r \geq 0$

(6.14) $$s(n, r+1) = (-1)^{n+r+1} \frac{(n-1)!}{r!} Y_r\left(H_{n-1}^{(1)}, -1! H_{n-1}^{(2)}, \ldots, (-1)^{r-1}(r-1)! H_{n-1}^{(r)}\right)$$

where $Y_r(x_1, \ldots, x_r)$ are the (exponential) complete Bell polynomials.

The above relationship was derived by Kölbig [54] and a slightly more direct proof of this important formula is shown in [30].

Substituting (6.14) in (6.9) results in

(6.15) $$\sum_{k=1}^{\infty} (-1)^{k+1} \frac{\varsigma(k+1)}{k} = \sum_{k=1}^{\infty} \frac{1}{k^2} \sum_{j=1}^{k} \frac{(-1)^{j-1}}{j!} Y_{j-1}\left(H_{k-1}^{(1)}, -1! H_{k-1}^{(2)}, \ldots, (-1)^{j-1}(j-2)! H_{k-1}^{(j-1)}\right)$$

which we have seen a different representation of this in (3.28).

## 7. Euler sums

Differentiating Nielsen's series (3.1) via the Leibniz rule results in

$$\sum_{k=0}^{r}\binom{r}{k}D^{k}\left[\psi(x)+\gamma\right]D^{r-k}\left[\psi(x)+\gamma\right]=\psi^{(r+1)}(x)-2(-1)^{r}r!\sum_{n=1}^{\infty}\frac{H_{n}}{(x+n)^{r+1}}$$

and dealing separately with the terms involving $k=0$ and $k=r$ this may be written as

$$\sum_{k=1}^{r-1}\binom{r}{k}\psi^{(k)}(x)\psi^{(r-k)}(x)+2\left[\psi(x)+\gamma\right]\psi^{(r)}(x)=\psi^{(r+1)}(x)-2(-1)^{r}r!\sum_{n=1}^{\infty}\frac{H_{n}}{(x+n)^{r+1}}$$

It is well known that [79, p.22]

(7.1) $\quad \psi^{(k)}(x)=(-1)^{k+1}k!\varsigma(k+1,x)$

where $\varsigma(s,x)$ is the Hurwitz zeta function. We therefore obtain

$$(7.2)\; 2\sum_{n=1}^{\infty}\frac{H_{n}}{(x+n)^{r+1}}=(r+1)\varsigma(r+2,x)-\sum_{k=1}^{r-1}\varsigma(k+1,x)\varsigma(r-k+1,x)+2\left[\psi(x)+\gamma\right]\varsigma(r+1,x)$$

The following relation was also recently derived by Eie [38, p.102] in a different manner for $x \in (0,1)$

$$(7.3)\; 2\sum_{n=0}^{\infty}\frac{H_{n}}{(x+n)^{r}}=r\varsigma(r+1,x)-\sum_{k=2}^{r-1}\varsigma(k,x)\varsigma(r-k+1,x)+2\left[\psi(x)+\gamma\right]\varsigma(r,x)$$

where the summation on the left-hand side starts at $n=0$. Since $H_{0}=0$, it is easy to show that (7.3) and (7.3) are equivalent.

Letting $r+1 \to r$ in (7.2) results in

$$(7.3.1)\quad 2\sum_{n=1}^{\infty}\frac{H_{n}}{(n+1)^{r}}=r\varsigma(r+1)-\sum_{k=1}^{r-2}\varsigma(k+1)\varsigma(r-k)$$

and letting $m=n+1$ we see that $\sum_{n=1}^{\infty}\frac{H_{n}}{(n+1)^{r}}=\sum_{m=2}^{\infty}\frac{H_{m-1}}{m^{r}}=\sum_{m=1}^{\infty}\frac{H_{m-1}}{m^{r}}$ since $H_{0}=0$.

Therefore we obtain

$$(7.4)\quad 2\sum_{n=1}^{\infty}\frac{H_{n-1}}{n^{r}}=r\varsigma(r+1)-\sum_{k=1}^{r-2}\varsigma(k+1)\varsigma(r-k)$$

which was first obtained by Euler in 1775 ( a considerably lengthier derivation appears in [79, p.103]).

With $r=1$ in (7.2) we have

$$2\sum_{n=1}^{\infty}\frac{H_n}{(x+n)^2} = 2\varsigma(3,x)+[\psi(x)+\gamma]\varsigma(2,x)$$

and with $x=1/2$ this becomes

$$4\sum_{n=1}^{\infty}\frac{H_n}{(2n+1)^2} = \varsigma\left(3,\frac{1}{2}\right)-\log 2\cdot\varsigma\left(2,\frac{1}{2}\right)$$

Since $\varsigma\left(s,\frac{1}{2}\right)=[2^s-1]\varsigma(s)$ we obtain

(7.5) $\quad 4\sum_{n=1}^{\infty}\dfrac{H_n}{(2n+1)^2} = 7\varsigma(3)-3\log 2\cdot\varsigma(2)$

which may be compared with the formula found by Bailey et al. [4]

(7.6) $\quad 4\sum_{n=1}^{\infty}\dfrac{H_n}{(2n-1)^2} = 7\varsigma(3)+\pi^2-\pi^2\log 2-8\log 2$

which displays a curious structure in that it lacks homogeneity of the weights of the various factors which is the norm for Euler sums (see [10, p.203]).

We also have another Euler sum due to Chen [14]

$$\sum_{n=1}^{\infty}\frac{H_{2n}}{(2n+1)^2} = \frac{7}{16}\varsigma(3)$$

Some recent results on Euler sums are contained in [37] and [62].

**8. Nielsen's beta function $\beta(x)$**

We now consider the function $\beta(x)$ which is defined as

(8.1) $\quad \beta(x) = \sum_{n=0}^{\infty}\dfrac{(-1)^n}{n+x}$

This is dealt with in some detail by Nielsen [67, pp.16, 51, 193, 231] and also by Bromwich [11, p.523] and we note that

$$\beta(1) = \log 2 \qquad \beta(1/2) = \frac{\pi}{2}$$

$$\beta(1+x) = \frac{1}{x} - \beta(x) \qquad \beta(x) = \psi(x) - \psi(x/2) - \log 2$$

$\beta(x)$ is referred to as the incomplete beta function in [8].

Differentiation results in

$$\beta^{(r)}(x) = (-1)^r r! \sum_{n=0}^{\infty} \frac{(-1)^n}{(n+x)^{r+1}}$$

and specifically we have

(8.2) $$\beta^{(r)}(1) = (-1)^r r! \sum_{n=0}^{\infty} \frac{(-1)^n}{(n+1)^{r+1}} = (-1)^r r! \varsigma_a(r+1)$$

where $\varsigma_a(r)$ is the alternating Riemann zeta function defined by

$$\varsigma_a(r) = \sum_{n=1}^{\infty} \frac{(-1)^{n+1}}{n^r} = (1 - 2^{1-r})\varsigma(r)$$

In view of (8.2) we have the Maclaurin expansion

(8.3) $$\beta(1+x) = \sum_{n=0}^{\infty} (-1)^n \varsigma_a(n+1) x^n$$

Nielsen [67] has shown that

(8.4) $$\beta^2(x) = \psi'(x) - 2h(x)$$

where

(8.5) $$h(x) = \sum_{n=1}^{\infty} \frac{(-1)^{n+1}}{n} H_n^{(1)}(x)$$

and $H_n^{(1)}(x)$ is the harmonic number function defined by

(8.6) $$H_n^{(1)}(x) = \sum_{j=0}^{n-1} \frac{1}{j+x}$$

and we note that $H_n^{(1)}(1) = H_n^{(1)}$.

With $x = 1$ in (8.4) we obtain

(8.7) $$\sum_{n=1}^{\infty} \frac{(-1)^{n+1} H_n^{(1)}}{n} = \frac{1}{2}\left[\varsigma(2) - \log^2 2\right]$$

For convenience of notation, we use the operator $D^r = \dfrac{d^r}{dx^r}$ and reference to (8.6) shows that

(8.8) $$D^r H_n^{(1)}(x) = (-1)^r r! \sum_{j=0}^{n-1} \frac{1}{(j+x)^{r+1}} = (-1)^r r! H_n^{(r+1)}(x)$$

Differentiation of (8.4) gives us

(8.9) $$2\beta(x)\beta'(x) = \psi''(x) + 2\sum_{n=1}^{\infty} \frac{(-1)^{n+1}}{n} H_n^{(2)}(x)$$

and with $x = 1$ we have

(8.10) $$\sum_{n=1}^{\infty} \frac{(-1)^{n+1}}{n} H_n^{(2)} = \varsigma(3) - \varsigma(2)\log 2$$

Applying (3.22) to (8.10) results in

$$\sum_{n=1}^{\infty} \frac{(-1)^{n+1}}{n} H_n^{(2)} = \sum_{n=1}^{\infty} \frac{1}{n^2} \sum_{k=n}^{\infty} \frac{(-1)^{k+1}}{k}$$

$$= \sum_{n=1}^{\infty} \frac{1}{n^2} \left[ \sum_{k=1}^{\infty} \frac{(-1)^{k+1}}{k} - \sum_{k=1}^{n-1} \frac{(-1)^{k+1}}{k} \right]$$

$$= \sum_{n=1}^{\infty} \frac{1}{n^2} \sum_{k=1}^{\infty} \frac{(-1)^{k+1}}{k} - \sum_{n=1}^{\infty} \frac{1}{n^2} \sum_{k=1}^{n-1} \frac{(-1)^{k+1}}{k}$$

$$= \varsigma(2)\log 2 - \sum_{n=1}^{\infty} \frac{1}{n^2} \sum_{k=1}^{n-1} \frac{(-1)^{k+1}}{k}$$

We see that

$$\sum_{n=1}^{\infty} \frac{1}{n^2} \sum_{k=1}^{n-1} \frac{(-1)^{k+1}}{k} = \sum_{n=1}^{\infty} \frac{1}{n^2} \sum_{k=1}^{n} \frac{(-1)^{k+1}}{k} - \sum_{n=1}^{\infty} \frac{(-1)^{n+1}}{n^3}$$

and thus we have

$$(8.11) \qquad \sum_{n=1}^{\infty} \frac{1}{n^2} \sum_{k=1}^{n} \frac{(-1)^{k+1}}{k} = 2\varsigma(2)\log 2 - \varsigma(3) + \varsigma_a(3)$$

Higher derivatives of (8.9) will result in $\sum_{k=1}^{\infty} \frac{(-1)^{n+1}}{n} H_n^{(r)}$ which, as shown above, may be transformed into $\sum_{n=1}^{\infty} \frac{1}{n^r} \sum_{k=1}^{n} \frac{(-1)^{k+1}}{k}$. Whilst such plain vanilla Euler sums are well documented in the literature (see for example [39], [61], [71]), the approach used in this section does however have the added advantage of presenting new formulae of the form

$$\sum_{j=0}^{p} \binom{p}{j} \beta^{(j)}(x) \beta^{(p-j)}(x) = \psi^{(p+1)}(x) - 2(-1)^p p! \sum_{n=1}^{\infty} \frac{(-1)^{n+1}}{n} H_n^{(p+1)}(x)$$

obtained by differentiating (8.4).

Using the following variant of (3.22)

$$(8.12) \qquad \sum_{n=1}^{\infty} a_n \sum_{k=0}^{n-1} b_k = \sum_{n=0}^{\infty} b_n \sum_{k=n+1}^{\infty} a_k$$

in the same manner as above we obtain

$$\sum_{k=1}^{\infty} \frac{(-1)^{n+1}}{n} H_n^{(p+1)}(x) = \sum_{n=0}^{\infty} \frac{1}{(n+x)^{p+1}} \sum_{k=n+1}^{\infty} \frac{(-1)^{k+1}}{k}$$

$$= \sum_{n=0}^{\infty} \frac{1}{(n+x)^{p+1}} \sum_{k=1}^{\infty} \frac{(-1)^{k+1}}{k} - \sum_{n=0}^{\infty} \frac{1}{(n+x)^{p+1}} \sum_{k=1}^{n} \frac{(-1)^{k+1}}{k}$$

and thus

$$(8.13) \qquad \sum_{k=1}^{\infty} \frac{(-1)^{n+1}}{n} H_n^{(p+1)}(x) = \varsigma(p+1, x)\log 2 - \sum_{n=0}^{\infty} \frac{1}{(n+x)^{p+1}} \sum_{k=1}^{n} \frac{(-1)^{k+1}}{k}$$

## 9. A tenuous connection with Dilcher's generalised gamma functions

We showed in Eq.(10.17) in [29] that

$$(9.1) \qquad \gamma_1(1+x) - \gamma_1 = -\sum_{n=1}^{\infty} (-1)^n \left[ \varsigma'(n+1) + H_n^{(1)} \varsigma(n+1) \right] x^n$$

so that

$$\gamma_1(1-x) - \gamma_1 = -\sum_{n=1}^{\infty}\left[\varsigma'(n+1) + H_n^{(1)}\varsigma(n+1)\right]x^n \quad (9.2)$$

This suggests that

$$\lim_{x \to 1}\left[\gamma_1(1-x) - \gamma_1 + \sum_{n=1}^{\infty} H_n^{(1)}\varsigma(n+1)x^n\right] = -\sum_{n=2}^{\infty}\varsigma'(n) \quad (9.3)$$

where the right-hand side corresponds with (3.10.3). It therefore appears that

$$\sum_{n=1}^{\infty} \frac{\log(n+1)}{n(n+1)} = \lim_{x \to 1}\left[\gamma_1(1-x) - \gamma_1 + \sum_{n=1}^{\infty} H_n^{(1)}\varsigma(n+1)x^n\right] \quad (9.4)$$

Since

$$\gamma_p(1+x) - \gamma_p(x) = -\frac{\log^p x}{x} \quad (9.5)$$

we have

$$\gamma_1(x) - \gamma_1 = \frac{\log x}{x} - \sum_{n=1}^{\infty}(-1)^n \varsigma'(n+1)x^n - \sum_{n=1}^{\infty}(-1)^n H_n \varsigma(n+1)x^n \quad (9.6)$$

We now recall

$$\gamma_1(x) - \gamma_1(1) = \sum_{n=0}^{\infty}\left[\frac{\log(n+x)}{n+x} - \frac{\log(n+1)}{n+1}\right] \quad (9.7)$$

and noting that Ramanujan [6, p.197] defined $\varphi(x)$ by

$$\varphi(x) = \sum_{n=1}^{\infty}\left[\frac{\log n}{n} - \frac{\log(n+x)}{n+x}\right] \quad (9.8)$$

we see that

$$\gamma_1(x) - \gamma_1(1) = -\sum_{n=0}^{\infty}\left[\frac{\log(n+1)}{n+1} - \frac{\log(n+x)}{n+x}\right]$$

$$= -\sum_{n=1}^{\infty}\left[\frac{\log(n+1)}{n+1} - \frac{\log(n+x)}{n+x}\right] + \frac{\log x}{x}$$

$$= \frac{\log x}{x} - \sum_{n=1}^{\infty} \left[ \frac{\log n}{n} - \frac{\log(n+x)}{n+x} + \frac{\log(n+1)}{n+1} - \frac{\log n}{n} \right]$$

We then have

(9.9) $$\gamma_1(x) - \gamma_1(1) = \frac{\log x}{x} - \sum_{n=1}^{\infty} \left[ \frac{\log n}{n} - \frac{\log(n+x)}{n+x} \right]$$

and thus

(9.10) $$\gamma_1(x) - \gamma_1 = \frac{\log x}{x} - \varphi(x)$$

We note that $\varphi(1) = 0$. It is not immediately obvious how this representation may be used to define $\varphi(x)$ for $x \leq 0$ but referring to (9.6)

$$\gamma_1(x) - \gamma_1 = \frac{\log x}{x} - \sum_{n=1}^{\infty} (-1)^n \varsigma'(n+1) x^n - \sum_{n=1}^{\infty} (-1)^n H_n \varsigma(n+1) x^n$$

it seems natural to write

(9.11) $$\varphi(x) = \sum_{n=1}^{\infty} (-1)^n \varsigma'(n+1) x^n + \sum_{n=1}^{\infty} (-1)^n H_n \varsigma(n+1) x^n$$

whereupon we have for negative values

(9.12) $$\varphi(-x) = \sum_{n=1}^{\infty} \varsigma'(n+1) x^n + \sum_{n=1}^{\infty} H_n \varsigma(n+1) x^n$$

Referring to (9.1) we therefore have

(9.13) $$\gamma_1(1-x) - \gamma_1 = -\varphi(-x)$$

With $x = 1/2$ we have

(9.14) $$\gamma_1\left(\frac{1}{2}\right) - \gamma_1 = -\left[ \sum_{n=1}^{\infty} \frac{\varsigma'(n+1)}{2^2} + \sum_{n=1}^{\infty} \frac{H_n \varsigma(n+1)}{2^2} \right]$$

Using

$$\gamma_1\left(\frac{1}{2}\right) = \gamma_1 - \log^2 2 - 2\gamma \log 2$$

we see that

$$(9.15) \qquad \log^2 2 + 2\gamma \log 2 = \left[ \sum_{n=1}^{\infty} \frac{\varsigma'(n+1)}{2^2} + \sum_{n=1}^{\infty} \frac{H_n \varsigma(n+1)}{2^2} \right]$$

which agrees with Ramanujan's result [6, p.199]

$$(9.16) \qquad \varphi\left(-\frac{1}{2}\right) = \log^2 2 + 2\gamma \log 2$$

□

Letting $-x \to x-1$ in (9.13) we obtain

$$(9.17) \qquad \gamma_1(x) - \gamma_1 = -\varphi(x-1)$$

and we then deduce that (see [6, p.200])

$$(9.18) \qquad \varphi(x-1) - \varphi(-x) = \gamma_1(1-x) - \gamma_1(x)$$

$$(9.19) \qquad \varphi(x-1) + \varphi(-x) = -[\gamma_1(1-x) + \gamma_1(x)] + 2\gamma_1$$

For example, we have with $x = 1/4$

$$\varphi(-3/4) + \varphi(-1/4) = -[\gamma_1(3/4) + \gamma_1(1/4)] + 2\gamma_1$$

and referring to (5.17) we see that

$$\varphi(-3/4) + \varphi(-1/4) = 7\log^2 2 + 6\gamma \log 2$$

which concurs with Ramanujan's result [6, p.199].

We also have

$$\varphi(-1/4) - \varphi(-3/4) = \gamma_1(1/4) - \gamma_1(3/4)$$

and referring to (5.17) we obtain

$$(9.20) \qquad \varphi(-1/4) - \varphi(-3/4) = -\pi \left[ \gamma + 4\log 2 + 3\log \pi - 4\log \Gamma\left(\frac{1}{4}\right) \right]$$

Using the definition (11.1) Ramanujan [6, p.198] showed that

$$\frac{1}{4}[\varphi(-1/4) - \varphi(-3/4)] = \sum_{n=1}^{\infty}\left[\frac{\log(n-3/4)}{4n-3} - \frac{\log(n-1/4)}{4n-1}\right]$$

$$= \sum_{n=0}^{\infty}\frac{(-1)^n \log(2n+1)}{2n+1} - \log 4 \sum_{n=0}^{\infty}\frac{(-1)^n}{2n+1}$$

$$= \sum_{n=0}^{\infty}\frac{(-1)^n \log(2n+1)}{2n+1} - \frac{1}{2}\pi \log 2$$

and using the known result

$$\sum_{n=0}^{\infty}\frac{(-1)^n \log(2n+1)}{2n+1} = -\left[\gamma + 2\log 2 + 3\log \pi - 4\log \Gamma\left(\frac{1}{4}\right)\right]$$

this confirms (9.20) above.

□

Dilcher [36] has defined generalised gamma functions for $x > 0$ by

(9.21) $$\Gamma_k(x) = \lim_{n \to \infty} \frac{\exp\left[\frac{x}{k+1}\log^{k+1} n\right]\prod_{j=1}^{n}\exp\left[\frac{1}{k+1}\log^{k+1} j\right]}{\prod_{j=0}^{n}\exp\left[\frac{1}{k+1}\log^{k+1}(j+x)\right]}$$

where we note that $\Gamma_0(x)$ corresponds with the classical gamma function

(9.22) $$\Gamma_0(x) = \lim_{n \to \infty}\frac{n^x n!}{x(x+1)\cdots(x+n)} = \Gamma(x)$$

and we have

$$\Gamma_k(1) = 1$$

$$\Gamma_k(x+1) = \exp\left[\frac{1}{k+1}\log^{k+1} x\right]\Gamma_k(x)$$

Is it possible to extend $\Gamma_k(x)$ by analytic continuation?

We showed in [29] that

(9.23) $$\sum_{n=1}^{\infty} \frac{(-1)^n}{n+1}[H_n \varsigma(n+1) + \varsigma'(n+1)]x^{n+1} = \sum_{n=1}^{\infty}\left[x\frac{\log n}{n} + \frac{1}{2}\log^2 n - \frac{1}{2}\log^2(n+x)\right]$$

and differentiation results in

(9.24) $$\sum_{n=1}^{\infty}(-1)^n [H_n \varsigma(n+1) + \varsigma'(n+1)]x^n = \sum_{n=1}^{\infty}\left[\frac{\log n}{n} - \frac{\log(n+x)}{n+x}\right]$$

which corresponds with (9.1).

Dilcher [36] showed that

(9.25) $$\log \Gamma_1(x+1) + \gamma_1 x = \sum_{n=1}^{\infty} \frac{(-1)^n}{n+1}[H_n \varsigma(n+1) + \varsigma'(n+1)]x^{n+1}$$

By definition we have

$$\log \Gamma_1(x) = \lim_{n\to\infty}\left[\frac{x}{2}\log^2 n + \frac{1}{2}\sum_{j=1}^{n}\log^2 j - \frac{1}{2}\sum_{j=0}^{n}\log^2(j+x)\right]$$

and we have

$$\log \Gamma_1(x+1) = \frac{1}{2}\log^2 x + \log \Gamma_1(x)$$

Therefore we have

$$\log \Gamma_1(x+1) + \gamma_1 x = \lim_{n\to\infty}\left[\frac{x}{2}\log^2 n + \frac{1}{2}\sum_{j=1}^{n}\log^2 j - \frac{1}{2}\sum_{j=1}^{n}\log^2(j+x)\right] + \lim_{n\to\infty}\sum_{j=1}^{n}\left[x\frac{\log j}{j} - \frac{x}{2}\log^2 n\right]$$

$$= \lim_{n\to\infty}\left[x\frac{\log j}{j} + \frac{1}{2}\sum_{j=1}^{n}\log^2 j - \frac{1}{2}\sum_{j=1}^{n}\log^2(j+x)\right]$$

and hence we see that (9.23) follows from (9.24).

It was also shown by Dilcher [36] that

$$\frac{1}{\Gamma_k(x)} = e^{\gamma_k x} \exp\left[\frac{x}{k+1}\log^{k+1} x\right] \prod_{n=1}^{\infty} \exp\left[-\frac{x}{n}\log^k n\right]\exp\left[\frac{1}{k+1}[\log^{k+1}(n+x) - \log^{k+1} n]\right]$$

whereupon it directly follows that

$$\log \Gamma_k(x+1) + \gamma_k x = \lim_{n \to \infty} \left[ x \frac{\log^k j}{j} + \frac{1}{k+1} \sum_{j=1}^{n} \log^k j - \frac{1}{k+1} \sum_{j=1}^{n} \log^k (j+x) \right]$$

□

We have the known integral expression for the harmonic numbers $H_n$ (see for example [22])

$$H_n = -n \int_0^1 (1-u)^{n-1} \log u \, du$$

and

$$\varsigma(n+1) = \frac{(-1)^n}{n!} \int_0^1 \frac{\log^n(1-v)}{v} dv$$

This gives us

$$\sum_{n=1}^{\infty} H_n^{(1)} \varsigma(n+1) x^n = -\sum_{n=1}^{\infty} \frac{(-1)^n n x^n}{n!} \int_0^1 \frac{\log^n(1-v)}{v} dv \int_0^1 (1-u)^{n-1} \log u \, du$$

$$= x \int_0^1 \frac{\log u}{v \log(1-u)} \int_0^1 \sum_{n=1}^{\infty} \frac{n[-x(1-u)\log(1-v)]^{n-1}}{n!} du \, dv$$

We see that

$$\sum_{n=0}^{\infty} \frac{[-x(1-u)\log(1-v)]^n}{n!} = \exp[-x(1-u)\log(1-v)]$$

and differentiating with respect to $x$ gives us

$$\sum_{n=0}^{\infty} \frac{n[-x(1-u)\log(1-v)]^{n-1}}{n!} = -(1-u)\log(1-v) \exp[-x(1-u)\log(1-v)]$$

$$= -\frac{(1-u)\log(1-v)}{(1-v)^{x(1-u)}}$$

Hence we have

$$\sum_{n=1}^{\infty} H_n^{(1)} \varsigma(n+1) x^n = -x \int_0^1 \frac{\log u}{v \log(1-u)} \int_0^1 \frac{(1-u)\log(1-v)}{(1-v)^{x(1-u)}} du \, dv$$

$$= -x \int_0^1 \int_0^1 \frac{(1-u)\log u \log(1-v)}{(1-v)^{x(1-u)} v \log(1-u)} \, du \, dv$$

However, it is not immediately obvious to me why the right-hand side diverges at $x = 1$.

## 10. An integral representation of Nielsen's function $\xi(x)$

This section was inspired by the approach recently adopted by Morales [65] (see also [27]).

The beta function $B(u, v)$ is defined for $\operatorname{Re}(u) > 0$ and $\operatorname{Re}(v) > 0$ by the Eulerian integral

$$B(u, v) = \int_0^1 t^{u-1} (1-t)^{v-1} \, dt$$

and it is well known that

$$B(u, v) = \frac{\Gamma(u) \Gamma(v)}{\Gamma(u+v)}$$

where $\Gamma(u)$ is the gamma function.

We see that

$$B(u, v+1) = \frac{\Gamma(u) \Gamma(v+1)}{\Gamma(u+v+1)} = \frac{v \Gamma(u) \Gamma(v)}{(u+v) \Gamma(u+v)} = \frac{v}{(u+v)} B(u, v)$$

and we immediately have

$$B(u, 1) = \frac{\Gamma(u) \Gamma(1)}{\Gamma(u+1)} = \frac{1}{u}$$

We have

$$\lim_{v \to 0} \frac{v}{(u+v)} B(u, v) = \lim_{v \to 0} \frac{1}{(u+v)} \lim_{v \to 0} [v B(u, v)] = \frac{1}{u} \lim_{v \to 0} [v B(u, v)]$$

but we also have $\lim_{v \to 0} B(u, v+1) = \frac{1}{u}$ and this implies that

$$\lim_{v \to 0} [v B(u, v)] = 1$$

We now consider

$$(10.1) \quad B(u,v) - \frac{1}{v} = \frac{vB(u,v)-1}{v} = \frac{(u+v)B(u,v+1)-1}{v}$$

and take the limit of the right-hand side

$$\lim_{v \to 0} \frac{(u+v)B(u,v+1)-1}{v} = \lim_{v \to 0}\left[(u+v)\frac{\partial}{\partial v}B(u,v+1) + B(u,v+1)\right]$$

$$= u\frac{\partial}{\partial v}B(u,v)\bigg|_{v=1} + \frac{1}{u}$$

where we have applied L'Hôpital's rule.

We accordingly obtain an important identity (refer to [27] for further applications)

$$\lim_{v \to 0}\left[B(u,v) - \frac{1}{v}\right] = u\frac{\partial}{\partial v}B(u,v)\bigg|_{v=1} + \frac{1}{u}$$

We now consider the integral

$$(10.2) \quad h(u,v) = -\frac{\partial}{\partial v}\int_0^1 (t^{u-1}-1)(1-t)^{v-1}\,dt$$

$$= -\int_0^1 (t^{u-1}-1)(1-t)^{v-1}\log(1-t)\,dt$$

and it is easily seen that

$$(10.3) \quad h(u,0) = \int_0^1 \frac{t^{u-1}-1}{t-1}\log(1-t)\,dt$$

We see that

$$h(u,v) = -\frac{\partial}{\partial v}\left[B(u,v) - \frac{1}{v}\right]$$

$$= -\frac{\partial}{\partial v}\frac{(u+v)B(u,v+1)-1}{v}$$

$$= -\frac{v[(u+v)\frac{\partial}{\partial v}B(u,v+1) + B(u,v+1)] - (u+v)B(u,v+1)+1}{v^2}$$

Differentiation gives us

(10.4) $$\frac{\partial}{\partial v} B(u,v) = \frac{\Gamma(u)\Gamma(v)}{\Gamma(u+v)}[\psi(v) - \psi(u+v)]$$

which we express as

$$B'(u,v) = B(u,v)[\psi(v) - \psi(u+v)]$$

We have

$$h(u,0) = -\lim_{v \to 0} \frac{v[(u+v)B'(u,v+1) + B(u,v+1)] - (u+v)B(u,v+1) + 1}{v^2}$$

Defining $g(v) = v[(u+v)B'(u,v+1) + B(u,v+1)] - (u+v)B(u,v+1) + 1$ we note that $g(0) = 0$ and hence we may apply L'Hôpital's rule. This gives us

$$h(u,0) = -\lim_{v \to 0} \frac{g'(v)}{2v}$$

We have

$$g'(v) = v[(u+v)B''(u,v+1) + 2B'(u,v+1)]$$

and since $g'(0) = 0$ we may apply L'Hôpital's rule again. This gives us

$$h(u,0) = -\lim_{v \to 0} \frac{g''(v)}{2}$$

We have

$$g''(v) = v[(u+v)B'''(u,v+1) + 3B''(u,v+1)] + [(u+v)B''(u,v+1) + 2B'(u,v+1)]$$

so that $g''(0) = uB''(u,1) + 2B'(u,1)$ and thus

$$h(u,0) = -\frac{1}{2}[uB''(u,1) + 2B'(u,1)]$$

We have

(10.5) $$\frac{\partial^2}{\partial v^2} B(u,v) = \frac{\Gamma(u)\Gamma(v)}{\Gamma(u+v)}\left([\psi(v) - \psi(u+v)]^2 + [\psi'(v) - \psi'(u+v)]\right)$$

and, in fact, because of (10.5), $\frac{\partial^n}{\partial v^n} B(u,v)$ may be expressed in terms involving the (exponential) complete Bell polynomials [30].

Since $u \dfrac{\partial^2}{\partial v^2} B(u,v)\Big|_{v=1} = [\psi(1)-\psi(1+u)]^2 + [\psi'(1)-\psi'(1+u)]$ we obtain

$$h(u,0) = -\frac{1}{2}\left\{[\gamma+\psi(1+u)]^2 + [\psi'(1)-\psi'(1+u)] + 2\frac{\psi(1)-\psi(1+u)}{u}\right\}$$

which may be written as

$$= \frac{1}{2}\left\{\psi'(u) - \varsigma(2) - [\psi(u)+\gamma]^2\right\}$$

Referring to (3.1)

$$2\xi(u) = \psi'(u) - \varsigma(2) - [\psi(u)+\gamma]^2$$

we then obtain

(10.6) $\qquad \displaystyle\int_0^1 \frac{t^{u-1}-1}{t-1}\log(1-t)\,dt = \xi(u)$

which was previously reported by Coffey [19].

It appears that Coffey's formula was obtained in the following manner:

Using the representation

$$\frac{1}{a} = \int_0^1 t^{a-1}\,dt$$

in (3.2) we have

$$\xi(u) = \sum_{n=1}^\infty H_n \left(\int_0^1 t^{u+n-1}\,dt - \int_0^1 t^n\,dt\right)$$

$$= \int_0^1 (t^{u-1}-1)\sum_{n=1}^\infty H_n t^n\,dt$$

where the interchange of summation and integration is justified by absolute convergence.

Then utilising the generating series for the harmonic numbers

$$-\frac{\log(1-t)}{1-t} = \sum_{n=1}^{\infty} H_n t^n$$

we deduce that

$$\xi(u) = \int_0^1 \frac{t^{u-1}-1}{t-1} \log(1-t)\, dt$$

Letting $u = 2$ in (10.6) results in

$$\int_0^1 \log(1-t)\, dt = \frac{1}{2}\left\{\varsigma(2,2) - \varsigma(2) - [\psi(2)+\gamma]^2\right\}$$

and, since $\varsigma(2,2) = \varsigma(2) - 1$, we obtain the expected result

$$\int_0^1 \log(1-t)\, dt = -1$$

We recall (3.2)

$$\xi(u) = \sum_{n=1}^{\infty} H_n \left(\frac{1}{u+n} - \frac{1}{n+1}\right)$$

and letting $u = 2$ gives us

$$\xi(2) = \sum_{n=1}^{\infty}\left(\frac{H_n}{n+2} - \frac{H_n}{n+1}\right)$$

$$= \sum_{n=1}^{\infty}\left(\frac{H_{n+2}}{n+2} - \frac{1}{(n+1)(n+2)} - \frac{1}{(n+2)^2} - \frac{H_n}{n+1}\right)$$

$$= \sum_{n=1}^{\infty}\left(\frac{H_{n+2}}{n+2} - \frac{1}{(n+2)^2} - \frac{H_n}{n+1}\right) - \sum_{n=1}^{\infty}\frac{1}{(n+1)(n+2)}$$

$$= -\frac{H_2}{2} + \sum_{n=1}^{\infty}\left(\frac{H_{n+1}}{n+1} - \frac{1}{(n+2)^2} - \frac{H_n}{n+1}\right) - \sum_{n=1}^{\infty}\frac{1}{(n+1)(n+2)}$$

$$= -\frac{H_2}{2} + \sum_{n=1}^{\infty}\left(\frac{H_n}{n+1} + \frac{1}{(n+1)^2} - \frac{1}{(n+2)^2} - \frac{H_n}{n+1}\right) - \sum_{n=1}^{\infty}\frac{1}{(n+1)(n+2)}$$

$$= -\frac{H_2}{2} + \frac{1}{2^2} - \sum_{n=1}^{\infty} \frac{1}{(n+1)(n+2)}$$

Employing the telescoping series

$$\sum_{n=1}^{\infty} \frac{1}{(n+1)(n+2)} = \sum_{n=1}^{\infty} \left( \frac{1}{n+1} - \frac{1}{n+2} \right) = \frac{1}{2}$$

we see that $\xi(2) = -1$

$\square$

Using integration by parts and the identities

$$\frac{d}{dt} Li_2(1-t) = -\frac{Li_1(1-t)}{1-t} = \frac{\log t}{1-t} \quad \text{and} \quad \frac{d}{dt} Li_3(1-t) = -\frac{Li_2(1-t)}{1-t}$$

which may be obtained from the series definition of the relevant polylogarithm, we get

$$\int \frac{\log(1-t) \log t}{1-t} dt = \log(1-t) Li_2(1-t) + \int \frac{Li_2(1-t)}{1-t} dt$$

$$= \log(1-t) Li_2(1-t) - Li_3(1-t)$$

We therefore have the known integral

$$\int_0^1 \frac{\log(1-t) \log t}{1-t} dt = \varsigma(3)$$

Differentiating (10.6) gives us

$$\xi'(u) = \int_0^1 \frac{t^{u-1} \log t \log(1-t)}{t-1} dt$$

so that

$$\xi'(1) = \int_0^1 \frac{\log t \log(1-t)}{t-1} dt = -\varsigma(3)$$

We also have

$$\xi'(u) = -\sum_{n=1}^{\infty} \frac{H_n}{(u+n)^2}$$

so that

$$\xi'(1) = -\sum_{n=1}^{\infty} \frac{H_n}{(n+1)^2}$$

$$= -\left[\sum_{n=1}^{\infty} \frac{H_{n+1}}{(n+1)^2} - \frac{1}{(n+1)^3}\right]$$

$$= -\left[\sum_{n=1}^{\infty} \frac{H_{n+1}}{(n+1)^2} - \frac{1}{(n+1)^3}\right]$$

$$= -\left[\sum_{m=2}^{\infty} \frac{H_m}{m^2} - \frac{1}{m^3}\right]$$

$$= -\left[\sum_{m=1}^{\infty} \frac{H_m}{m^2} - \frac{1}{m^3}\right]$$

$$= -\sum_{m=1}^{\infty} \frac{H_m}{m^2} + \varsigma(3)$$

This gives us the well-known Euler sum

$$\sum_{n=1}^{\infty} \frac{H_n}{n^2} = 2\varsigma(3)$$

Differentiating (10.6) gives us

$$\xi^{(p)}(u) = \int_0^1 \frac{t^{u-1} \log^p t \log(1-t)}{t-1} dt$$

We also have

$$\xi^{(p)}(u) = (-1)^p p! \sum_{n=1}^{\infty} \frac{H_n}{(u+n)^{p+1}}$$

so that

$$\int_0^1 \frac{t^{u-1} \log^p t \log(1-t)}{t-1} dt = (-1)^p p! \sum_{n=1}^{\infty} \frac{H_n}{(u+n)^{p+1}}$$

and hence we see that

$$2(-1)^p p! \sum_{n=1}^{\infty} \frac{H_n}{(u+n)^{p+1}} = \frac{d^p}{du^p}\left(\psi'(u) - \varsigma(2) - [\psi(u)+\gamma]^2\right)$$

□

A similar procedure enables us to evaluate the family of integrals

$$\int_0^1 \frac{t^{u-1} \log^p t \log^q (1-t)}{t-1} dt \qquad \int_0^1 \frac{(t^{u-1}-1) \log^q (1-t)}{t-1} dt$$

and a very direct derivation follows from the next proposition.

**Proposition**

More generally, we consider the function $F(u,v)$ defined by

(10.10) $$F(u,v) = \frac{G(u,v) - 1}{v}$$

where the partial derivatives $G^{(n)}(u,v) := \frac{\partial^n}{\partial v^n} G(u,v)$ for all $n \geq 1$ exist for all $v \in [0,a]$. We also require that $G(u,0) = 1$.

Then we have

(10.11) $$F^{(n)}(u,0) = \frac{G^{(n+1)}(u,0)}{n+1}$$

where $F^{(n)}(u,v) := \frac{\partial^n}{\partial v^n} F(u,v)$ and $F^{(n)}(u,0) := F^{(n)}(u,v)\big|_{v=0}$.

**Proof**

We have
$$\lim_{v \to 0} F(u,v) = \lim_{v \to 0} \frac{G(u,v) - 1}{v}$$

and since $G(u,0) = 1$ we may apply L'Hôpital's rule to obtain

$$\lim_{v \to 0} F(u,v) = \lim_{v \to 0} G^{(1)}(u,v)$$

so that

$$F(u,0) = G^{(1)}(u,0)$$

We write (10.10) in the form

$$vF(u,v) = G(u,v) - 1$$

and, applying the Leibniz rule for differentiation, we obtain for all $n \geq 1$

(10.12) $$vF^{(n+1)}(u,v)+(n+1)F^{(n)}(u,v)=G^{(n+1)}(u,v)$$

With $n=0$ we have

$$F^{(1)}(u,v)=\frac{G^{(1)}(u,v)-F(u,v)}{v}$$

and, since $G^{(1)}(u,0)-F(u,0)=0$, we may apply L'Hôpital's rule to obtain

$$F^{(1)}(u,0)=G^{(2)}(u,0)-F^{(1)}(u,0)$$

so that

$$F^{(1)}(u,0)=\frac{G^{(2)}(u,0)}{2}$$

Hence (10.11) is true for $n=1$.

We see from (10.12) that

$$F^{(n+1)}(u,v)=\frac{G^{(n+1)}(u,v)-(n+1)F^{(n)}(u,v)}{v}$$

and if (10.11) is valid for the number $n$ we have

$$G^{(n+1)}(u,0)-(n+1)F^{(n)}(u,0)=0$$

Then we have using L'Hôpital's rule

$$F^{(n+1)}(u,0)=G^{(n+2)}(u,0)-(n+1)F^{(n+1)}(u,0)$$

or equivalently

$$F^{(n+1)}(u,0)=\frac{G^{(n+2)}(u,0)}{n+2}$$

Therefore, assuming that (10.11) is valid for $n$, we deduce that it is also valid for $n+1$. Since we have shown that it is true for $n=1$ the proof by mathematical induction is complete.

□

We now consider an application of this proposition.

Defining $G(u,v)$ for $u > 0$ by

$$G(u,v) = \frac{\Gamma(u)\Gamma(1+v)}{\Gamma(u+v)}$$

we see that $G(u,0) = 1$ and that the other conditions of the above proposition are satisfied.

Differentiation results in

(10.13) $\qquad G^{(1)}(u,v) = G(u,v)[\psi(1+v) - \psi(u+v)]$

and, as mentioned above, because of (10.13), $\dfrac{\partial^n}{\partial v^n} G(u,v)$ may be expressed as follows in terms involving the (exponential) complete Bell polynomials [30]

(10.14) $\qquad G^{(n)}(u,v) = G(u,v) Y_n\left(g(u,v), g^{(1)}(u,v),..., g^{(n-1)}(u,v)\right)$

where
$$g(u,v) = \psi(1+v) - \psi(u+v)$$
and
$$g^{(j)}(u,v) = \frac{\partial^j}{\partial v^j}[\psi(1+v) - \psi(u+v)]$$

This gives us

(10.15) $\qquad G^{(n)}(u,0) = Y_n\left(-\gamma - \psi(u), \psi^{(1)}(u),..., \psi^{(n-1)}(u)\right)$

We recall that
$$h(u,v) = -\frac{\partial}{\partial v}\left[B(u,v) - \frac{1}{v}\right]$$
and since
$$B(u,v) - \frac{1}{v} = \left[\frac{\Gamma(u)\Gamma(1+v)}{\Gamma(u+v)} - 1\right]/v$$
we have
$$h(u,v) = -\frac{\partial}{\partial v} F(u,v)$$

Therefore we have for $n \geq 2$

$$h^{(n-2)}(u,v) = -F^{(n-1)}(u,v)$$

so that

$$h^{(n-2)}(u,0) = -F^{(n-1)}(u,0)$$

Hence using (10.11) gives us

$$h^{(n-2)}(u,0) = -\frac{G^{(n)}(u,0)}{n}$$

Accordingly we obtain

$$h^{(n-2)}(u,0) = \int_0^1 \frac{t^{u-1}-1}{t-1} \log^{n-1}(1-t)\, dt$$

so that

(10.16) $\displaystyle\int_0^1 \frac{t^{u-1}-1}{t-1} \log^{n-1}(1-t)\, dt = -\frac{1}{n} Y_n\left(-\gamma - \psi(u), \psi^{(1)}(u), \ldots, \psi^{(n-1)}(u)\right)$

We see that

(10.17) $\displaystyle\frac{d^p}{du^p} h^{(n-2)}(u,0) = \int_0^1 \frac{t^{u-1} \log^p t \log^{n-1}(1-t)}{t-1}\, dt$

This generalises the approach adopted by Kölbig [53] in 1982 where he considered integrals of the form $\int_0^1 t^{-1} \log^{n-1} t \log^p(1-t)\, dt$. See also the recent papers by Laurenzi [57] and Sofo and Cvijović [75].

□

Integrating (10.6) gives us

(10.18) $\displaystyle\int_0^1 \xi(u)\, du = \int_0^1 \left[\frac{1}{t \log t} - \frac{1}{t-1}\right] \log(1-t)\, dt$

which concurs with (2.16.2) and (3.5).

## 11. Open access to our own work

This paper contains references to 59 other papers and, rather surprisingly, more than 84% of them are currently freely available on the internet. Surely now is the time that <u>all</u> of <u>our</u> work should be freely accessible by <u>all</u>. The mathematics community should lead the way on this by publishing <u>everything</u> on arXiv, or in an equivalent open access repository. We think it, we write it, so why hide it? You know it makes sense.

Wessex House,
Devizes Road,
Upavon,
Pewsey,
Wiltshire SN9 6DL